\documentstyle[amscd,amssymb,xypic,verbatim,11pt]{amsart}
\xyoption{all}
\newcommand{\nc}{\newcommand}
\newcommand{\rnc}{\renewcommand}

\nc{\exto}[1]{\stackrel{#1}{\longrightarrow}}
\nc{\dlim}{{\mathop{\lim\limits_{\longrightarrow}}}}
\nc{\lan}{\big\langle}
\nc{\ran}{\big\rangle}

\nc{\kk}{{\mathsf{k}}}
\nc{\ix}{{\mathsf{i}}}
\nc{\jx}{{\mathsf{j}}}

\nc{\C}{{\mathbb{C}}}
\nc{\HH}{{\mathbb{H}}}
\nc{\PP}{{\mathbb{P}}}
\nc{\QQ}{{\mathbb{Q}}}
\nc{\ZZ}{{\mathbb{Z}}}

\nc{\CA}{{\mathcal{A}}}
\nc{\CB}{{\mathcal{B}}}
\nc{\CC}{{\mathcal{C}}}
\nc{\D}{{\mathcal{D}}}
\nc{\CE}{{\mathcal{E}}}
\nc{\CF}{{\mathcal{F}}}
\nc{\CG}{{\mathcal{G}}}
\nc{\CH}{{\mathcal{H}}}
\nc{\CJ}{{\mathcal{J}}}
\nc{\CK}{{\mathcal{K}}}
\nc{\CL}{{\mathcal{L}}}
\nc{\CM}{{\mathcal{M}}}
\nc{\CN}{{\mathcal{N}}}
\nc{\CO}{{\mathcal{O}}}
\nc{\CQ}{{\mathcal{Q}}}
\nc{\CR}{{\mathcal{R}}}
\nc{\CS}{{\mathcal{S}}}
\nc{\CT}{{\mathcal{T}}}
\nc{\CU}{{\mathcal{U}}}
\nc{\CV}{{\mathcal{V}}}
\nc{\CW}{{\mathcal{W}}}
\nc{\CX}{{\mathcal{X}}}
\nc{\CY}{{\mathcal{Y}}}
\nc{\CZ}{{\mathcal{Z}}}
\nc{\CMo}{{\mathcal{M}^\circ}}
\nc{\Co}{{{C}^\circ}}

\nc{\BY}{{\overline{Y}}}
\nc{\BYD}{{\overline{Y}{}^{|D|}}}
\nc{\OZ}{{\overline{Z}}}
\nc{\bg}{{\bar{g}}}

\nc{\bq}{{\mathbf{q}}}
\nc{\BD}{{\mathbf{D}}}
\nc{\BG}{{\mathbf{G}}}
\nc{\BM}{{\mathbf{M}}}
\nc{\BP}{{\mathbf{P}}}
\nc{\BZ}{{\mathbf{Z}}}
\nc{\BPr}{{\mathsf{P}}}
\nc{\BR}{{\mathbf{R}}}
\nc{\BRO}[1]{{{\mathbf{R}}^{\circ}_{#1}}}
\nc{\BRD}[1]{{{\mathbf{R}}^{|D|}_{#1}}}
\nc{\BRP}[1]{{{\mathbf{R}}^{1}_{#1}}}
\nc{\BRTP}[1]{{{\mathbf{\tilde{R}}}{}^{1}_{#1}}}
\nc{\BS}{{\mathbf{S}}}
\nc{\BMS}{{{\mathbf{M}}^{{s}}}}
\nc{\BMSS}{{{\mathbf{M}}^{{ss}}}}
\nc{\BMZ}{{\mathbf{M}^{\circ}}}
\nc{\BCL}{{\mathbf{L}}}

\nc{\PCC}{{{}^\perp\CC}}

\nc{\Cl}{{\mathsf{Cliff}}}
\nc{\Clev}{{\mathop{\mathsf{Cliff}}^{\circ}}}
\nc{\FA}{{\mathfrak{A}}}
\nc{\FB}{{\mathfrak{B}}}
\nc{\FF}{{\mathfrak{F}}}
\nc{\FI}{{\mathfrak{I}}}
\nc{\FZ}{{\mathfrak{Z}}}

\nc{\TFA}{{\tilde{\mathfrak{A}}}}
\nc{\TFB}{{\tilde{\mathfrak{B}}}}

\nc{\fa}{{\mathfrak{a}}}
\nc{\fg}{{\mathfrak{g}}}
\nc{\fp}{{\mathfrak{p}}}
\nc{\FD}{{\mathfrak{D}}}
\nc{\FE}{{\mathfrak{E}}}
\nc{\FL}{{\mathfrak{L}}}
\nc{\FM}{{\mathfrak{M}}}
\nc{\FR}{{\mathfrak{R}}}
\nc{\FS}{{\mathsf{S}}}

\nc{\sfc}{{\mathsf{c}}}
\nc{\sfch}{{\mathsf{ch}}}
\nc{\sfh}{{\mathsf{h}}}

\nc{\SK}{{\mathsf{K}}}
\nc{\SO}{{\mathsf{O}}}
\nc{\SQ}{{\mathsf{Q}}}
\nc{\SPV}{{\mathsf{S}^+\mathsf{V}}}
\nc{\SMV}{{\mathsf{S}^-\mathsf{V}}}
\nc{\SPMV}{{\mathsf{S}^\pm\mathsf{V}}}
\nc{\SX}{{S_X}}
\nc{\SY}{{S_Y}}
\nc{\phipsi}{{q}}
\nc{\eps}{\varepsilon}

\nc{\pim}{{\pi_-}}
\nc{\pip}{{\pi_+}}

\nc{\BE}{{{\mathbf E}}}
\nc{\TD}{{\widetilde{\D}}}
\nc{\TFD}{{\widetilde{\FD}}}
\nc{\TE}{{\tilde{\CE}}}
\nc{\TQ}{{\tilde{Q}}}
\nc{\TCF}{{\tilde{\CF}}}
\nc{\TCG}{{\tilde{\CG}}}
\nc{\TCH}{{\tilde{\CH}}}
\nc{\TCL}{{\tilde{\CL}}}
\nc{\TF}{{\tilde{F}}}
\nc{\TW}{{\tilde{W}}}
\nc{\TCB}{{\widetilde{\CB}}}
\nc{\TCC}{{\tilde{\CC}}}
\nc{\TCX}{{\tilde{\CX}}}
\nc{\TCY}{{\tilde{\CY}}}
\nc{\TPhi}{{\tilde{\Phi}}}
\nc{\OPhi}{{\bar{\Phi}}}
\nc{\txi}{{\tilde{\xi}}}
\nc{\tp}{{\tilde{p}}}
\nc{\tq}{{\tilde{q}}}
\nc{\tzeta}{{\tilde{\zeta}}}
\nc{\tpi}{{\tilde{\pi}}}

\nc{\HCB}{{\widehat{\CB}}}
\nc{\HCU}{{\widehat{\CU}}}
\nc{\HE}{{\widehat{\CE}}}
\nc{\HS}{{\widehat{S}}}
\nc{\HX}{{\hat{X}}}
\nc{\HY}{{\hat{Y}}}
\nc{\HZ}{{\hat{Z}}}
\nc{\hxi}{{\hat{\xi}}}

\nc{\UH}{{\mathcal{H}}}

\nc{\TM}{{\widetilde{M}}}
\nc{\TCM}{{\widetilde{\CM}}}
\nc{\TS}{{\widetilde{S}}}
\nc{\TT}{{\widetilde{T}}}
\nc{\TU}{{\widetilde{U}}}
\nc{\TX}{{\widetilde{X}}}
\nc{\TY}{{\widetilde{Y}}}
\nc{\TZ}{{\widetilde{Z}}}
\nc{\TYO}{{{\widetilde{Y}}^\circ}}
\nc{\barf}{{\bar{f}}}
\nc{\te}{{\tilde{e}}{}}
\nc{\tf}{{\tilde{f}}}
\nc{\tg}{{\tilde{g}}}
\nc{\ti}{{\tilde{\imath}}}
\nc{\tj}{{\tilde{\jmath}}}
\nc{\ty}{{\tilde{y}}}
\nc{\tphi}{{\tilde{\phi}}}
\nc{\hf}{{\hat{f}}}

\nc{\urho}{{\underline{\rho}}}

\nc{\LRA}{\Leftrightarrow}
\nc{\RA}{\Rightarrow}
\nc{\lotimes}{\mathbin{\mathop{\otimes}\limits^{\mathbb{L}}}}
\nc{\CEnd}{\mathop{\mathcal{E}\mathit{nd}}\nolimits}
\nc{\CExt}{\mathop{\mathcal{E}\mathit{xt}}\nolimits}
\nc{\CHom}{\mathop{\mathcal{H}\mathit{om}}\nolimits}
\nc{\RH}{\mathop{{\mathsf{R}}\Gamma}\nolimits}
\nc{\RGamma}{\mathop{{\mathsf{R}}\Gamma}\nolimits}
\nc{\RHom}{\mathop{\mathsf{RHom}}\nolimits}
\nc{\RCHom}{\mathop{\mathsf{R}\mathcal{H}\mathit{om}}\nolimits}
\nc{\RG}{\mathop{\mathsf{R\Gamma}}\nolimits}
\nc{\Hom}{\mathop{\mathsf{Hom}}\nolimits}
\nc{\Ext}{\mathop{\mathsf{Ext}}\nolimits}
\nc{\End}{\mathop{\mathsf{End}}\nolimits}
\nc{\Tor}{\mathop{\mathsf{Tor}}\nolimits}
\nc{\Tordim}{\mathop{\mathsf{Tor}\text{\rm-}\mathsf{dim}}\nolimits}
\nc{\Hilb}{\mathop{\mathsf{Hilb}}\nolimits}
\nc{\Spec}{\mathop{\mathsf{Spec}}\nolimits}
\nc{\Proj}{\mathop{\mathsf{Proj}}\nolimits}
\nc{\Pic}{\mathop{\mathsf{Pic}}\nolimits}
\renewcommand{\Im}{\mathop{\mathsf{Im}}\nolimits}
\nc{\Tw}{\mathop{\mathsf{Tw}}\nolimits}
\nc{\Ker}{\mathop{\mathsf{Ker}}\nolimits}
\nc{\Coker}{\mathop{\mathsf{Coker}}\nolimits}
\nc{\codim}{\mathop{\mathsf{codim}}\nolimits}
\nc{\sing}{{\mathsf{sing}}}
\nc{\supp}{\mathop{\mathsf{supp}}}
\nc{\perf}{{\mathsf{perf}}}
\nc{\rank}{\mathop{\mathsf{rank}}}
\nc{\Pf}{{\mathsf{Pf}}}
\nc{\Gr}{{\mathsf{Gr}}}
\nc{\OGr}{{\mathsf{OGr}}}
\nc{\Flag}{{\mathsf{Fl}}}
\nc{\Kosz}{{\mathsf{Kosz}}}
\nc{\LGr}{{\mathsf{LGr}}}
\nc{\GTGr}{{\mathsf{G_2Gr}}}
\nc{\GTF}{{\mathsf{G_2F}}}
\nc{\OF}{{\mathsf{OF}}}
\nc{\Fl}{{\mathsf{Fl}}}
\nc{\Bl}{{\mathsf{Bl}}}
\nc{\GL}{{\mathsf{GL}}}
\nc{\PGL}{{\mathsf{PGL}}}
\nc{\SL}{{\mathsf{SL}}}
\nc{\SP}{{\mathsf{Sp}}}
\nc{\Spin}{{\mathsf{Spin}}}
\nc{\Tot}{{\mathsf{Tot}}}
\nc{\ev}{{\mathsf{ev}}}
\nc{\od}{{\mathsf{odd}}}
\nc{\coev}{{\mathsf{coev}}}
\nc{\id}{{\mathsf{id}}}
\nc{\opp}{{\mathsf{opp}}}
\nc{\PS}{{{\PP^3}}}
\nc{\Qu}{{{Q^3}}}
\nc{\tdim}{\mathop{\Tor\dim}}
\nc{\ecart}{{\fbox{$\scriptstyle\mathsf{EC}$}}}
\nc{\ad}{{\mathop{\mathsf ad}}}
\nc{\gr}{{\mathop{\mathsf gr}}}
\nc{\qgr}{{\mathop{\mathsf qgr}}}
\nc{\tor}{{\mathop{\mathsf tor}}}
\rnc{\mod}{{\mathop{\mathsf mod}}}
\nc{\Mod}{{\mathop{\mathsf Mod}}}
\nc{\Coh}{{\mathop{\mathsf Coh}}}
\nc{\Ab}{{\mathop{\mathcal{A}\mathit{b}}}}
\nc{\QCoh}{{\mathop{\mathsf QCoh}}}

\nc{\AAV}{{\mathcal{AAV}}}

\nc{\Rep}{{\mathsf{Rep}}}

\nc{\Cubics}{{{\mathcal{S}}_3}}
\nc{\VFT}{{{\mathcal{S}}_{14}}}
\nc{\VFTE}{{{\mathcal{N}}_{\mathrm{reg,sm}}}}
\nc{\MX}{{\CM_X}}
\nc{\MY}{{\CM_Y}}
\nc{\MYE}{{\CM_{Y,\CE}}}
\nc{\Yd}{{Y_d}}
\nc{\Yfive}{{Y_5}}
\nc{\Xg}{{X_{2g-2}}}
\nc{\Xtt}{{X_{22}}}
\nc{\Xst}{{X_{16}}}
\nc{\Xtw}{{X_{12}}}
\nc{\Xe}{{X_{8}}}
\nc{\Xf}{{X_{4}}}

\nc{\git}{{/\!\!/\!{}_\chi}}

\theoremstyle{plain}

\newtheorem{theo}{Theorem}[]
\newtheorem{conj}[theo]{Conjecture}
\newtheorem{theorem}{Theorem}[section]

\newtheorem{lemma}[theorem]{Lemma}
\newtheorem{proposition}[theorem]{Proposition}
\newtheorem{corollary}[theorem]{Corollary}

\theoremstyle{definition}

\newtheorem{definition}[theorem]{Definition}

\theoremstyle{remark}

\newtheorem{remark}[theorem]{Remark}

\newenvironment{proof}{\noindent{\sf Proof:}}{\qed\medskip}

\title[Homological projective duality for Grassmannians of lines]%
{Homological projective duality for Grassmannians of lines}
\author{Alexander Kuznetsov}
\address{\sloppy
\parbox{0.9\textwidth}{
Algebra Section, Steklov Mathematical Institute,
8 Gubkin str., Moscow 119991 Russia
\hfill\\[5pt]
The Poncelet Laboratory, Independent University of Moscow
\hfill\\[5pt]
}}
\email{akuznet@@mi.ras.ru}
\date{}
\thanks{I was partially supported by RFFI grant 05-01-01034,
Russian Presidential grant for young scientists No. MK-6122.2006.1,
CRDF Award No. RUM1-2661-MO-05 and gratefully acknowledge the support
of the Pierre Deligne fund based on his 2004 Balzan prize in mathematics.}

\begin{document}

\begin{abstract}
We show that homologically projectively dual varieties for Grassmannians $\Gr(2,6)$
and $\Gr(2,7)$ are given by certain noncommutative resolutions
of singularities of the corresponding Pfaffian varieties.
As an application we describe the derived categories of linear sections
of these Grassmannians and Pfaffians.
In particular, we show that

(1) the derived category of a Pfaffian cubic 4-fold admits a semiorthogonal
decompositions consisting of 3 exceptional line bundles, and of the derived category
of a K3-surface;

(2) mutually orthogonal Calabi-Yau linear sections of $\Gr(2,7)$ and
of the corresponding Pfaffian variety are derived equivalent.

We also conjecture a rationality criterion for cubic 4-folds
in terms of their derived categories.
\end{abstract}

\maketitle

\section{Introduction}

Derived categories of coherent sheaves on algebraic varieties and their
semiorthogonal decompositions lately attract a lot of interest.
The most powerful method to produce such decompositions
was introduced in~\cite{K2}. It allows to describe derived categories
of all complete linear sections of a given algebraic variety $X$
if its Homologically Projectively Dual variety is known.
In this paper we find Homological Projectively Dual varieties
for the Grassmannians $\Gr(2,6)$ and $\Gr(2,7)$ and produce
the corresponding semiorthogonal decompositions for
their linear sections.

The Homological Projective Duality (HP-duality for short) is a homological extension
of the classical notion of projective duality. To a smooth (noncommutative) algebraic
variety $X$ with a map $X \to \PP(V)$ to a projective space it
associates a smooth (noncommutative) algebraic variety $Y$ with a map
$Y \to \PP(V^*)$ into the dual projective space, depending on
a choice of a semiorthogonal decomposition of $\D^b(X)$
of a specific form (a Lefschetz decomposition).
All necessary definitions can be found in section~\ref{prelim}.
Now we describe the HP-dual varieties for $\Gr(2,6)$ and $\Gr(2,7)$.

Let $W$ be a vector space, $\dim W = n$.
Consider the projective space $\PP(\Lambda^2W^*)$ of skew-forms on~$W$.
For each $0 \le t \le \lfloor n/2\rfloor$ we consider the following
closed subset of $\PP(\Lambda^2W^*)$
$$
\Pf(2t,n) = \Pf(2t,W^*) := \PP(\{ \omega \in \Lambda^2W^*\ |\ \rank(\omega) \le 2t \}),
$$
where $\rank(\omega)$ is the rank of $\omega$ (the dimension of the image
of the map $W \to W^*$ induced by $\omega$). This subsets form a filtration
of $\PP(\Lambda^2W^*)$. It is clear that $\Pf(2\lfloor n/2\rfloor,W^*) = \PP(\Lambda^2W^*)$.
The biggest proper component of this filtration is $\Pf(2\lfloor n/2\rfloor - 2,W^*)$,
we call it {\sf the Pfaffian variety}, or simply the {\sf Pfaffian}.
If $n = \dim W$ is even then the Pfaffian variety $\Pf(n - 2,W^*)$
is a hypersurface of degree $n/2$ (its equation is the Pfaffian polynomial of a general skew-form),
and if $n = \dim W$ is odd then the Pfaffian variety $\Pf(n-3,W^*)$ has codimension $3$
(its ideal is generated by Pfaffians of principal minors of a general skew-form).
Other varieties $\Pf(2t,W^*)$ will be called {\sf generalized Pfaffian varieties}.

It is a classical result that the Pfaffian variety $Y = \Pf(2\lfloor n/2\rfloor - 2,W^*)$
is (classically) projectively dual to the Grassmannian $X = \Gr(2,W)$.
However, it is singular along $\Pf(2\lfloor n/2\rfloor - 4,W^*)$, so it cannot be HP-dual to $X$.
By some speculations based on the properties of HP-duality one can argue
that the HP-dual of $X$ should be given by a (noncommutative)
resolution of singularities of $Y$. Usual (commutative) resolutions
of $Y$ turn out to be too big, so noncommmutative resolutions
come into focus.

A noncommutative resolution of singularities of an algebraic variety $Y$
is a coherent sheaf $\CR$ of $\CO$-algebras on $Y$ which is a matrix
algebra at the generic point of $Y$ (birationality, up to Morita equivalence),
and which has finite homological dimension (smoothness).
Such noncommutative resolutions of Pfaffian varieties $\Pf(4,6)$ and $\Pf(4,7)$ were constructed in~\cite{K5}.
The main result of the present paper is the following

\begin{theo}\label{th1}
Let $W$ be a vector space, $\dim W = 6$ or $\dim W = 7$.
Let $(Y,\CR)$ be the noncommutative resolution of singularities of the Pfaffian variety $Y = \Pf(4,W^*)$
constructed in~\cite{K5}. Then $(Y,\CR)$ is Homologically Projectively Dual to the Grassmannian $X = \Gr(2,W)$.
\end{theo}

As we already mentioned above, whenever we have a pair of HP-dual varieties,
there are semiorthogonal decompositions for their linear sections.
In our case we obtain semiorthogonal decompositions for linear sections
of the Grassmannians $\Gr(2,6)$, $\Gr(2,7)$ and of the corresponding Pfaffians $\Pf(4,6)$, $\Pf(4,7)$.
The following particular case of these decompositions seems to be especially interesting.

\begin{theo}\label{th2}
Let $Y'$ be a smooth Pfaffian cubic $4$-fold in $\PP^5$.
Then there is a semiorthogonal decomposition
$$
\D^b(Y') = \langle \CO(-3),\CO(-2),\CO(-1),\D^b(X') \rangle
$$
where $X'$ is a smooth K$3$-surface of degree $14$.
\end{theo}

By definition, a Pfaffian cubic fourfold is an intersection
of the Pfaffian $Y = \Pf(4,6) \subset \PP^{14}$ with
a linear subspace $\PP^5 \subset \PP^{14}$
(not every cubic 4-fold is isomorphic to a Pfaffian cubic,
the Pfaffian cubics form a divisor in the moduli space of all cubic 4-folds).
The corresponding K3-surface $X'$ is then the intersection of the orthogonal
linear subspace $\PP^8 \subset (\PP^{14})^\vee$ with the Grassmannian $\Gr(2,6)$.

A relation of cubic 4-folds to K3-surfaces has been known for a long time.
For example, they have similar primitive Hodge structures in the middle cohomology \cite{H2}.
Also, it is known that the Fano variety of lines on a cubic 4-fold
is a deformation of the Hilbert scheme of length 2 subschemes
on a K3 surface~\cite{BD}. For Pfaffian cubics, this relation is more explicit.
The primitive Hodge structure of $X'$ is a substructure of the primitive Hodge structure of $Y'$,
and the Fano variety of $Y'$ is isomorphic to the Hilbert scheme of~$X'$~\cite{BD}.
Moreover, the Pfaffian cubics are known to be rational, and
the birational transformation from~$\PP^4$ to $Y'$ includes a blow-up
of a K3-surface isomorphic to $X'$.

Note that the line bundles $\CO(-3),\CO(-2),\CO(-1)$ form an exceptional collection
on any cubic 4-fold $Y'$ (not necessarily Pfaffian). So, for every $Y'$ we can consider
a triangulated category
$$
\CC' = {}^\perp\langle \CO(-3),\CO(-2),\CO(-1) \rangle \subset \D^b(Y'),
$$
which by Theorem~\ref{th2} is equivalent to the derived category of a K3-surface if $Y'$ is Pfaffian,
and thus for general $Y'$, being a deformation of the derived category of a K3-surface,
can be considered as the derived category of a noncommutative K3-surface.
These speculations suggest the following

\begin{conj}
A smooth cubic $4$-fold $Y'$ is rational if and only if the triangulated category
$$
\CC' = {}^\perp\langle \CO(-3),\CO(-2),\CO(-1) \rangle \subset \D^b(Y')
$$
is equivalent to the derived category of a usual K$3$-surface.
\end{conj}

This conjecture is supported by the following observation.
Assume that $Y'$ is a smooth cubic 4-fold, containing a plane $\PP^2 \subset Y' \subset \PP^5$.
The projection from this plane makes $Y'$ into a family of 2-dimensional quadrics
parameterized by another $\PP^2$. The degenerate quadrics in the family are parameterized
by a sextic curve, and the sheaf of even parts of Clifford algebras corresponding to this
family of quadrics comes from a sheaf of quaterinionic algebras $\CQ$ on the double
cover of $\PP^2$ branched in the sextic curve, which is a K3-surface.
This K3-surface ringed with a sheaf of agebras $\CQ$ can be considered
as a twisted K3-surface, which is a particular class of noncommutative K3-surfaces.
Using results of~\cite{K3} one can show that $\CC'$ is equivalent to the derived
category of this twisted K3. On the other hand, one can check that the class of $\CQ$
in the Brauer group vanishes (and noncommutativity of the K3-surface vanishes as well)
precisely when the corresponding family of 2-dimensional quadrics has a rational section,
which, as was noticed by B.Hasset in~\cite{H1}, implies rationality of the cubic~$Y'$.

\medskip

Another interesting application of Theorem~\ref{th1} is the following.
Let $\dim W = 7$, consider the Grassmannian $X = \Gr(2,7) = \Gr(2,W) \subset \PP(\Lambda^2W)$ and
the corresponding Pfaffian variety $Y = \Pf(4,7) = \Pf(4,W^*) \subset \PP(\Lambda^2W^*)$.
Let $L \subset \Lambda^2W^*$ be a 7-dimensional subspace and
let $L^\perp \subset \Lambda^2W$ be its orthogonal subspace (14-dimensional).
Denote by $X_L = X \cap \PP(L^\perp)$ and $Y_L = Y \cap \PP(L)$ the corresponding
linear sections. Note that the expected dimension of both $X_L$ and $Y_L$ is $3$.

\begin{theo}\label{th3}
If $\dim X_L = \dim Y_L = 3$ and $\PP(L) \subset \PP(\Lambda^2W^*)$
doesn't intersect $\Gr(2,W^*) \subset \PP(\Lambda^2W^*)$, then there is
an equivalence of categories
$$
\D^b(X_L) \cong \D^b(Y_L).
$$
\end{theo}

If the condition $\dim X_L = \dim Y_L = 3$ holds then both $X_L$ and $Y_L$ are Calabi-Yau.
So, the above Theorem gives an example of derived-equivalence between non-birational Calabi-Yau 3-folds.
This equivalence was predicted by E.R{\o}dland~\cite{R}, who compared solutions of the Picard--Fuchs
equations corresponding to these Calabi-Yau families and argued that they have the same mirror.
Another argument from the point of view of string theory was given recently by K.Hori and D.Tong \cite{HT}.
Theorem~\ref{th3} has been independently proved in a particular case of smooth $X_L$ and $Y_L$
by L.Borisov and A.C\u{a}ld\u{a}raru in~\cite{BC} by a direct calculation.
The equivalence in~\cite{BC} is given by the same functor as in our proof.

\medskip

Now let us say a few words about possible generalizations.
First of all, it is natural to consider the case of $\dim W > 7$.
Then we expect the following

\begin{conj}\label{conj2}
Let $X = \Gr(2,W) \subset \PP(\Lambda^2W)$ and
let $Y = \Pf(2\lfloor n/2\rfloor - 2,W^*) \subset \PP(\Lambda^2W^*)$
be the corresponding Pfaffian variety. There exists a noncommutative resolution
of singularities $(Y,\CR)$ of $Y$ which is Homologically Projectively Dual to $X$.
\end{conj}

We expect the noncommutative resolution $(Y,\CR)$ of $Y = \Pf(2\lfloor n/2\rfloor - 2,W^*)$ for $\dim W > 7$
to be constructed in a similar way as for $\dim W = 6,7$, see section~\ref{ncrpf}.
Moreover, we expect that our proof of Theorem~\ref{th1} should work
in the general situation as well.
Among other applications, Conjecture~\ref{conj2} would lead to a description of
derived categories of Pfaffian hypersurfaces of all degrees.

Another direction of generalization is to consider generalized Pfaffian varieties.

\begin{conj}
For every $0\le t\le \lfloor n/2 \rfloor$ there exist noncommutative resolutions
of singularities of the generalized Pfaffian varieties $X = \Pf(2t,W)$ and
$Y = \Pf(2\lfloor n/2\rfloor - 2t,W^*)$ which are Homologically Projectively Dual.
\end{conj}

Now let us describe the structure of the paper.
In section~\ref{prelim} we introduce the necessary background,
reminding the notions of semiorthogonal and Lefschetz decompositions
and giving a brief overview of Homological Projective Duality.
In section~\ref{ncrpf} we describe the noncommutative resolutions
of singularities of the Pfaffian varieties $\Pf(4,6)$ and $\Pf(4,7)$.
In section~\ref{s_main} we give a precise formulation and a plan of the proof of Theorem~\ref{th1}.
The detailed proof takes sections~\ref{s_lc}--\ref{uhps}.
Finally, in sections~\ref{s_gr26} and~\ref{s_gr27} we list the corollaries
of the Homological Projective Duality in cases $\dim W = 6$
and $\dim W = 7$.
In particular, we prove Theorem~\ref{th2} in section~\ref{s_gr26} and
Theorem~\ref{th3} in section~\ref{s_gr27}.

%\bigskip

{\bf Acknowledgements:}
I am grateful to A.~Bondal, D.~Kaledin, and D.~Orlov for very helpful discussions.
Also I would like to thank K.~Hori who pointed out to me the R{\o}dland's work.

\section{Preliminaries}\label{prelim}

We start this section with a brief reminder of the notions of
Lefschetz decompositions and Homological Projective Duality.
The general reference for this is~\cite{K2}.

\subsection{Notation}

All algebraic varieties are assumed to be of finite type over
an algebraically closed field~$\kk$.
% of zero characteristic.
For an algebraic variety $X$, we denote by $\D^b(X)$
the bounded derived category of coherent sheaves on $X$,
and by $\D^-(X)$ the unbounded from below derived category
of coherent sheaves.
For $F,G \in \D^-(X)$, we denote by $\RCHom(F,G)$ the local $\RCHom$-complex
and by $F \otimes G$ the derived tensor product.
Similarly, for a map $f:X \to Y$, we denote by
$f_*$ the derived pushforward functor and by $f^*$ the derived pullback.
Finally, $f^!$ stands for the twisted pullback functor.

\subsection{Semiorthogonal decompositions}

If $\CA$ is a full subcategory of $\CT$ then the {\sf right orthogonal}\/ to $\CA$ in $\CT$
(resp.\ the {\sf left orthogonal}\/ to $\CA$ in $\CT$) is the full subcategory
$\CA^\perp$ (resp.\ ${}^\perp\CA$) consisting of all objects $T\in\CT$
such that $\Hom_\CT(A,T) = 0$ (resp.\ $\Hom_\CT(T,A) = 0$) for all $A \in \CA$.

For any sequence of subcategories $\CA_1,\dots,\CA_n$ in $\CT$ we denote
by $\lan\CA_1,\dots,\CA_n\ran$ the minimal triangulated subcategory of $\CT$
containing $\CA_1$, \dots, $\CA_n$.

\begin{definition}[\cite{BK,BO1,BO2}]
A sequence $\CA_1,\dots,\CA_n$ of full triangulated subcategories in a triangulated
category $\CT$ is called {\sf semiorthogonal collection}\/ if $\Hom_{\CT}(\CA_i,\CA_j) = 0$ for $i > j$.
A semiorthogonal collection $\CA_1,\dots,\CA_n$ is a {\sf semiorthogonal decomposition}\/ of $\CT$
if for every object $T \in \CT$ there exists a chain of morphisms
$0 = T_n \to T_{n-1} \to \dots \to T_1 \to T_0 = T$ such that
the cone of the morphism $T_k \to T_{k-1}$ is contained in $\CA_k$
for each $k=1,2,\dots,n$. In other words, if there exists a diagram
\begin{equation*}\label{tower}
\vcenter{
\xymatrix@C-7pt{
0 \ar@{=}[r] & T_n \ar[rr]&& T_{n-1} \ar[dl] \ar[rr]&& \quad\dots\quad \ar[rr]&& T_2 \ar[rr]&& T_1 \ar[dl] \ar[rr]&& T_0 \ar[dl] \ar@{=}[r] & T \\
&& A_n \ar@{..>}[ul] &&& \dots &&& A_2 \ar@{..>}[ul]&& A_1 \ar@{..>}[ul]&&
}}
\end{equation*}
where all triangles are distinguished (dashed arrows have degree $1$) and $A_k \in \CA_k$.
\end{definition}

Thus, every object $T\in\CT$ admits a decreasing ``filtration''
with factors in $\CA_1$, \dots, $\CA_n$ respectively.
Semiorthogonality implies that this filtration is unique and functorial.

If $\CT = \lan\CA_1,\dots,\CA_n\ran$ is a semiorthogonal decomposition then
$\CA_i = {}^\perp\lan\CA_1,\dots,\CA_{i-1}\ran \cap \lan\CA_{i+1},\dots,\CA_n\ran{}^\perp$.

\begin{definition}[\cite{BK,B}]
A full triangulated subcategory $\CA$ of a triangulated category $\CT$ is called
{\sf right admissible}\/ if the inclusion functor $i:\CA \to \CT$ has
a right adjoint $i^!:\CT \to \CA$, and
{\sf left admissible}\/ if it has a left adjoint $i^*:\CT \to \CA$.
Subcategory $\CA$ is called {\sf admissible}\/ if it is both right and left admissible.
\end{definition}

\begin{lemma}[\cite{B}]\label{sod_adm}
If $\CT = \lan\CA,\CB\ran$ is a semiorthogonal decomposition then
$\CA$ is left amissible and $\CB$ is right admissible.
If\/ $\CA_1,\dots,\CA_n$ is a semiorthogonal collection in $\CT$
such that $\CA_1,\dots,\CA_k$ are left admissible and
$\CA_{k+1},\dots,\CA_n$ are right admissible then
$$
\lan\CA_1,\dots,\CA_k,
{}^\perp\lan\CA_1,\dots,\CA_k\ran \cap \lan\CA_{k+1},\dots,\CA_n\ran{}^\perp,
\CA_{k+1},\dots,\CA_n\ran
$$
is a semiorthogonal decomposition.
\end{lemma}

Let $f:X \to S$ be a morphism of algebraic varieties.
Triangulated subcategory $\CA \subset \D^b(X)$ is called {\sf $S$-linear}\/ \cite{K1}
if it is stable with respect to tensoring by pull-backs of vector bundles on $S$:
$$
\CA \otimes f^* F \subset \CA
\qquad\text{for any vector bundle $F$ on $S$}.
$$

If $E_1,E_2,\dots,E_n$ is a sequence of objects of $\D^b(X)$
of finite $\Tor$-dimension over $S$ we denote by
$\langle E_1,E_2,\dots,E_n\rangle_S$ the $S$-linear triangulated subcategory
of $\D^b(X)$ generated by $E_1,E_2,\dots,E_n$. In other words,
$\langle E_1,E_2,\dots,E_n\rangle_S =
\langle E_1\otimes f^*\D^b(S),E_2\otimes f^*\D^b(S),\dots,E_n\otimes f^*\D^b(S)\rangle$.

A semiorthogonal collection $\lan \CA_1,\CA_2,\dots,\CA_n\ran \subset \D^b(X)$
is {\sf $S$-linear}\/ if all subcategories $\CA_k$ are $S$-linear.
If $\CA$ is an $S$-linear admissible subcategory in $\D^b(X)$ then
both $\CA^\perp$ and ${}^\perp\CA$ are $S$-linear and semiorthogonal decompositions
$\D^b(X) = \langle\CA,{}^\perp\CA\rangle = \langle\CA^\perp,\CA\rangle$ are $S$-linear.

Let $f:X \to S$ and $g:Y \to S$ be morphisms to $S$.
A functor $\Phi:\D^b(X) \to \D^b(Y)$ is called {\sf $S$-linear}\/ \cite{K1},
if there is given a bifunctorial isomorphism
$$
\Phi(G \otimes f^*F) \cong \Phi(G) \otimes g^*F,
$$
where $G \in \D^b(X)$ and $F$ is a vector bundle on $S$.
If an $S$-linear functor $\Phi:\D^b(X) \to \D^b(Y)$
has an (either left or right) adjoint functor $\Psi:\D^b(Y) \to \D^b(X)$
then $\Psi$ is also $S$-linear.

\subsection{Kernel functors}\label{kf}

Let $X$ and $Y$ be smooth projective varieties and
let $p_X:X\times Y \to X$ and $p_Y:X\times Y \to Y$ denote the projections.
Take any $K\in\D^b(X\times Y)$ and define a functor
$$
\Phi_K : \D^b(Y) \to \D^b(X),
\qquad
\Phi_K(F) := {P_X}_*(p_Y^*F \otimes K).
$$
We call $\Phi_K$ the {\sf kernel functor}\/ with kernel $K$.

\begin{lemma}[\cite{BO1}]\label{ladjoint}
The functor $\Phi_K:\D^b(Y) \to \D^b(X)$ admits a left adjoint functor $\Phi_K^*$ which is isomorphic
to a kernel functor $\Phi_{K^\#}:\D^b(X) \to \D^b(Y)$ with the kernel
$K^\#:= \RCHom(K,\omega_X[\dim X])$.
\end{lemma}

Consider kernels
$K \in \D^b(X\times Y)$,
$L \in \D^b(Y\times Z)$.
Denote by $p_{XY}$, $p_{YZ}$ and $p_{XZ}$ the projections
of $X\times Y\times Z$ to $X\times Y$, $Y\times Z$ and $X\times Z$ respectively.
The {\sf convolution}\/ of kernels is defined as follows
$$
K\circ L := {p_{XZ}}_*(p_{XY}^*K\otimes p_{YZ}^*L).
$$
It is well known \cite{BO1,BO2} that the composition of kernel functors is given
by the convolution of their kernels
$$
\Phi_K\circ\Phi_L \cong \Phi_{K\circ L},
\qquad
\Phi_L^*\circ\Phi_K^* \cong \Phi_{K\circ L}^*.
$$

Let $f:X \to S$ and $g:Y \to S$ be morphisms to $S$. Let $j:X\times_S Y \to X\times S$ be the embedding.
If $K \in \D^b(X\times_S Y)$ then the functors
$\Phi_{j_*K}:\D^b(Y) \to \D^b(X)$ and $\Phi_{j_*K}^*:\D^b(X) \to \D^b(Y)$ are $S$-linear.

If $\CA \subset \D^b(X)$ is a triangulated subcategory and $S$ is a scheme we denote
by $\CA \boxtimes \D^b(S)$ the minimal triangulated subcategory of $\D^b(X\times S)$
containing all objects of the form $A\boxtimes F$, where $A\in \CA$, $G \in \D^b(S)$.

\begin{lemma}\label{boxtimes}
If $X$ is a smooth projective variety, $\CA \subset \D^b(X)$ is an admissible subcategory
equivalent to $\D^b(Y)$ for some smooth projective variety $Y$ and $S$ is any scheme then
$\CA \boxtimes \D^b(S) \cong \D^b(Y\times S)$ is an admissible $S$-linear
subcategory in $\D^b(X\times S)$.
\end{lemma}
\begin{proof}
By theorem of D.Orlov \cite{O1,O2} any equivalence $\D^b(Y) \to \CA \subset \D^b(X)$ is isomorphic
to a kernel functor $\Phi_K:\D^b(Y) \to \D^b(X)$, $K \in \D^b(Y\times X)$.
Let $\Delta:Y\times S\times X \to Y\times S\times X\times S$ be the embedding
induced by the diagonal embedding of $S$. Then it is easy to see that $\Delta_*(K\boxtimes\CO_S)$ gives
an $S$-linear fully faithful functor $\D^b(Y\times S) \to \D^b(X\times S)$ inducing an equivalence
$\D^b(Y\times S) \to \CA\boxtimes \D^b(S)$.
\end{proof}

\subsection{Lefschetz decompositions}

Let $X$ be an algebraic variety and let $\CO_X(1)$ be a line bundle on~$X$.

\begin{definition}[\cite{K2}]
A {\sf Lefschetz decomposition} of the derived category $\D^b(X)$ is a semiorthogonal decomposition
of the form
$$
\D^b(X) = \lan \CA_0,\CA_1(1),\dots,\CA_{m-1}(m-1)\ran,
\qquad\text{where $0 \subset \CA_{m-1} \subset \dots \subset \CA_1 \subset \CA_0 \subset \D^b(X)$.}
$$
Similarly,
a {\sf dual Lefschetz decomposition}\/ of $\D^b(X)$ is a semiorthogonal decomposition
of the form
$$
\D^b(X) = \lan \CA_{m-1}(1-m),\dots,\CA_1(-1),\CA_0\ran,
\qquad\text{where $0 \subset \CA_{m-1} \subset \dots \subset \CA_1 \subset \CA_0 \subset \D^b(X)$.}
$$
\end{definition}

Actually, these notions are equivalent. Given a Lefschetz decomposition
one can canonically construct a dual Lefschetz decomposition with
the same category $\CA_0$ and vice versa.

Let $\D^b(X) = \lan \CA_0,\CA_1(1),\dots,\CA_{m-1}(m-1)\ran$ be a Lefschetz decomposition.
Let $\fa_k$ denote the right orthogonal to $\CA_{k+1}$ in $\CA_k$, so that
we have a semiorthogonal decomposition
$$
\CA_k = \langle\fa_k,\fa_{k+1},\dots,\fa_{m-1}\rangle.
$$
The categories $\fa_0,\fa_1,\dots,\fa_{m-1}$ are called {\sf primitive}
categories of the Lefschetz decomposition.

Denote by $\alpha_0:\CA_0 \to \D^b(X)$ the embedding functor and
let $\alpha_0^*:\D^b(X) \to \CA_0$ be its left adjoint.

\begin{lemma}[\cite{K2}]
The functor $\alpha_0^*$ is fully faithful on subcategories $\fa_k(k+1) \subset \D^b(X)$
and there is a semiorthogonal decomposition
$\CA_0 = \langle \alpha_0^*(\fa_0(1)),\alpha_0^*(\fa_1(2)),\dots,\alpha_0^*(\fa_{m-1}(m))\rangle$.
\end{lemma}
We call this semiorthogonal decomposition of $\CA_0$ the {\sf dual primitive decomposition}.

The simplest example of a Lefschetz decomposition is given by the Beilinson exceptional collection
on the projective space $\PP^n$:
$$
\D^b(\PP^n) = \langle \CO,\CO(1),\dots,\CO(n) \rangle.
$$
Here $\CA_0 = \CA_1 = \dots = \CA_n = \langle \CO \rangle$.
The primitive categories are $\fa_0 =\dots = \fa_{n-1} = 0$ and $\fa_n = \langle \CO \rangle$.

More relevant for the present paper are the following Lefschetz decompositions of the derived
categories of Grassmannians $\Gr(2,W)$ of lines in a vector space $W$.
Let $\CU$ denote the tautological rank $2$ vector bundle on $\Gr(2,W)$.
Then
$
\D^b(\Gr(2,W)) = \langle \CA_0,\CA_1(1),\dots,\CA_{m-1}(m-1) \rangle,
$
where $m = \dim W$ and
\begin{equation*}\label{ldgr}
\hspace{-25pt}
\begin{array}{ll}
\CA_0 = \CA_1 = \dots = \CA_{m-1} = \langle S^{k-1}\CU,\dots,\CU,\CO \rangle, & \text{if $m = 2k + 1$}\\
\CA_0 = \dots = \CA_{k-1} = \langle S^{k-1}\CU,\dots,\CU,\CO \rangle,\quad
\CA_k = \dots = \CA_{2k-1} = \langle S^{k-2}\CU,\dots,\CU,\CO \rangle, & \text{if $m = 2k$.}
\end{array}
\hspace{-20pt}
\end{equation*}
These decompositions were constructed in~\cite{K4}.
The primitive subcategories here are
$$
\begin{array}{ll}
\fa_0 = \fa_1 = \dots = \fa_{m-2} = 0, \quad
\fa_{m-1} = \langle S^{k-1}\CU,\dots,\CU,\CO \rangle, & \text{if $m = 2k + 1$}\\
\fa_0 = \dots = \fa_{k-2} = \fa_k = \dots = \fa_{2k-2} = 0, \
\fa_{k-1} = \langle S^{k-1}\CU \rangle, \
\fa_{2k-1} = \langle S^{k-2}\CU,\dots,\CU,\CO \rangle, & \text{if $m = 2k$.}
\end{array}
$$
The dual primitive decomposition takes form
$$
\CA_0 = \CA_0\quad\text{for $m = 2k + 1$},\qquad
\CA_0 = \langle \Lambda^{k-1}(W/\CU),{}^\perp_{\CA_0}(\Lambda^{k-1}(W/\CU)) \rangle\quad\text{for $m = 2k$.}
$$

\subsection{Homological projective duality}

Fix a smooth projective variety $X$ and a Lefschetz decomposition
$\D^b(X) = \langle \CA_0,\CA_1(1),\dots,\CA_{\ix-1}(\ix-1) \rangle$
with respect to a line bundle $\CO_X(1)$.
Let $f:X \to \PP(V)$ be a morphism into a projective space such that
$f^*(\CO_{\PP(V)}(1)) \cong \CO_X(1)$ and let $\CX \subset X\times\PP(V^*)$
be the universal hyperplane section of $X$
(i.e. the canonical divisor of bidegree $(1,1)$ in $X\times\PP(V^*)$).

\begin{definition}[\cite{K2}]
An algebraic variety $Y$ with a projective morphism $g:Y\to\PP(V^*)$
is called {\sf Homologically Projectively Dual} to $f:X\to \PP(V)$
with respect to the given Lefschetz decomposition, if
there exists an object $\CE\in\D^b(\CX\times_{\PP(V^*)} Y)$ such that
the kernel functor
$\Phi = \Phi_\CE:\D^b(Y) \to \D^b(\CX)$ is fully faithful
and gives the following semiorthogonal decomposition
\begin{equation*}\label{dbx1}
\D^b(\CX) = \langle \Phi(\D^b(Y)),\CA_1(1)\boxtimes\D^b(\PP(V^*)),\dots,
\CA_{\ix-1}(\ix-1)\boxtimes\D^b(\PP(V^*))\rangle.
\end{equation*}
\end{definition}

For every linear subspace $L \subset V^*$ we consider
the corresponding linear sections of $X$ and $Y$:
$$
X_L = X\times_{\PP(V)}\PP(L^\perp),
\qquad
Y_L = Y\times_{\PP(V^*)}\PP(L),
$$
where $L^\perp \subset V$ is the orthogonal subspace to $L\subset V^*$.
Let $N = \dim V$.

The main property of Homologically Projectively Dual varieties is the following

\begin{theorem}[\cite{K2}]\label{hp}
If\/ $Y$ is Homologically Projectively Dual to $X$ then

\noindent$(i)$
$Y$ is smooth and $\D^b(Y)$ admits a dual Lefschetz decomposition
$$
\D^b(Y) = \lan \CB_{\jx-1}(1-\jx),\dots,\CB_{1}(-1),\CB_0\ran,\qquad
0 \subset \CB_{\jx-1} \subset \dots \subset \CB_1 \subset \CB_0 \subset \D^b(Y)
$$
with the same set of primitive subcategories:
$\CB_k = \langle\fa_0,\dots,\fa_{N-k-2}\rangle$;

\noindent$(ii)$
for any linear subspace $L\subset V^*$, $\dim L = r$, such that we have
$\dim X_L = \dim X - \dim L$, and $\dim Y_L = \dim Y + \dim L - N$\/
there exist a triangulated category $\CC_L$ and semiorthogonal decompositions
$$
\begin{array}{lll}
\D^b(X_L) &=& \langle \CC_L,\CA_{r}(1),\dots,\CA_{\ix-1}(\ix-r)\rangle\\
\D^b(Y_L) &=& \langle \CB_{\jx-1}(N-r-\jx),\dots,\CB_{N-r}(-1),\CC_L\rangle.
\end{array}
$$
\end{theorem}

We will need below the following necessary and sufficient condition
for an algebraic variety $Y$ to be Homologically Projectively Dual to $X$.

Let $g:Y \to \PP(V^*)$ be a regular map.
Note that the map
$\CX\times_{\PP(V^*)} Y \subset (X\times\PP(V^*))\times_{\PP(V^*)} Y \cong X\times Y$
identifies $\CX\times_{\PP(V^*)} Y$ with a divisor of bidegree $(1,1)$ in $X\times Y$
which we call the {\sf incidence quadric} and denote by $Q(X,Y)$. Let $j$ denote
the embedding $Q(X,Y) \cong \CX\times_{\PP(V^*)} Y \to \CX\times Y$.
Assume that we are given an object $\CE \in \D^b(Q(X,Y))$ such that
the kernel functor $\Phi_{j_*\CE}:\D^b(Y) \to \D^b(\CX)$ is a fully faithful
embedding into a subcategory
$$
\CC =
[\langle \CA_1(1)\boxtimes\D^b(\PP(V^*)),\dots,
\CA_{\ix-1}(\ix-1)\boxtimes\D^b(\PP(V^*))\rangle]^\perp.
$$
In~\cite{K2} there was constructed the following dual Lefschetz collection
in the category $\CC$:
$$
\langle \CB'_{\jx-1}(1-\jx),\CB'_{\jx-2}(2-\jx),\dots,\CB'_1(-1),\CB'_0\rangle \subset \CC,
$$
where $\jx = N - 1 - \max\{ k\ |\ \CA_k = \CA_0\}$ and
\begin{equation}\label{defcb}
\CB'_k = \gamma^*\pi^*
\Big\langle\alpha_0^*(\fa_0(1)),\dots,\alpha_0^*(\fa_{N-k-2}(N-k-1))\Big\rangle
 \subset \CB'_0 = \gamma^*\pi^*\CA_0,
\end{equation}
the image of the part $\langle\alpha_0^*(\fa_0(1)),\dots,\alpha_0^*(\fa_{N-k-2}(N-k-1))\rangle$
of the dual primitive decomposition of the category $\CA_0$ under the functor
$\gamma^*\pi^*:\CA_0 \subset \D^b(X) \to \D^b(\CX) \to \CC$,
where $\gamma^*:\D^b(\CX) \to \CC$ is the left adjoint functor
to the embedding functor $\gamma:\CC \to \D^b(\CX)$, and
$\pi:\CX \to X$ is the projection.

\begin{theorem}\label{hpd_crit}
Assume that the functor $\Phi_{j_*\CE}:\D^b(Y) \to \D^b(\CX)$ induces a fully faithful
embedding $\D^b(Y) \to \CC$. Assume additionally that the functor
$\Phi_{j_*\CE}^*\pi^*:\D^b(X) \to \D^b(Y)$ is fully faithful
on the components $\alpha_0^*(\fa_k(k+1)) \subset \CA_0$
of the dual primitive decomposition of $\CA_0$ and that the categories
$$
\CB_k = \Phi_{j_*\CE}^*\left(\pi^*
\Big\langle\alpha_0^*(\fa_0(1)),\dots,\alpha_0^*(\fa_{N-k-2}(N-k-1))\Big\rangle
\right) \subset \CB_0 = \Phi_{j_*\CE}^*(\pi^*(\CA_0)) \subset \D^b(Y),
$$
form a dual Lefschetz collection
$$
\langle \CB_{\jx-1}(1-\jx),\dots,\CB_1(-1),\CB_0\rangle \subset \D^b(Y).
$$
Then the functor $\Phi_{j_*\CE}:\D^b(Y) \to \CC$ is an equivalence and the above collection
generates $\D^b(Y)$. In particular, $Y$ is Homologically Projectively Dual to $X$.
\end{theorem}
\begin{proof}
Actually, the proof of Theorem~6.3 in section~6 of~\cite{K2}
uses nothing but the assumptions of the present theorem. On the other hand,
at the output of the proof of Theorem~6.3 in~\cite{K2} we obtain an equivalence
$\D^b(Y) \cong \CC$, and the fullness of the above collection in $\D^b(Y)$.
\end{proof}

%\begin{remark}\label{phispis}
%It is easy to see that
%$\Phi_{j_*\CE}^*(\CB_k) =
%\Phi_{j_*\CE}^*(\gamma^*(\pi^*(\langle\alpha_0^*(\fa_0(1)),\dots,\alpha_0^*(\fa_{N-k-2}(N-k-1))\rangle))) =
%\Phi_{j_*\CE}^*\pi^*(\langle\alpha_0^*(\fa_0(1)),\dots,\alpha_0^*(\fa_{N-k-2}(N-k-1))\rangle)$.
%\end{remark}

\subsection{The Borel--Bott--Weil Theorem}\label{ss_bbw}

The Borel--Bott--Weil Theorem computes the cohomology of line bundles
on the flag variety of a semisimple Lie group. We use it to compute
the cohomology of equivariant vector bundles on Grassmannians.
We restrict here to the case of the group $\GL(V)$.

Let $V$ be a vector space of dimension $n$. The standard identification
of the weight lattice of the group $\GL(V)$ with $\ZZ^n$ takes the $k$-th
fundamental weight $\pi_k$ (the heighest weight of the representation $\Lambda^kV$)
to the vector $(1,1,\dots,1,0,0,\dots,0) \in \ZZ^n$ (the first $k$ entries
are 1, and the last $n-k$ are 0). Under this identification
the cone of dominant weights of $\GL(V)$ gets identified with the set of nonincreasing
sequences $\alpha = (a_1,a_2,\dots,a_n)$ of integers. For such $\alpha$
we denote by $\Sigma^\alpha V = \Sigma^{a_1,a_2,\dots,a_n} V$ the corresponding
representation of $\GL(V)$.
Note that $\Sigma^{1,1,\dots,1}V = \det V$.

Similarly, given a vector bundle $E$ of rank $n$ on a scheme $S$
we consider the corresponding principal $\GL(n)$-bundle on $S$ and
denote by $\Sigma^\alpha E$ the vector bundle associated
with the $\GL(n)$-representation of highest weight $\alpha$.

The group $\BS_n$ of permutations acts naturally on the weight lattice $\ZZ^n$.
Denote by $\ell:\BS_n \to \ZZ$ the standard length function.
Note that for every $\alpha \in \ZZ^n$ there exists a permutation $\sigma \in \BS_n$
such that $\sigma(\alpha)$ is nonincreasing. If all entries of $\alpha$
are distinct then such $\sigma$ is unique and $\sigma(\alpha)$ is strictly decreasing.

Let $X$ be the flag variety of $\GL(V)$. Let $L_\alpha$ denote the line bundle
on $X$ corresponding to the weight~$\alpha$ (so that $L_{\pi_k}$ is the pullback
of $\CO_{\PP(\Lambda^kV)}(1)$ under the natural projection $X \to \PP(\Lambda^kV)$).

Denote by
$$
\rho = (n,n-1,\dots,2,1)
$$
half the sum of the positive roots of $\GL(V)$.
The corresponding line bundle $L_\rho$ is the square root of the anticanonical line bundle.

The Borell-Bott-Weil Theorem computes the cohomology of line bundles $L_\alpha$ on $X$.
\begin{theorem}[\cite{D}]\label{bbw}
Assume that all entries of $\alpha + \rho$ are distinct.
Let $\sigma$ be the unique permutation such that $\sigma(\alpha+\rho)$ is strictly decreasing.
Then
$$
H^k(X,L_\alpha) = \begin{cases}
\Sigma^{\sigma(\alpha+\rho) - \rho}V^*, & \text{if $k = \ell(\sigma)$}\\
0, & \text{otherwise}
\end{cases}
$$
If not all entries of $\alpha + \rho$ are distinct then $H^\bullet(X,L_\alpha) = 0$.
\end{theorem}

Now consider a Grassmannian $\BG = \Gr(k,V)$. Let $\CU \subset V\otimes\CO_\BG$ denote
the tautological subbundle of rank $k$. Denote by $W/\CU$ the corresponding quotient bundle
and by $\CU^\perp$ its dual, so that we have the following (mutually dual) exact sequences
$$
0 \to \CU \to W\otimes\CO_\BG \to W/\CU \to 0,
\qquad
0 \to \CU^\perp \to W^*\otimes\CO_\BG \to \CU^* \to 0.
$$
Note that $\Sigma^{1,1,\dots,1}\CU^* \cong \Sigma^{-1,-1,\dots,-1}\CU^\perp$ is the positive
generator of $\Pic\BG$.
Let $\pi:X \to \BG$ denote the canonical projection from the flag variety to the Grassmannian.

\begin{proposition}[\cite{Ka}]\label{equb}
Let $\beta \in \ZZ^k$ and $\gamma \in \ZZ^{n-k}$ be nonincreasing sequences.
Let $\alpha = (\beta,\gamma) \in \ZZ^n$ be their concatenation.
Then we have $\pi_* L_\alpha \cong \Sigma^\beta\CU^*\otimes\Sigma^\gamma\CU^\perp$.
\end{proposition}

\begin{corollary}\label{bbwg}
If $\beta \in \ZZ^k$ and $\gamma \in \ZZ^{n-k}$ are nonincreasing sequences
and $\alpha = (\beta,\gamma) \in \ZZ^n$ then
$$
H^\bullet(\BG,\Sigma^\beta\CU^*\otimes\Sigma^\gamma\CU^\perp) \cong H^\bullet(X,L_\alpha).
$$
\end{corollary}

Note that every $\GL(V)$-equivariant vector bundle on $\BG$ is isomorphic to
$\Sigma^\beta\CU^*\otimes\Sigma^\gamma\CU^\perp$ for some nonincreasing $\beta\in\ZZ^k$,
$\gamma\in\ZZ^{n-k}$. Thus a combination of corollary~\ref{bbwg} with the Borel--Bott--Weil
Theorem allows to compute the cohomology of any equivariant vector bundle on $\BG$.

\section{The Pfaffian varieties and their noncommutative resolutions}\label{ncrpf}

Let $W$ be a vector space over $\kk$, $\dim W = n$.
Consider the projective space $\BP = \PP(\Lambda^2W^*)$ of skew-forms on $W$.
For each $0 \le t \le \lfloor n/2\rfloor$ we consider the following
closed subset of $\BP = \PP(\Lambda^2W^*)$
$$
\Pf(2t,n) = \Pf(2t,W^*) = \PP(\{ \omega \in \Lambda^2W^*\ |\ \rank(\omega) \le 2t \}),
$$
where $\rank(\omega)$ is the rank of $\omega$ (the dimension of the image
of the map $W \to W^*$ induced by $\omega$).
We call $\Pf(2\lfloor n/2\rfloor-2,W^*)$ the {\sf Pfaffian variety},
and other varieties $\Pf(2t,W^*)$ are called the {\sf generalized Pfaffian varieties}.
The ideal of the generalized Pfaffian variety $\Pf(2t,W^*)$ is generated
by the Pfaffians of all diagonal $(2t+2)\times(2t+2)$-minors of a skew-form,
hence the name.

It is clear that the following generalized Pfaffian varieties
$$
\Pf(0,W^*) = \emptyset,\qquad
\Pf(2,W^*) = \Gr(2,W^*),\qquad
\Pf(2\lfloor n/2\rfloor,W^*) = \BP = \PP(\Lambda^2W^*)
$$
are smooth. However, for $t \ne 0,1,\lfloor n/2\rfloor$ the Pfaffian variety
$\Pf(2t,W^*)$ is singular, the singular locus being the previous Pfaffian
$\Pf(2t-2,W^*) \subset \Pf(2t,W^*)$.

In this section we describe a noncommutative resolution of singularities
of the generalized Pfaffian variety $\Pf(4,W^*)$ for $n = \dim W \ge 6$.
So, put
$$
Y = \Pf(4,W^*),
\qquad
Z = \mathop{{\sf Sing}}(Y) = \Pf(2,W^*) = \Gr(2,W^*) = \Gr(n-2,W).
$$
Note that all skew-forms in $Y \setminus Z$ are of rank $4$, hence their kernels are $(n - 4)$-dimensional.
Similarly, all skew-forms in $Z$ are of rank $2$, and their kernels are $(n - 2)$-dimensional.
Let $\TY$ be the set of all pairs $(\omega,K)$, where $K$ is an $(n - 4)$-dimensional subspace
in $W$ and $\omega$ is a skew-form containing $K$ in the kernel.
More precisely, $\TY = \PP_\BG(\Lambda^2\CK^\perp)$,
where $\BG = \Gr(n - 4,W)$ is a Grassmannian, $\CK \subset W\otimes\CO_\BG$ is the tautological subbundle of rank $n - 4$,
and $\CK^\perp \subset W^*\otimes\CO_\BG$ is the orthogonal to $\CK$ subbundle of rank~$4$.
Note that $\TY$ is smooth.

The embedding of vector bundles $\Lambda^2\CK^\perp \to \Lambda^2W^*\otimes\CO_\BG$
induces a projection $g:\TY \to \PP(\Lambda^2W^*) = \BP$. It is clear that it factors through
a map to the Pfaffian variety $g_Y:\TY \to Y \subset \BP$.
The map $g_Y$ is one-to-one over $Y\setminus Z$ and a $\Gr(n-4,n-2)$-fibration over $Z$.
Indeed, let $\TZ = g_Y^{-1}(Z) \subset \TY$.
Note that the bundle of kernels of skew-forms on $Z \cong \Gr(2,W^*) \cong \Gr(n-2,W)$
can be identified with the tautological subbundle $\CK_{n-2} \subset W\otimes\CO_Z$ of rank $n-2$,
hence $\TZ \cong \Gr_Z(n-4,\CK_{n-2})$, the relative Grassmannian.
In the other words, $\TZ \cong \Fl(n-4,n-2;W)$, the partial flag variety.

Let $\zeta:\TY \to \BG$ and $g_Z:\TZ \to Z$ be the projections, and
let $i_\TY:\TZ \to \TY$ and $i_Y:Z \to Y$ be the embeddings.
Let $\eta:Y \to \BP$ be the embedding, so that $g = \eta\circ g_Y$.
Then we have the following commutative diagram of varieties and maps:
$$
\xymatrix{
&& \TY \ar[dl]_-\zeta \ar[dd]_g \ar[dr]^-{g_Y} &&
\TZ \ar[ll]_{i_\TY} \ar[dr]^-{g_Z} &
\Fl(n-4,n-2;W) \ar@{=}[l] \ar[dr] \\
\Gr(n-4,W) \ar@{=}[r] &
\BG &&
Y \ar[dl]^\eta &&
Z \ar[ll]^{i_Y} &
\Gr(n-2,W) \ar@{=}[l] \\
& \PP(\Lambda^2 W^*) \ar@{=}[r] &
\BP
}
%\xymatrix{
%\Gr(2,W) &
%\TY \ar[l]_\zeta \ar[dl]_g \ar[d]^{g_Y} &
%\TZ \ar[l]_{i_\TY} \ar[d]^{g_Z} &
%\Fl(2,4;W) \ar@{=}[l] \ar[d] \\
%\PP(\Lambda^2 W^*) &
%Y \ar[l]^\eta &
%Z \ar[l]^{i_Y} &
%\Gr(4,W) \ar@{=}[l]
%}
$$

We denote by $H_G$ the divisor class of a hyperplane section of $\BG$ and by $H_Y$
the divisor class of a hyperplane section of $Y$. The pullbacks of these classes to $\TY$, $\TZ$
and other varieties are denoted by the same letters.

\begin{lemma}
We have a linear equivalence $\TZ \sim 2H_Y - H_G$ on $\TY$.
In particular, we have an exact sequence
\begin{equation}\label{ohz}
0 \to \CO_{\TY}(H_G-2H_Y) \to \CO_\TY \to i_{\TY*}\CO_\TZ \to 0,
\end{equation}
and an isomorphism
$\CN_{\TZ/\TY} \cong \CO_\TZ(2H_Y-H_G)$.
\end{lemma}
\begin{proof}
Since $\TY = \PP_\BG(\Lambda^2\CK^\perp)$, the Picard group
of $\TY$ is generated by $H_G$ and $H_Y$, hence $\TZ \sim \lambda H_G + \mu H_Y$
for some $\lambda,\mu \in \ZZ$. Moreover, it easy to see that $K_\TY = - (n-3)H_G - 6H_Y$.
Similarly, since $\TZ = \Fl(n-4,n-2;W)$, the Picard group
of $\TZ$ is generated by $H_G$ and $H_Y$, and $K_\TZ = - (n-2)H_G - 4H_Y$.
By adjunction formula we find $\lambda = -1$, $\mu = 2$.
\end{proof}

From now on we restrict ourselves to the cases $\dim W = 6,7$.
See~\cite{K5} for the description of a noncommutative resolution of singularities
of the generalized Pfaffian varieties $\Pf(4,W^*)$ for all $\dim W$.

In~\cite{K5} we have constructed a $Z$-linear dual Lefschetz decomposition with respect to $\CN^*_{\TZ/\TY}$
\begin{equation}\label{hz}
\D^b(\TZ) = \langle \CC_{n-3}((3-n)H_G),\dots,\CC_1(-H_G),\CC_0 \rangle,
%\begin{array}{ll}
%\D^b(\TZ) = \langle \CC_3(-3H_G),\CC_2(-2H_G),\CC_1(-H_G),\CC_0 \rangle & \text{if $n = 6$}\\
%\D^b(\TZ) = \langle \CC_4(-4H_G),\CC_3(-3H_G),\CC_2(-2H_G),\CC_1(-H_G),\CC_0 \rangle & \text{if $n = 7$}
%\end{array}
\end{equation}
where $n = \dim W$ and
\begin{equation}\label{lg24}
\begin{array}{ll}
\CC_0 = \CC_1 = \langle \CO,\CK^* \rangle_Z,\qquad
\CC_2 = \CC_3 = \langle\CO\rangle_Z, & \text{if $n = 6$}\\[2pt]
\CC_0 = \CC_1 = \CC_2 = \CC_3 = \CC_4 = \langle\CO,\CK^*\rangle_Z, & \text{if $n = 7$.}
\end{array}
\end{equation}

\begin{theorem}[\cite{K5}]\label{NCR}
The pushforward functor $i_{\TY*}:\D^b(\TZ) \to \D^b(\TY)$ is full and faithful on subcategories
$\CC_{n-3}((3-n)H_G)$, \dots,  $\CC_1(-H_G)$ of $\D^b(\TZ)$
and there is a $\BP$-linear semiorthogonal decomposition
\begin{equation}\label{ty}
\D^b(\TY) =
\langle
i_{\TY*}(\CC_{n-3}((3-n)H_G)),
\dots,
i_{\TY*}(\CC_1(-H_G)),
\TD
\rangle,
\end{equation}
where
\begin{equation}\label{td}
\TD = \{ F \in \D^b(\TY)\ |\ i_\TY^* F \in \CC_0 \}.
\end{equation}
Moreover, the pushforward
\begin{equation}\label{ca}
\CR = g_{Y*}\CEnd(\CO_\TY \oplus \CK)
\end{equation}
is a pure sheaf,
the category $\Coh(Y,\CR)$ of coherent sheaves of right $\CR$-modules on $Y$
has finite homological dimension, and
there is an equivalence of categories
\begin{equation}\label{tda}
\TD \cong \D^b(Y,\CR),
\end{equation}
where $\D^b(Y,\CR)$ is the bounded derived category of $\Coh(Y,\CR)$.
The equivalence is given by the functors
$$
\begin{array}{l}
\rho_*:\D^b(\TY) \to \D^b(Y,\CR),\qquad
F \mapsto g_{Y*}(F \otimes (\CO_\TY \oplus \CK)),\qquad\text{and}\\[2pt]
\rho^*:\D^b(Y,\CR) \to \D^b(\TY),\qquad
G \mapsto g_Y^{-1}G \otimes_{g_Y^{-1}(\CR)} (\CO_\TY \oplus \CK^*).
\end{array}
$$
Finally, the category $\TD$ admits a Serre functor~(see~\cite{BK}) $\FS_\TD:\TD \to \TD$
and we have
$$
\begin{array}{ll}
\FS_\TD (F) \cong F(-12H_Y)[13], & \text{if $n = 6$ and $i^*F \in \langle\CO_\TZ\rangle_Z$,}\\[2pt]
\FS_\TD (F) \cong F(-14H_Y)[17], & \text{if $n = 7$ and any $F \in \TD$.}
\end{array}
$$
\end{theorem}

\begin{remark}
Note that by~\eqref{ca} the restriction of the sheaf of algebras $\CR$
to $Y \setminus Z$ is isomorphic to a matrix algebra.
This allows to consider the noncommutative variety $\D^b(Y,\CR)$
as a noncommutative resolution of singularities of $Y$.
\end{remark}

Below we will need also a description of the derived category
of left $\CR$-modules (or, equivalently, of right $\CR^\opp$-modules) on $Y$,
and of the derived category of $\CR$-bimodules (equivalently, of right
$\CR\boxtimes\CR^\opp$-modules) on $Y\times Y$ as subcategories
of $\D^b(\TY)$ and $\D^b(\TY\times\TY)$ respectively.

This will be done as follows. We take
\begin{equation}\label{tds}
\arraycolsep = 2pt
\begin{array}{lll}
\TD_\opp &:= \TD^* &=
\{F \in \D^b(\TY)\ |\ F^* \in \TD \} =
\{ F \in \D^b(\TY)\ |\ i_\TY^* F \in (\CC_0)^* =
\langle \CO_\TZ,\CK \rangle_Z  \}\\[5pt]
\TD_\natural &:= \TD\boxtimes\TD_\opp & \subset \D^b(\TY\times\TY).
\end{array}
\end{equation}
The functors
$$
\begin{array}{ll}
\rho_*^\opp:\TD_\opp \to \D^b(Y,\CR^\opp), &
F \mapsto g_{Y*}(F\otimes(\CO_\TY\otimes\CK^*)),\\[2pt]
\rho^*_\opp:\D^b(Y,\CR^\opp) \to \TD_\opp, &
G \mapsto g_Y^{-1}G \otimes _{g_Y^{-1}\CR^\opp}(\CO_\TY\otimes\CK),\quad\text{and}\\[5pt]
\rho_*^\natural:\TD_\natural \to \D^b(Y\times Y,\CR\boxtimes\CR^\opp), &
F \mapsto (g_Y\times g_Y)_*(F\otimes(\CO_\TY\oplus\CK)\boxtimes(\CO_\TY\oplus\CK^*)),\\[2pt]
\rho^*_\natural:\D^b(Y\times Y,\CR\boxtimes\CR^\opp) \to \TD_\natural, &
G \mapsto (g_Y\times g_Y)^{-1}G \otimes_{(g_Y\times g_Y)^{-1}(\CR\boxtimes\CR^\opp)}((\CO_\TY\oplus\CK^*)\boxtimes(\CO_\TY\oplus\CK)),
\end{array}
$$
give equivalences
$$
\TD_\opp \cong \D^b(Y,\CR^\opp),
\qquad
\TD_\natural \cong \D^b(Y\times Y,\CR\boxtimes\CR^\opp).
$$
Equivalences of $\TD$, $\TD_\opp$ and $\TD_\natural$ with $\D^b(Y,\CR)$, $\D^b(Y,\CR^\opp)$ and $\D^b(Y\times Y,\CR\boxtimes\CR^\opp)$
provide these triangulated categories with a t-structure $(\TD^{\le0},\TD^{\ge0})$
\begin{equation}\label{tdtstr}
\begin{array}{l}
\TD^{\le0} = \{ F \in \TD\ |\ g_{Y*}(F\otimes(\CO_\TY \oplus \CK)) \in \D^{\le 0}(Y) \},\\[5pt]
\TD^{\ge0} = \{ F \in \TD\ |\ g_{Y*}(F\otimes(\CO_\TY \oplus \CK)) \in \D^{\ge 0}(Y) \},
\end{array}
\end{equation}
and similarly defined t-structures $(\TD_\opp^{\le0},\TD_\opp^{\ge0})$ and
$(\TD_\natural^{\le0},\TD_\natural^{\ge0})$ in $\TD_\opp$ and $\TD_\natural$ respectively.

\begin{lemma}\label{rhophirho}
Let $K \in \D^b(Y\times Y,\CR\boxtimes\CR^\opp)$. Then
$\rho^*\circ\Phi_K\circ\rho_* \cong \Phi_{\rho_\natural^*K}:\D^b(\TY) \to \D^b(\TY)$.
\end{lemma}
\begin{proof}
Since $\rho_\natural^*K \in \TD_\natural = \TD\boxtimes\TD_\opp$ we have $\Phi_{\rho_\natural^*K}(\D^b(\TY)) \subset \TD$.
Since $\rho_*:\TD \to \D^b(Y,\CR)$ is an equivalence, it suffices to check that
$\Phi_K\circ\rho_* \cong \rho_*\circ\Phi_{\rho_\natural^*K}:\D^b(\TY) \to \D^b(Y,\CR)$.
Denote $\CO_\TY \oplus \CK = E$ for brevity and take any $F \ in \D^b(\TY)$.
Then we have
\begin{multline*}
\rho_*(\Phi_{\rho_\natural^*K}(F)) \cong
g_{Y*}(p_{1*}(p_2^*F \otimes ((g_Y\times g_Y)^{-1} K\otimes_{(g_Y\times g_Y)^{-1}\CR\boxtimes\CR^\opp}(E^*\boxtimes E)))\otimes E) \cong
\\ \cong
p_{1*}(g_Y\times g_Y)_*(((E\otimes E^*)\boxtimes (F\otimes E)) \otimes_{(g_Y\times g_Y)^{-1}\CR\boxtimes\CR^\opp} (g_Y\times g_Y)^{-1}K) \cong
\\ \cong
p_{1*}((g_Y\times g_Y)_*((E\otimes E^*)\boxtimes (F\otimes E)) \otimes_{\CR\boxtimes\CR^\opp} K) \cong
\\ \cong
p_{1*}((\CR\boxtimes \rho_*(F)) \otimes_{\CR\boxtimes\CR^\opp} K) \cong
p_{1*}(p_2^*\rho_*(F) \otimes_{\CR\boxtimes\CR^\opp} K) \cong
\Phi_K(\rho_*(F))
\end{multline*}
(we have used a version of the projection formula in the third isomorphism).
\end{proof}

We also will need relative analogues of these categories. Let $S$ be a smooth algebraic variety. Consider the subcategories
$$
\begin{array}{lll}
\TD_S & := \TD\boxtimes\D^b(S) & \subset \D^b(\TY\times S),\\[2pt]
\TD_{\opp S} & := \D^b(S)\boxtimes\TD_\opp & \subset \D^b(S\times \TY),\quad\text{and}\\[2pt]
\TD_{\natural S} & := \TD\boxtimes\D^b(S)\boxtimes\TD_\opp & \subset \D^b(\TY\times S\times\TY).
\end{array}
$$
The functors
$$\arraycolsep=0pt
\begin{array}{ll}
\rho_*^{S}:\TD_{S} \to \D^b(Y\times S,\CR\boxtimes\CO_S),
&
F \mapsto (g_Y\!\times\! \id_S)_*(F\otimes[(\CO_\TY\oplus\CK)\boxtimes\CO_S]),\\[2pt]
\rho_*^{\opp S}:\TD_{\opp S} \to \D^b(S\times Y,\CO_S\boxtimes\CR^\opp),
&
F \mapsto (\id_S\!\times\! g_Y)_*(F\otimes[\CO_S\boxtimes(\CO_\TY\oplus\CK^*)]),\\[2pt]
\rho_*^{\natural S}:\TD_{\natural S} \to \D^b(Y\times S\times Y,\CR\boxtimes\CO_S\boxtimes\CR^\opp),
\ &
F \mapsto (g_Y\!\times\! \id_S\!\times\! g_Y)_*(F\otimes[(\CO_\TY\oplus\CK)\boxtimes\CO_S\boxtimes(\CO_\TY\oplus\CK^*)]),\\[2pt]
\end{array}
$$
and similarly defined adjoint functors $\rho^*_S$, $\rho^*_{\opp S}$ and $\rho^*_{\natural S}$ by Theorem~6.4 of~\cite{K5}
give equivalences
$$
\begin{array}{lll}
\TD_{S} &\cong& \D^b(Y\times S,\CR\boxtimes\CO_S),\\[2pt]
\TD_{\opp S} &\cong& \D^b(S\times Y,\CO_S\boxtimes\CR^\opp),\\[2pt]
\TD_{\natural S} &\cong& \D^b(Y\times S\times Y,\CR\boxtimes\CO_S\boxtimes\CR^\opp).
\end{array}
$$
Moreover, for any map $\phi:S \to T$ the functors
$\rho_*,\rho^\opp_*,\rho^\natural_*$ and their adjoints
commute with the pushforward $\phi_*$ and the pullback $\phi^*$ functors.

We denote by $(\TD_{S}^{\le 0},\TD_{S}^{\ge 0})$, $(\TD_{\opp S}^{\le 0},\TD_{\opp S}^{\ge 0})$ and
$(\TD_{\natural S}^{\le 0},\TD_{\natural S}^{\ge 0})$
the t-structures on $\TD_S$, $\TD_{\opp S}$ and $\TD_{\natural S}$ induced by the equivalences
$\rho^*_S$, $\rho^*_{\opp S}$ and $\rho^*_{\natural S}$.

Further, $\TD^0 = \TD^{\le 0}\cap \TD^{\ge 0}$ denotes the heart of the t-structure.
Similarly, we write $\TD^{[a,b]} = \TD^{\le b}\cap \TD^{\ge a}$ for all $a\le b$.
The $k$-th cohomology functor with respect to the t-structure $(\TD^{\le 0},\TD^{\ge 0})$
is denote by $\TCH^k$. A similar notation is used for t-structures in categories
$\TD_\opp$, $\TD_\natural$, $\TD_S$, $\TD_{\opp S}$ and $\TD_{\natural S}$.

These t-structures are related as follows.

\begin{lemma}\label{tensort}
Assume that $S$ is smooth.
If $F \in \TD_{S}^{[a,b]} \subset \D^b(\TY\times S)$ and $F' \in \TD_{\opp S}^{[a',b']} \subset \D^b(S\times \TY)$ then
$p_{12}^*F\otimes p_{23}^*F' \in \TD_{\natural S}^{[a+a'-\dim S,b+b']}$, where $p_{12}$ and $p_{23}$ are the projections
of $\TY\times S\times\TY$ to $\TY\times S$ and $S\times\TY$ respectively.
\end{lemma}
\begin{proof}
First of all note that $\rho_{\natural S}^*(p_{12}^*G\otimes p_{23}^*G') \cong p_{12}^*\rho_{S}^*(G)\otimes p_{23}^*\rho_{\opp S}^*(G')$.
On the other hand, for any objects $G \in \D^{[a,b]}(Y\times S,\CR\boxtimes\CO_S)$, $G' \in \D^{[a',b']}(S\times Y,\CO_S\boxtimes\CR^\opp)$
it is easy to see that
$p_{12}^*G\otimes p_{23}^*G' \in \D^{[a+a'-\dim S,b+b']}(Y\times S\times Y,\CR\boxtimes\CO_S\boxtimes\CR^\opp)$
since $S$ is smooth.
\end{proof}

\begin{lemma}\label{mapt}
Let $\phi:S \to T$ be a morphism with fibers of dimension not exceeding
$n$.
Then $F \in \TD_{S}^{[a,b]}$ implies $\phi_* F \in \TD_{T}^{[a,b+n]}$.
Similarly for $\TD_{\opp S}$ and $\TD_{\natural S}$.
\end{lemma}
\begin{proof}
Indeed, $F\in \D^{[a,b]}(Y\times S,\CR\otimes\CO_S)$
implies $\phi_*F\in \D^{[a,b+n]}(Y\times T,\CR\otimes\CO_T)$
since the cohomological dimension of $\phi_*$ is $n$.
\end{proof}

\begin{lemma}\label{finitet}
Let $\phi:S \to T$ be a morphism and assume that $F \in \TD_S \subset \D^b(\TY\times S)$
has finite support over $\TY\times T$. Then
$F \in \TD_{S}^{[a,b]}$ is equivalent to $\phi_* F \in \TD_{T}^{[a,b]}$
and $\TCH^k(\phi_*F) \cong \phi_*\TCH^k(F)$ for all $k$.
Similarly for $\TD_{\opp S}$ and $\TD_{\natural S}$.
\end{lemma}
\begin{proof}
Since $F$ has finite support over $\TY\times T$ it follows that
$\CH^k(\rho_*^T(\phi_*F)) =
\CH^k(\phi_*(\rho_*^SF)) =
\phi_*\CH^k(\rho_*^SF)$,
which means that $\TCH^k(\phi_*F) \cong \phi_*\TCH^k(F)$.
\end{proof}

%The functor
%$$
%\rho^\opp_*:\D^b(\TY) \to \D^b(Y,\CA^\opp_{Y}),\qquad
%F \mapsto g_{Y*}(F \otimes (\CO_\TY \oplus \CK^*))
%$$
%gives an equivalence of $\TD^\opp$ with the derived category $\D^b(Y,\CR^\opp)$
%of sheaves of modules over the opposite algebra $\CR^\opp = g_{Y*}\CEnd(\CO_\TY\oplus\CK^*)$
%such that the following diagram commutes
%$$
%\xymatrix{
%\TD \ar[rrr]^{\RCHom(-,\CO_\TY)} \ar[d]_{\rho_*} &&&
%\TD^\opp \ar[d]^{\rho^\opp_*} \\
%\D^b(Y,\CR) \ar[rrr]^{\RCHom(-,\CO_Y)} &&&
%\D^b(Y,\CA^\opp_{Y})
%}
%$$

%\begin{comment}

\section{The Main Theorem}\label{s_main}

Let $W$ be a vector space, $\dim W = 6$ or $\dim W = 7$.
Consider the Grassmannian $X = \Gr(2,W)$ and the Pfaffian variety $Y = \Pf(4,W^*)$.
Let $(Y,\CR)$ be the noncommutative resolution of singularities of $Y$ constructed
in the previous section. We consider $\D^b(Y,\CR)$ as a subcategory of $\D^b(\TY)$
via the equivalence~\eqref{tda}. Then the map $g = \eta\circ g_Y:\TY \to Y \to \BP = \PP(\Lambda^2W^*)$
can be considered as a morphism of the noncommutative variety $(Y,\CR)$ to the projective space
$\BP = \PP(\Lambda^2W^*)$, which we denote by the same letter $g$.
So, consider the maps $f:X \to \BP^\vee = \PP(\Lambda^2 W)$ (the Pl\"ucker embedding)
and $g:(Y,\CR) \to \BP$.

The main result of the paper, Theorem~\ref{th1} from the Introduction, claims the Homological
Projective Duality between $X$ and $Y$. To make this statement precise we have to specify
the involved Lefschetz decompositions of $D^b(X)$ and $\D^b(Y,\CR)$ (or at least one of them).

Consider the following exceptional triple on $X$:
\begin{equation}\label{e012}
E_0 = \CO_X,\qquad\qquad
E_1 = \CU,\qquad\qquad
E_2 = S^2\CU.
\end{equation}
Here $\CU \subset W\otimes\CO_X$ is the tautological subbundle of rank $2$ on $X$.
It was shown in~\cite{K4} that this triple generates the following
Lefschetz decomposition for $\D^b(X)$
\begin{equation}\label{ldx}
\begin{array}{ll}
\D^b(X) = \langle\CA_0,\CA_1(H_X),\CA_2(2H_X),\CA_3(3H_X),\CA_4(4H_X),\CA_5(5H_X)\rangle,\\[2pt]
\CA_0 = \CA_1 = \CA_2 = \langle E_2,E_1,E_0 \rangle,
\quad
\CA_3 = \CA_4 = \CA_5 = \langle E_1,E_0 \rangle, & \text{if $n = 6$}\\[5pt]
\D^b(X) = \langle\CA_0,\CA_1(H_X),\CA_2(2H_X),\CA_3(3H_X),\CA_4(4H_X),\CA_5(5H_X),\CA_6(6H_X)\rangle,\\[2pt]
\CA_0 = \CA_1 = \CA_2 = \CA_3 = \CA_4 = \CA_5 = \CA_6 = \langle E_2,E_1,E_0 \rangle, & \text{if $n = 7$}
\end{array}
\end{equation}
where $H_X$ is the divisor class of a hyperplane section of $X$.

On the other hand, consider the following triple of bundles on $\TY$
\begin{equation}\label{f012def}
F_0 = \Lambda^2(W/\CK) / \CO_\TY(H_G - H_Y),\qquad\qquad
F_1 = W/\CK,\qquad\qquad F_2 = \CO_\TY,
\end{equation}
where $\CK \subset W\otimes\CO_\TY$ is the pullback of the tautological rank $n-4$ subbundle on
$\BG = \Gr(n-4,W)$ via the projection $\zeta:\TY \to \BG$,
$H_Y$ is the divisor class of a hyperplane section of $Y$,
$H_G$ is the divisor class of a hyperplane section of $\BG$,
and the embedding
$\CO_\TY(H_G - H_Y) \to \Lambda^2(W/\CK)$ is obtained by a $\CO_\TY(H_G)$-twist from the embedding
$\CO_\TY(-H_Y) \to \Lambda^2\CK^\perp = \Lambda^2 (W/\CK) \otimes \CO_\TY(-H_G)$
inducing the map $g:\TY = \PP_\BG(\Lambda^2\CK^\perp) \to \BP = \PP(\Lambda^2W^*)$.

We will show in proposition~\ref{ldcb} below that the dual triple of~\eqref{f012def} generates
the following dual Lefschetz collection in~$\TD = \D^b(Y,\CR)$:
\begin{equation}\label{ldtd}
\begin{array}{ll}
\D^b(Y,\CR) = \TD \supset \langle \CB_{11}(-11H_Y),\CB_{10}(-10H_Y),\dots,\CB_1(-H_Y),\CB_0 \rangle,\\[2pt]
\CB_0 = \CB_1 = \dots = \CB_8 = \langle F_0^*,F_1^*,F_2^* \rangle,
\qquad
\CB_9 = \CB_{10} = \CB_{11} = \langle F_2^* \rangle, & \text{if $n=6$}\\[5pt]
\D^b(Y,\CR) = \TD \supset \langle \CB_{13}(-13H_Y),\CB_{12}(-12H_Y),\dots,\CB_1(-H_Y),\CB_0 \rangle,\\[2pt]
\CB_0 = \CB_1 = \dots = \CB_{13} = \langle F_0^*,F_1^*,F_2^* \rangle, & \text{if $n=7$}
\end{array}
\end{equation}

Now we give a precise statement of Theorem~\ref{th1}.

\begin{theorem}\label{themain}
If $\dim W = 6$ or $\dim W = 7$ then noncommutative resolution of singularities $(Y,\CR)$
of the Pfaffian variety $Y = \Pf(4,W^*)$ is Homologically Projectively Dual
to the Grassmannian $X = \Gr(2,W)$ with respect to the Lefschetz decomposition~\eqref{ldx}.
The corresponding Lefschetz decomposition of $\D^b(Y,\CR)$ is given by~\eqref{ldtd}.
\end{theorem}

The proof of theorem~\ref{themain} takes sections~\ref{s_comp}--\ref{uhps}.
Let us briefly describe the principal steps.

First of all we have to construct an object $\CE \in \D^b(\CX\times_{\BP}\TY)$
which gives a fully faithful functor $\D^b(Y,\CR) = \TD \subset \D^b(\TY) \stackrel{\Phi_\CE}\to \D^b(\CX)$,
where $\CX \subset X \times \BP$ is the universal hyperplane section of $X$.
For this we will show that on $X \times \TY$ there is a natural complex of vector bundles
$$
E_2\boxtimes F_2 \to E_1\boxtimes F_1 \to E_0\boxtimes F_0,
$$
and that this complex is quasiisomorphic to a coherent sheaf $\CE$
supported on the incidence quadric $Q(X,\TY) \subset X\times\TY$
($Q(X,\TY)$ is the preimage of the usual incidence quadric
$Q \subset \BP^\vee\times\BP$ under the projection
$f\times g:X\times\TY \to \BP^\vee\times\BP$).
The canonical isomorphism $\CX\times_{\BP}\TY \cong Q(X,\TY)$
allows to regard $\CE$ as a kind of object we need.

Let $j$ denote the embedding $Q(X,\TY) \cong \CX\times_\BP\TY \to \CX \times \TY$.
The most difficult part of the proof is to verify that the functor
$\Phi_{j_*\CE}:\D^b(\TY) \to \D^b(\CX)$ induces a fully faithful embedding
of $\TD$ into $\D^b(\CX)$. A straightforward way is to compute the composition
$\Phi_{j_*\CE}^* \circ \Phi_{j_*\CE} : \D^b(\TY) \to \D^b(\TY)$.
We do this using the following trick.

Let $\alpha:\CX \to X\times\BP$ be the embedding.
Note that the functor $\alpha_*\circ\Phi_{j_*\CE}:\D^b(\TY) \to \D^b(X\times\BP)$
is given by the kernel $i_*\CE \in \D^b(X\times\TY) = \D^b((X\times\BP)\times_\BP\TY)$,
where $i:Q(X,\TY) \to X\times\TY$ is the embedding.
On the other hand, since $\alpha$ is a divisorial embedding,
we have a distinguished triangle of functors $\D^b(\CX) \to \D^b(\CX)$
$$
\alpha^*\alpha_* \to \id \to \CO_\CX(-H_X-H_P)[2]
$$
(the last term denotes the functor of the $\CO_\CX(-H_X-H_P)$-twisting followed
by the $[2]$-shift). Composing this with the functor $\Phi_{j_*\CE}$ on the right
and with the functor $\Phi_{j_*\CE}^*$ on the left we obtain a distinguished triangle
of functors $\D^b(\TY) \to \D^b(\TY)$
$$
(\Phi_{j_*\CE}^*\alpha^*)\circ(\alpha_*\Phi_{j_*\CE}) \to
\Phi_{j_*\CE}^*\circ\Phi_{j_*\CE} \to
\Phi_{j_*\CE}^*\circ\Phi_{j_*\CE(-H_X-H_P)}[2].
$$
Taking any $t \in \ZZ$ and twisting by $\CO(tH_X + tH_P)$ we obtain the following distinguished triangles
$$
(\Phi_{j_*\CE}^*\alpha^*)\circ(\alpha_*\Phi_{j_*\CE(tH_X+tH_P)}) \to
\Phi_{j_*\CE}^*\circ\Phi_{j_*\CE(tH_X+tH_P)} \to
\Phi_{j_*\CE}^*\circ\Phi_{j_*\CE((t-1)H_X+(t-1)H_P)}[2].
$$
The point is that the first term of these triangles can be computed quite easily
using the resolution $i_*\CE \cong \{E_2\boxtimes F_2 \to E_1\boxtimes F_1 \to E_0\boxtimes F_0\}$.
In particular, it is easy to see that
$$
(\Phi_{j_*\CE}^*\alpha^*)\circ(\alpha_*\Phi_{j_*\CE(tH_X+tH_P)}) = 0
\qquad
\text{for $1\le t\le 6$}
$$
in the case $\dim W = 7$ to which we restrict from this moment in this short explanation of the proof
(in the case $\dim W = 6$ the arguments are slightly different but of the same spirit).
It follows immediately that
$$
\Phi_{j_*\CE}^*\circ\Phi_{j_*\CE(6H_X+6H_P)} \cong
\Phi_{j_*\CE}^*\circ\Phi_{j_*\CE(5H_X+5H_P)}[2] \cong \dots \cong
\Phi_{j_*\CE}^*\circ\Phi_{j_*\CE}[12],
$$
hence we have a distinguished triangle
$$
(\Phi_{j_*\CE}^*\alpha^*)\circ(\alpha_*\Phi_{j_*\CE(7H_X+7H_P)}) \to
\Phi_{j_*\CE}^*\circ\Phi_{j_*\CE(7H_X+7H_P)} \to
\Phi_{j_*\CE}^*\circ\Phi_{j_*\CE}[14].
$$
On the other hand, we can find an estimate for the set of $k \in \ZZ$ such that
the $k$-the cohomology $\TCH^k$ of the kernel of the functor
$(\Phi_{j_*\CE}^*\alpha^*)\circ(\alpha_*\Phi_{j_*\CE(7H_X+7H_P)})$
is nonzero (let us call this set {\sf the cohomology support interval}),
and a uniform (in $t$) estimate of the cohomology support intervals of the kernels of the functors
$\Phi_{j_*\CE}^*\circ\Phi_{j_*\CE(tH_X+tH_P)}$. The $[14]$-shift in the above triangle
makes the cohomology support intervals of the kernels of the functors in the last triangle intersect
only at one point, which means in particular that the kernel of the functor
$\Phi_{j_*\CE}^*\circ\Phi_{j_*\CE}$ is a pure object isomorphic to the $(-13)$-th cohomology
of the kernel of the functor $(\Phi_{j_*\CE}^*\alpha^*)\circ(\alpha_*\Phi_{j_*\CE(7H_X+7H_P)})$.
A direct computation allows to identify this cohomology with the object in $\D^b(\TY\times\TY)$
inducing the projection functor to the subcategory $\TD \subset \D^b(\TY)$.
This finally shows that $\Phi_{j_*\CE}$ is the composition of the projection
$\D^b(\TY) \to \TD$ and a fully faithful embedding $\TD \to \D^b(\CX)$.

It is worth emphasizing that in the above arguments we always use the t-structure $(\TD_\natural^{\le 0},\TD_\natural^{\ge 0})$
in the category $\TD_\natural \subset \D^b(\TY\times\TY)$ constructed in the previous section.
The standard t-structure of $\D^b(\TY\times\TY)$ is not sufficiently sharp and doesn't work here.

The final step in the proof uses theorem~\ref{hpd_crit}.
According to this theorem it remains to check that the functor
$\Phi_{j_*\CE}^*\circ\pi^*:\D^b(X) \to \D^b(\TY)$
is fully faithful on the subcategory $\CA_0 \subset \D^b(X)$
and that its image $\CB_0 \subset \D^b(\TY)$ generates a dual Lefschetz
collection in $\D^b(\TY)$. But the functor $\Phi_{j_*\CE}^*\circ\pi^*$
is the left adjoint of $\pi_*\circ\Phi_{j_*\CE}:\D^b(\TY) \to \D^b(X)$
which is a kernel functor with the kernel $i_*\CE \in \D^b(X\times\TY)$.
Again, using the resolution $i_*\CE \cong \{E_2\boxtimes F_2 \to E_1\boxtimes F_1 \to E_0\boxtimes F_0\}$
it is easy to perform all required verifications.

The proof is spread between sections~\ref{s_lc}--\ref{uhps} as follows.
In section~\ref{s_lc} we show that~\eqref{ldtd} is a Lefschetz collection in $\TD$.
%and compute the cohomological amplitude of the objects $F_k,F_k^*$.
In section~\ref{s_comp} we compute the pushforwards of some objects on $\TY$
to $\BP$. These computations are used later in section~\ref{s_fd}
to identify a cohomology of some object in $\D^b(\TY\times\TY)$
with the kernel of the projection functor $\D^b(\TY) \to \TD$.
In section~\ref{s_ker} we construct the kernel $\CE$ by showing that there exists
a natural complex $E_2\boxtimes F_2 \to E_1\boxtimes F_1 \to E_0\boxtimes F_0$
and checking that it is quasiisomorphic to $i_*\CE$ for some~$\CE$.
Finally, in section~\ref{uhps} we finish the proof.

\section{A Lefschetz collection for the Pfaffian varieties}\label{s_lc}

Recall that we have defined in~\eqref{f012def} a triple of vector bundles $(F_0,F_1,F_2)$ on $\TY$.

\begin{lemma}\label{fintd}
We have $F_0,F_1,F_2 \in \TD_\opp$, $F_0^*,F_1^*,F_2^* \in \TD$.
%The objects $F_0,F_1,F_2$ are contained in the subcategory $\TD^\opp$ of $\D^b(\TY)$.
\end{lemma}
\begin{proof}
By definition~\eqref{tds} of the category $\TD_\opp$ it suffices to verify the first inclusion.
So, by~\eqref{tds} we just have to check that the restrictions of $F_0$, $F_1$ and $F_2$ to $\TZ = \Fl(2,4;W)$
are contained in the subcategory $\CC_0^* = \langle \CO_\TZ,\CK\rangle_Z$.
For $F_2 = \CO_\TY$ this is evident.
For $F_1$ the restriction to $\TZ$ of the exact sequence
$$
0 \to \CK \to W\otimes\CO_\TY \to F_1 \to 0
$$
also gives the desired embedding. It remains to consider $F_0$.
By definition we have an exact sequence
\begin{equation}\label{f0}
0 \to \CO_\TY(H_G-H_Y) \to \Lambda^2(W/\CK) \to F_0 \to 0.
\end{equation}
Consider its restriction to the divisor $\TZ = \Fl(n-4,n-2;W)$.
Let $\CK_{n-2}$ denote (the pullbacks to $\TZ$ of) the tautological
subbundle in $W\otimes\CO_{\Gr(n-2,W)}$ of rank $n-2$. Since we have
$\CO_\TZ(-H_Y) = \Lambda^2\CK_{n-2}^\perp$, $\Lambda^2(W/\CK) = \Lambda^2\CK^\perp\otimes\CO_\TZ(H_G)$,
the restriction takes form
$$
0 \to \Lambda^2\CK_{n-2}^\perp \otimes \CO_\TZ(H_G) \to
\Lambda^2\CK^\perp \otimes \CO_\TZ(H_G) \to F_{0|\TZ} \to 0.
$$
The first map here is induced by the embedding $\CK_{n-2}^\perp \subset \CK^\perp$ on $\TZ = \Fl(n-4,n-2;W)$.
Therefore, we have the following exact sequence
$$
0 \to
\CK_{n-2}^\perp\otimes(\CK^\perp/\CK_{n-2}^\perp)\otimes\CO_\TZ(H_G) \to
F_{0|\TZ} \to
\Lambda^2(\CK^\perp/\CK_{n-2}^\perp)\otimes\CO_\TZ(H_G) \to 0.
$$
It remains to note that
$$
\begin{array}{l}
\Lambda^2(\CK^\perp/\CK_{n-2}^\perp) \cong
\det\CK^\perp \otimes (\det\CK_{n-2}^\perp)^{-1} \cong
\CO_\TZ(-H_G)\otimes\CO_\TZ(H_Y),\\[2pt]
\CK^\perp/\CK_{n-2}^\perp \cong
(\CK^\perp/\CK_{n-2}^\perp)^*\otimes\det(\CK^\perp/\CK_{n-2}^\perp) \cong
(\CK_{n-2}/\CK)\otimes\CO_\TZ(-H_G)\otimes\CO_\TZ(H_Y)
\end{array}
$$
and the claim follows.
\end{proof}

Our next goal is to show that the triple $(F_0^*,F_1^*,F_2^*)$ in $\TD$ is exceptional
and to describe the subcategory of $\TD$ generated by this triple.
We start with the following

\begin{lemma}
The quadruple
$$
(\CO_\TY(H_G-H_Y),\CO_\TY,W/\CK,\Lambda^2(W/\CK))
$$
in $\D^b(\TY)$ is exceptional.
Moreover, for $0\le k\le l\le 2$ we have $\Hom^\bullet(\Lambda^k(W/\CK),\Lambda^l(W/\CK)) \cong \Lambda^{l-k}W$
and
$$
\Hom^\bullet(\CO_\TY(H_G-H_Y),\CO_\TY) =
\Hom^\bullet(\CO_\TY(H_G-H_Y),W/\CK) = 0,\quad
\Hom^\bullet(\CO_\TY(H_G-H_Y),\Lambda^2(W/\CK)) = \kk.
%\Hom(\CO_\TY(H_G-H_Y),\Lambda^k(W/\CK)) = \begin{cases}
%0, & \text{for $k=0,1$}\\
%\kk, & \text{for $k=2$}\\
%\end{cases}
$$
In other words, the algebra of endomorphisms of this exceptional quadruple is the path algebra of the quiver
$$
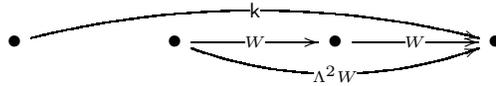
\xymatrix@C=4pc{
\bullet \ar@/^1pc/[rrr]|{\kk} &
\bullet \ar[r]|{\scriptscriptstyle W} \ar@/_1pc/[rr]|{\scriptscriptstyle \Lambda^2W} &
\bullet \ar[r]|{\scriptscriptstyle W} &
\bullet
}
$$
\end{lemma}
\begin{proof}
The triple $(\CO_\TY,W/\CK,\Lambda^2(W/\CK))$ is exceptional since it is exceptional in $\D^b(\BG)$
and the pullback functor $\zeta^*:\D^b(\BG) \to \D^b(\TY)$ is fully faithful. Further, for all $k$ we have
$$
\Hom^\bullet(\Lambda^k(W/\CK),\CO_\TY(H_G-H_Y)) =
\Hom^\bullet(\Lambda^k(W/\CK),\zeta_*(\CO_\TY(H_G-H_Y))) = 0
$$
since $\zeta_*(\CO_\TY(H_G-H_Y)) = 0$.
Finally,
\begin{multline*}
\Hom^\bullet(\CO_\TY(H_G-H_Y),\Lambda^k(W/\CK)) \cong
\Hom^\bullet(\CO_\TY,\Lambda^k(W/\CK)(H_Y-H_G)) \cong
\\ \cong
H^\bullet(\Gr(2,W),\zeta_*(\Lambda^k(W/\CK)(H_Y-H_G))) \cong
H^\bullet(\Gr(2,W),\Lambda^k(W/\CK)\otimes\Lambda^2(W/\CK)\otimes\CO(-H_G))) \cong
\\ \cong
H^\bullet(\Gr(2,W),\Lambda^k(W/\CK)\otimes\Lambda^2\CK^\perp) \cong
\Hom^\bullet(\Lambda^2(W/\CK),\Lambda^k(W/\CK))
\end{multline*}
and the last claim follows.
\end{proof}

\begin{corollary}
The triple $(F_2,F_1,F_0)$ in $\TD_\opp$ is exceptional. Moreover,
$\Hom(F_2,F_1) \cong \Hom(F_1,F_0) \cong W$, $\Hom(F_2,F_0) \cong \Lambda^2W$
and the composition $\Hom(F_1,F_0)\otimes\Hom(F_2,F_1) \to \Hom(F_2,F_0)$
coincides with the canonical projection $W\otimes W \to \Lambda^2W$.
\end{corollary}
\begin{proof}
Follows from exact sequence~\eqref{f0} combined with the previous lemma.
\end{proof}

\begin{proposition}\label{ldcb}
The triple $(F_0^*,F_1^*,F_2^*)$ in $\TD$ is exceptional.
Moreover, the collection of subcategories~\eqref{ldtd}
is a dual Lefschetz collection in $\TD$.
\end{proposition}
\begin{proof}
The first claim follows immediately from the previous lemma by duality.
Moreover, since $\TD$ is a $\BP$-linear subcategory in $\D^b(\TY)$ by theorem~\ref{NCR},
it is stable under $\CO_\TY(-tH_Y)$-twists. Therefore, the whole collection~\eqref{ldtd}
is contained in $\TD$. So, it remains to check semiorthogonality of components of~\eqref{ldtd}.

To check that~\eqref{ldtd} is a Lefschetz collection in case $n=6$ we should check that
\begin{equation}\label{fkl}
\hspace{-10pt}
\begin{array}{lll}
\Hom^\bullet(F_l^*,F_k^*(-tH_Y)) =
%\Hom^\bullet(F_k,F_l(-tH_Y)) =
0, &
\text{for $k,l \in \{0,1,2\}$}, &  1\le t\le 8,
\quad\text{and}\quad\\[2pt]
\Hom^\bullet(F_l^*,F_2^*(-tH_Y)) =
%\Hom^\bullet(F_2,F_l(-tH_Y)) =
0, &
\text{for $l \in \{0,1,2\}$}, & 9\le t\le 11.
\end{array}
\end{equation}
Note also, that by Theorem~\ref{NCR} the Serre functor of $\TD$ acts on $F_2^*$ as
$\FS_\TD(F_2^*) \cong F_2^*(-12H_Y)[13]$, therefore
$$
\Hom^\bullet(F_l^*,F_2^*(-tH_Y)) =
\Hom^\bullet(F_2^*((12-t)H_Y),F_l^*[13])^* =
\Hom^\bullet(F_2^*,F_l^*((t-12)H_Y)[13])^*,
$$
hence the second line of~\eqref{fkl} follows from the first.

So, we must verify the first line of~\eqref{fkl}.
It will be convenient to reformulate it slightly. By duality we have
$\Hom^\bullet(F_l^*,F_k^*(-tH_Y)) = \Hom^\bullet(F_k,F_l(-tH_Y))$, so we must check that
\begin{equation}\label{fkl1}
\Hom^\bullet(F_k,F_l(-tH_Y)) = 0,
\qquad
\text{for $k,l \in \{0,1,2\}$ and $1\le t\le 8$}.
\end{equation}
Since $F_k$ is closely related to $\Lambda^{2-k}(W/\CK)$ we start by noting
$$
\Hom^\bullet(\Lambda^{2-k}(W/\CK),\Lambda^{2-l}(W/\CK)(-tH_Y)) \cong
\Hom^\bullet(\Lambda^{2-k}(W/\CK),\Lambda^{2-l}(W/\CK)\otimes \zeta_*\CO_\TY(-tH_Y)).
$$
Since $\TY = \PP_\BG(\Lambda^2\CK^\perp)$ we have
$$
\zeta_*\CO_\TY(-tH_Y) \cong
\begin{cases}
0, & \text{for $t = 1,2,3,4,5$}\\
\CO_\BG(-3H_G)[-5], & \text{for $t = 6$}\\
\Lambda^2\CK^\perp\otimes\CO_\BG(-3H_G)[-5], & \text{for $t = 7$}\\
\Sigma^{2,2}\CK^\perp\otimes\CO_\BG(-3H_G)[-5] \oplus \CO(-4H_G)[-5], & \text{for $t = 8$}
\end{cases}
$$
where we use the notation introduced in subsection~\ref{ss_bbw}.
It follows from the Borel--Bott--Weil Theorem (theorem~\ref{bbw}) that
$$
\Hom^\bullet(\Lambda^{2-k}(W/\CK),\Lambda^{2-l}(W/\CK)(-tH_Y)) = 0
$$
for $k,l \in \{0,1,2\}$ and $1\le t\le 8$.
In particular, we have~\eqref{fkl1} for $k,l \in \{1,2\}$.

Further, we have
$$
\Hom^\bullet(\CO_\TY(H_G-H_Y),\Lambda^{2-l}(W/\CK)(-tH_Y)) =
\Hom^\bullet(\CO(H_G),\Lambda^{2-l}(W/\CK)((1-t)H_Y))
$$
and the same arguments as above show that this is zero for $l \in \{0,1,2\}$ and $1\le t\le 8$.
Using~\eqref{f0} we see that~\eqref{fkl1} is satisfied also for $k = 0$ and $l \in \{1,2\}$.

Further, twisting~\eqref{f0} by $\CO_\TY(-tH_Y)$ and pushing forward to $\BG$ we compute
$$
\zeta_*F_0(-tH_Y) = \begin{cases}
0, & \text{for $t = 1,2,3,4$ and $t = 6$}\\
\CO_\TY(-2H_G)[-4], & \text{for $t = 5$}\\
%0, & \text{for $t = 6$}\\
\Sigma^{2,1,1}\CK^\perp(-2H_G)[-5], & \text{for $t=7$}\\
\Sigma^{3,2,1}\CK^\perp(-2H_G)[-5] \oplus \Lambda^2\CK^\perp(-3H_G)[-5], & \text{for $t=8$}
%\Sigma^{4,3,1}\CK^\perp(-2H_G)[-5] \oplus \Sigma^{2,2}\CK^\perp(-3H_G)[-5] \oplus \Sigma^{2,1,1}\CK^\perp(-3H_G)[-5], & \text{for $t=9$}\\
%\Sigma^{5,4,1}\CK^\perp(-2H_G)[-5] \oplus \Sigma^{3,3}\CK^\perp(-3H_G)[-5] \oplus \Sigma^{3,2,1}\CK^\perp(-3H_G)[-5] \\
%\hphantom{\Sigma^{5,4,1}\CK^\perp(-2H_G)[-5]} \oplus \Lambda^2\CK^\perp(-4H_G)[-5], & \text{for $t=10$}\\
%\Sigma^{6,5,1}\CK^\perp(-2H_G)[-5] \oplus \Sigma^{4,4}\CK^\perp(-3H_G)[-5] \oplus \Sigma^{4,3,1}\CK^\perp(-3H_G)[-5] \\
%\hphantom{\Sigma^{6,5,1}\CK^\perp(-2H_G)[-5]} \oplus \Sigma^{2,2}\CK^\perp(-4H_G)[-5] \oplus \Sigma^{2,1,1}\CK^\perp(-4H_G)[-5], & \text{for $t=11$}
\end{cases}
$$
Using the Borel--Bott-Weil Theorem again we see that
$$
\Hom^\bullet(\Lambda^{2-k}(W/\CK),F_0(-tH_Y)) = 0
$$
for $k \in \{0,1,2\}$ and $1\le t\le 8$. In particular we have~\eqref{fkl1} for $k \in \{1,2\}$ and $l = 0$.

Finally,
$$
\Hom^\bullet(\CO_\TY(H_G-H_Y),F_0(-tH_Y)) =
\Hom^\bullet(\CO(H_G),F_0((1-t)H_Y))
$$
and the same arguments as above show that this is zero for $1\le t\le 8$.
Using~\eqref{f0} we see at last that~\eqref{fkl1} is satisfied for $k = l = 0$.

Similarly, in case $n=7$ we must check that
\begin{equation}\label{fkl7}
\hspace{-10pt}
\Hom^\bullet(F_l^*,F_k^*(-tH_Y)) = 0
%\Hom^\bullet(F_k,F_l(-tH_Y)) = 0
\quad\text{for $k,l \in \{0,1,2\}$ and $1\le t\le 13$.}
\end{equation}
For $1\le t\le 7$ the same arguments as in the case $n = 6$ prove~\eqref{fkl7}.
On the other hand, for $8\le t\le 13$ we can use the Serre functor of $\TD$,
which by theorem~\ref{NCR} takes $F_l^*$ to $F_l^*(-14H_Y)[17]$. So we have
$$
\Hom^\bullet(F_l^*,F_k^*(-tH_Y)) =
\Hom^\bullet(F_k^*(-tH_Y),F_l^*(-14H_Y)[17])^* =
\Hom^\bullet(F_k^*,F_l^*((t-14)H_Y))^*[-17]
$$
and since $1\le 14-t\le 6$ for $8\le t\le 13$, we conclude that~\eqref{fkl7} holds for $8\le t\le 13$ as well.
\end{proof}

Consider the following natural complexes on $\TY$:
\begin{equation}\label{fprime}
F_2' := F_2,
\qquad
F_1' := \{W\otimes F_2 \to F_1\},
\qquad
F_0' := \{S^2W\otimes F_2 \to W\otimes F_1 \to F_0\}.
\end{equation}
Note that $F'_k \in \TD_\opp$ since $\TD_\opp$ is a triangulated subcategory of $\D^b(\TY)$.

\begin{remark}
In terms of~\cite{B} the triple $(F'_0,F'_1,F'_2)$ is the left mutation of the triple $(F_2,F_1,F_0)$.
\end{remark}

\begin{lemma}\label{ffst}
We have $F_0,F_1,F_2,F'_2,F'_1 \in \TD_\opp^{\le 0}$ and $F'_0 \in \TD_\opp^{\le 1}$.
\end{lemma}
\begin{proof}
To check that an object $F\in\TD_\opp$ is in $\TD_\opp^{\le t}$ it suffices to show that
$H^{> t}(Y,\rho^\opp_*(F)\otimes\CO_Y(nH_Y)) = 0$ for $n \gg 0$. But
\begin{multline*}
H^\bullet(Y,\rho^\opp_*(F)\otimes\CO_Y(nH_Y)) =
H^\bullet(Y,g_{Y*}(F\otimes(\CO_\TY \oplus \CK^*))(nH_Y)) = \\ =
H^\bullet(\TY,F(nH_Y)\otimes(\CO_\TY \oplus \CK^*)) =
H^\bullet(\BG,\zeta_*(F(nH_Y))\otimes(\CO_{\BG} \oplus \CK^*)),
\end{multline*}
so we should investigate the cohomology of $\zeta_*(F(nH_Y))\otimes(\CO_{\BG} \oplus \CK^*)$ on $\BG = \Gr(n-4,W)$.
Recall that we have $F_2 = \CO_\TY$, $F_1 = W/\CK$, $0 \to \CO_\TY(H_G - H_Y) \to \Lambda^2(W/\CK) \to F_0 \to 0$,
whereof we deduce $F'_2 = \CO_\TY$, $F'_1 = \CK$, and $\CO_\TY(H_G-H_Y)[1] \to F_0 \to S^2\CK$.
Further, we have
$$
\begin{array}{lll}
\zeta_*(\CO_\TY(nH_Y)) &=& S^n\Lambda^2(W/\CK),\\[2pt]
\zeta_*(\CK(nH_Y)) &=& \CK\otimes S^n\Lambda^2(W/\CK),\\[2pt]
\zeta_*(W/\CK(nH_Y)) &=& W/\CK\otimes S^n\Lambda^2(W/\CK),\\[2pt]
\zeta_*(S^2\CK(nH_Y)) &=& S^2\CK\otimes S^n\Lambda^2(W/\CK),\\[2pt]
\zeta_*(\Lambda^2(W/\CK)(nH_Y)) &=& \Lambda^2(W/\CK)\otimes S^n\Lambda^2(W/\CK),\\[2pt]
\zeta_*(\CO_\TY(H_G+(n-1)H_Y)) &=& S^{n-1}\Lambda^2(W/\CK)\otimes\CO_\BG(H_G)
\end{array}
$$
Applying Borel--Bott-Weil theorem we deduce all the claims.
\end{proof}

\section{Some computations}\label{s_comp}

In this section we compute $\rho_*(F_k^*)$ and $\rho_*({F'_k}^*)$,
where $F_k$ and $F'_k$ were defined in~\eqref{f012def} and in~\eqref{fprime}.
%Below we will need to know the pushforwards to $\BP = \PP(\Lambda^2W^*)$ of some vector bundles on $\TY$.
%We compute them in this section.

Consider the projectivization $\PP_\TY(\CK)$.
It is clear that $\PP_\TY(\CK) = \PP_{\Fl(1,n-4;W)}(\Lambda^2\CK^\perp)$.
Let us denote the line bundle $\CO_{\PP(W)}(-1)$ by $\CK_1$ (as well as all its pullbacks).
Then we have an embedding $\CK^\perp \subset \CK_1^\perp$ inducing a map
$\phi:\PP_{\Fl(1,n-4;W)}(\Lambda^2\CK^\perp) \to \PP_{\PP(W)}(\Lambda^2\CK_1^\perp)$.
Thus we have the following commutative diagram
\begin{equation}\label{diagr}
\vcenter{\xymatrix{
& \PP_\TY(\CK) \ar[d]_\phi \ar[dr]^q \ar[dl]_p \\
\Fl(1,n-4;W) \ar[d] & \PP_{\PP(W)}(\Lambda^2\CK_1^\perp) \ar[dl]_p \ar[dr]^q &
\TY \ar[d]^g \\
\PP(W) &&
\qquad\BP\qquad
}}
\end{equation}

\begin{lemma}\label{kosres}
The map $\PP_{\PP(W)}(\Lambda^2\CK_1^\perp) \to \PP(W)\times\BP$ induced by the projections
$p$ and $q$ is a closed embedding, and its image is the zero locus of a regular section
of the vector bundle $\CK_1^\perp(H_W) \boxtimes \CO_\BP(H_P)$.
In particular, we have the following Koszul resolution
$$
\dots \!\to\!
\Lambda^2(W/\CK_1)(-2H_W)\boxtimes\CO_\BP(-2H_P) \!\to\!
(W/\CK_1)(-H_W)\boxtimes\CO_\BP(-H_P) \!\to\!
\CO_{\PP(W)}\boxtimes\CO_\BP \!\to\!
\CO_{\PP_{\PP(W)}(\Lambda^2\CK_1^\perp)} \!\to\! 0,
$$
where $H_W$ is the divisor class of a hyperplane in $\PP(W)$, so that
$\CO_{\PP(W)}(-H_W) = \CK_1$.
\end{lemma}
\begin{proof}
The short exact sequence $0 \to \CK_1^\perp \to W^*\otimes\CO_{\PP(W)} \to \CO_{\PP(W)}(H_W) \to 0$
gives (by taking $\Lambda^2$) an exact sequence
$0 \to \Lambda^2\CK_1^\perp \to \Lambda^2W^*\otimes\CO_{\PP(W)} \to \CK_1^\perp(H_W) \to 0$.
Since the map $\PP_{\PP(W)}(\Lambda^2\CK_1^\perp) \to \PP(W)\times\BP$
is induced by the above embedding of the vector bundles
$\Lambda^2\CK_1^\perp \to \Lambda^2W^*\otimes\CO_{\PP(W)}$
it follows that it is a closed embedding and its image is the zero locus
of a section of $\CK_1^\perp(H_W)\boxtimes\CO_\BP(H_P)$.
Comparing the codimension of the image and the rank of the bundle
we conclude that the section is regular.
\end{proof}

\begin{lemma}\label{phisok6}
If $n = 6$ then the map $\phi:\PP_\TY(\CK) \to \PP_{\PP(W)}(\Lambda^2\CK_1^\perp)$ is birational and we have
an isomorphism
$$
\phi_*\CO_{\PP_\TY(\CK)} \cong \CO_{\PP_{\PP(W)}(\Lambda^2\CK_1^\perp)}
$$
and an exact sequence
$$
0 \to \CO_{\PP_{\PP(W)}(\Lambda^2\CK_1^\perp)}(-H_W) \to
\phi_*\CK \to
\CO_{\PP_{\PP(W)}(\Lambda^2\CK_1^\perp)}(H_W-2H_Y) \to 0.
$$
\end{lemma}
\begin{proof}
Note that $\PP_\TY(\CK) = \PP_{\Fl(1,n-4;W)}(\Lambda^2\CK^\perp)$
is the space of all triples $(K_1,K_{n-4},\omega)$, where
$K_1 \subset K_{n-4} \subset W$ is a flag of dimension $(1,n-4)$ and
$\omega$ is a skew form, such that $K_{n-4} \subset \Ker\omega$.
Similarly, $\PP_{\PP(W)}(\Lambda^2\CK_1^\perp)$
is the space of all pairs $(K_1,\omega)$,
such that $K_1 \subset \Ker\omega$. The map $\phi$
forgets $K_{n-4}$.  If $n = 6$ and $\omega$ has a nontrivial kernel
($K_1 \subset \Ker\omega$ then $r(\omega) \le 4$ hence there exists
(unique if $r(\omega) = 4$) subspace $K_2 \subset \Ker\omega$ such that
$K_1 \subset K_2$. This shows that $\phi$ is birational and proves
the isomorphism.

Further, consider
a short exact sequence $0 \to \CK_1 \to \CK \to \CK/\CK_1 \to 0$ on $\PP_\TY(\CK)$
and note that $\CK/\CK_1 \cong \CO_{\PP_\TY(\CK)}(H_W-H_G)$.
Tensoring the pullback of~\eqref{ohz} to $\PP_\TY(\CK)$ by
$\CO_{\PP_\TY(\CK)}(H_W-H_G)$, then pushing forward
to $\PP_{\PP(W)}(\Lambda^2\CK_1^\perp)$, and taking into account
that the sheaf $\CO_{\PP_\TY(\CK)}(H_W-H_G)$ is acyclic on the fibers
of $\PP_\TZ(\CK) \subset \PP_\TY(\CK)$ over $\PP_{\PP(W)}(\Lambda^2\CK_1^\perp)$
(the fibers are $\PP^2$ and the sheaf restricts to $\CO(-1)$),
we deduce that
$$
\phi_*(\CO_{\PP_\TY(\CK)}(H_W-H_G)) \cong
\phi_*(\CO_{\PP_\TY(\CK)}(H_W-2H_Y)) \cong
\CO_{\PP_{\PP(W)}(\Lambda^2\CK_1^\perp)}(H_W-2H_Y).
$$
Combining with an isomorphism
$\phi_*\CK_1 =
\phi_*\phi^*\CO_{\PP_{\PP(W)}(\Lambda^2\CK_1^\perp)}(-H_W) \cong
\CO_{\PP_{\PP(W)}(\Lambda^2\CK_1^\perp)}(-H_W)$
we obtain the required exact sequence.
\end{proof}

\begin{lemma}\label{phisok7}
If $n = 7$ then the map $\phi$ is birational onto a divisor
in $\PP_{\PP(W)}(\Lambda^2\CK_1^\perp)$ and we have
exact sequences
$$
\begin{array}{c}
0 \to
\CO_{\PP_{\PP(W)}(\Lambda^2\CK_1^\perp)}(H_W-3H_Y) \to
\CO_{\PP_{\PP(W)}(\Lambda^2\CK_1^\perp)} \to
\phi_*\CO_{\PP_\TY(\CK)} \to 0,\\[5pt]
0 \to
W\otimes\CO_{\PP_{\PP(W)}(\Lambda^2\CK_1^\perp)}(H_W-3H_Y) \to
\CM \to
\phi_*\CK \to 0,
\end{array}
$$
where $\CM$ is the unique {\rm(}non-split{\rm)} extension
$0 \to \CO_{\PP_{\PP(W)}(\Lambda^2\CK_1^\perp)}(-H_W) \to
\CM \to \CK_1^\perp(H_W - 2H_Y) \to 0$.
\end{lemma}
\begin{proof}
The image of $\phi$ is the space of all $(K_1,\omega)$ such that
$K_1 \subset \Ker\omega$ and $\omega$ is degenerate on $W/K_1$.
This actually means that the fiber of $\phi$ over $K_1 \in \PP(W)$ is the Pfaffian cubic
in $\PP(\Lambda^2(W/K_1))$. The Pfaffian is naturally an element of
$\Lambda^6\CK_1^\perp\otimes\CO_{\PP_{\PP(W)}(\Lambda^2\CK_1^\perp)}(3H_Y) \cong
\CO_{\PP_{\PP(W)}(\Lambda^2\CK_1^\perp)}(3H_Y - H_W)$, hence
$\Im \phi \sim 3H_Y - H_W$ and the first sequence follows.

Consider the map $(W/\CK_1)(-H_Y) \to \CK_1^\perp$
induced by the embedding $\CO_{\PP_{\PP(W)}(\Lambda^2\CK_1^\perp)}(-H_Y) \to \Lambda^2\CK_1^\perp$.
It is clear that its kernel is zero, and its cokernel is
$\phi_*(\CK_1^\perp/\CK^\perp) \cong \phi_*(\CK/\CK_1)^*$,
so that we have an exact sequence
$$
0 \to (W/\CK_1)(-H_Y) \to \CK_1^\perp \to \phi_*(\CK/\CK_1)^* \to 0.
$$
On the other hand, it is clear that
$$
\CK/\CK_1 \cong \det(\CK/\CK_1) \otimes (\CK/\CK_1)^* \cong
(\CK/\CK_1)^*(H_W-H_G).
$$
Tensoring the pullback of~\eqref{ohz} to $\PP_\TY(\CK)$ by
$(\CK/\CK_1)^*(H_W-H_G)$, then pushing forward
to $\PP_{\PP(W)}(\Lambda^2\CK_1^\perp)$, and taking into account
that the sheaf $(\CK/\CK_1)^*(H_W-H_G)$ is acyclic on the fibers
of $\PP_\TZ(\CK) \subset \PP_\TY(\CK)$ over $\PP_{\PP(W)}(\Lambda^2\CK_1^\perp)$
(the fibers are $\Gr(2,4)$ and the sheaf restricts to the tautological rank 2 subbundle),
we deduce that
$\phi_*(\CK/\CK_1) \cong
\phi_*((\CK/\CK_1)^*(H_W-H_G)) \cong
\phi_*((\CK/\CK_1)^*(H_W-2H_Y))$,
hence by the projection formula we have an exact sequence
$$
0 \to (W/\CK_1)(H_W-3H_Y) \to \CK_1^\perp(H_W-2H_Y) \to \phi_*(\CK/\CK_1) \to 0.
$$
On the other hand,
$\phi_*\CK_1 =
\phi_*\phi^*\CO_{\PP_{\PP(W)}(\Lambda^2\CK_1^\perp)}(-H_W) \cong
(\phi_*\CO_{\PP_\TY(\CK)})(-H_W)$,
hence we have an exact sequence
$$
0 \to
\CO_{\PP_{\PP(W)}(\Lambda^2\CK_1^\perp)}(-3H_Y) \to
\CO_{\PP_{\PP(W)}(\Lambda^2\CK_1^\perp)}(-H_W) \to
\phi_*\CK_1 \to 0.
$$
Now applying $\phi_*$ to the exact sequence
$0 \to \CK_1 \to \CK \to \CK/\CK_1 \to 0$
we deduce the second claim of the lemma.
\end{proof}

\begin{proposition}\label{gsk}
Let $n = 6$.
We have the following resolutions
$$
\arraycolsep=1.5pt
\begin{array}{rl}
0 \to \CO_{\BP}(-3) \to \CO_{\BP} \to & g_*F_2^* \to 0,\\[2pt]
0 \to W^*\otimes\CO_{\BP}(-3) \to W\otimes\CO_{\BP}(-1) \to & g_*F_1^* \to 0,\\[2pt]
0 \to \Lambda^2W^*\otimes\CO_{\BP}(-3) \to \Lambda^2W\otimes\CO_{\BP}(-2) \to & g_*F_0^* \to 0,\\[5pt]
0 \to W\otimes\CO_{\BP}(-3) \to W^*\otimes\CO_{\BP}(-2) \to & g_*(F_2^*\otimes\CK) \to 0,\\[2pt]
0 \to \CO_{\BP}(-6) \to \Lambda^2W\otimes\CO_{\BP}(-4) \to \Lambda^2W^*\otimes\CO_{\BP}(-2) \to \CO_{\BP} \to & g_*(F_1^*\otimes\CK)[1] \to 0,\\[2pt]
0 \to W^*\otimes\CO_{\BP}(-6) \to \Lambda^3W^*\otimes\CO_{\BP}(-5) \to \Lambda^3W\otimes\CO_{\BP}(-2) \to W\otimes\CO_{\BP}(-1) \to & g_*(F_0^*\otimes\CK)[1] \to 0.
\end{array}
$$
In particular, we have
$$
\eta_*\rho_*(F_0^*),\eta_*\rho_*(F_1^*) \in \D^{\le 1}(\BP), \qquad
\eta_*\rho_*(F_2^*) \in \D^{\le 0}(\BP),
$$
and for all $l = 0,1,2$ we have $g^*\eta_*\rho_*(F_l^*) \in \D^{\ge -1}(\TY)$.
\end{proposition}
\begin{proof}
Consider the diagram~\eqref{diagr}. We have
$$
g_*S^l\CK^* \cong
g_*q_*\phi^*p^*\CO_{\PP(W)}(lH_W) \cong
q_*\phi_*\phi^*p^*\CO_{\PP(W)}(lH_W) \cong
q_*(p^*\CO_{\PP(W)}(lH_W)\otimes\phi_*\CO_{\PP_\TY(\CK)})
$$
and similarly
$$
g_*(S^l\CK^*\otimes\CK) \cong
g_*q_*(\phi^*p^*\CO_{\PP(W)}(lH_W)\otimes\CK) \cong
q_*\phi_*(\phi^*p^*\CO_{\PP(W)}(lH_W)\otimes\CK) \cong
q_*(p^*\CO_{\PP(W)}(lH_W)\otimes\phi_*\CK).
$$
Using resolutions of lemma~\ref{phisok6} we reduce the computations
of $g_*S^l\CK^*$ and $g_*(S^l\CK^*\otimes\CK)$ for $0\le l\le 2$ to the computation of
$q_*(p^*\CO_{\PP(W)}(tH_W))$ for $-1\le t\le 3$.
Further, using the resolution of lemma~\ref{kosres}
we reduce these to the computation of
$H^\bullet(\PP(W),\Lambda^s(W/\CK_1)((t-s)H_W)$
for $0\le s\le 5$, $-1\le t\le 3$, which can be easily done
by the Borel--Bott--Weil Theorem.

Further, we use the following evident resolutions
$$
0 \to \CK^\perp \to W^*\otimes\CO_\TY \to \CK^* \to 0,
\qquad
0 \to \Lambda^2\CK^\perp \to \Lambda^2W^*\otimes\CO_\TY \to W^*\otimes\CK^* \to S^2\CK^* \to 0
$$
to compute $g_*(\Lambda^l\CK^\perp)$ and $g_*(\Lambda^l\CK^\perp\otimes\CK)$.
Since
$0 \to F_0^* \to \Lambda^2\CK^\perp \to \CO_\TY(H_Y - H_G) \to 0$,
$F_1^* = \CK^\perp$, and $F_2^* = \CO_\TY$,
it remains to compute $g_*(\CO_\TY(H_Y - H_G))$ and $g_*(\CK(H_Y - H_G))$.
For this we tensor~\eqref{ohz} by $\CO_\TY(H_Y - H_G)$ and $\CK(H_Y - H_G)$
respectively and note that
$$
g_*(\CO_\TY(H_Y - H_G)\otimes\CO_\TZ) \cong
\eta_*g_{Y*}i_{\TY*}\CO_\TZ(H_Y - H_Z) \cong
\eta_*i_{Y*}g_{Z*}\CO_\TZ(H_Y - H_Z)
$$
(and similarly $g_*(\CK(H_Y - H_G)\otimes\CO_\TZ) \cong \eta_*i_{Y*}g_{Z*}(\CK_{|\TZ}(H_Y - H_Z))$).
Since
$$
g_{Z*}\CO_\TZ(H_Y - H_Z) = g_{Z*}(\CK_{|\TZ}(H_Y - H_Z)) = 0
$$
(the fibers of $g_Z$ are Grassmannians $\Gr(2,4)$ and the sheaves in question
restrict to the fibers as $\CO(-1)$ and $\CK(-1)$ which are acyclic by the Borel--Bott--Weil Theorem),
we conclude that
$$
g_*(\CO_\TY(H_Y - H_G)) \cong g_*\CO_\TY(-H_Y),
\qquad
g_*(\CK(H_Y - H_G)) \cong g_*\CK(-H_Y),
$$
which we already have computed above. This way we obtain the desired resolutions for $g_*F_l^*$
and $g_*(F_l^*\otimes\CK)$.

Further, $\eta_*\rho_*(F_l^*) = g_*F_l^* \oplus g_*(F_l^*\otimes\CK)$ by definition of the functor $\rho_*$,
and the second claim follows. For the last claim it suffices to check that
the pullbacks via $g^*$ of our resolutions lie in the subcategory $\D^{\ge -1}(\TY)$.
For the first four resolutions this is evident. For the last two we have to check
that the maps $\CO_\TY(-6H_Y) \to \Lambda^2W\otimes\CO_\TY(-5H_Y)$ and
$W^*\otimes\CO_\TY(-6H_Y) \to \Lambda^3W^*\otimes\CO_\TY(-5H_Y)$ are embeddings.
This can be checked in a generic point of $\TY$ and corresponds to the fact
that for a skew-form $\omega$ of rank $4$ the maps
$\xymatrix@1{\kk \ar[r]^-{\omega\wedge\omega} & \Lambda^2W}$ and
$\xymatrix@1{W^* \ar[r]^-{-\wedge\omega} & \Lambda^3W^*}$ are embeddings, which is clear.
\end{proof}

\begin{proposition}\label{gsk7}
Let $n = 7$. We have the following resolutions
$$
\arraycolsep = 1.5pt
\begin{array}{rl}
0 \to \CO_{\BP}(-7) \to W\otimes\CO_{\BP}(-4) \to W^*\otimes\CO_{\BP}(-3) \to \CO_{\BP} \to & g_*F_2^* \to 0,\\[2pt]
0 \to W^*\otimes\CO_{\BP}(-7) \to \Lambda^2W\otimes\CO_{\BP}(-5) \to \Lambda^2W^*\otimes\CO_{\BP}(-3) \to W\otimes\CO_{\BP}(-1) \to & g_*F_1^* \to 0,\\[2pt]
0 \to \Lambda^2W^*\otimes\CO_{\BP}(-7) \to \Lambda^3W\otimes\CO_{\BP}(-6) \to \Lambda^3W^*\otimes\CO_{\BP}(-3) \to \Lambda^2W\otimes\CO_{\BP}(-2) \to & g_*F_0^* \to 0,\\[5pt]
%
%------------------------------------------------------------------------------------------------------------------------
%
0 \to W\otimes\CO_{\BP}(-7) \to S^2W\otimes\CO_{\BP}(-4) \oplus W^*\otimes\CO_{\BP}(-6) \to
\quad \\[2pt] \to
W\otimes W^*\otimes\CO_{\BP}(-3) \to \Lambda^2W^*\otimes\CO_{\BP}(-2) \to & g_*(F_2^*\otimes\CK) \to 0,\\[2pt]
%------------------------------------------------------------------------------------------------------------------------
0 \to (W^*\otimes W/\kk)\otimes\CO_{\BP}(-7) \to
\hspace{225pt} \\[2pt] \to
\CO_{\BP}(-7) \oplus W^*\otimes W^*\otimes\CO_{\BP}(-6) \oplus (W\otimes\Lambda^2W/\Lambda^3W)\otimes\CO_\BP(-5) \to
\quad \\[2pt] \to
W\otimes\CO_{\BP}(-4) \oplus W^*\otimes\Lambda^2W\otimes\CO_{\BP}(-4) \to
\quad \\[2pt] \to
W^*\otimes\CO_{\BP}(-3) \oplus \Lambda^3W^*\otimes\CO_{\BP}(-2) \to \CO_{\BP} \to & g_*(F_1^*\otimes\CK)[1] \to 0,\\[2pt]
%------------------------------------------------------------------------------------------------------------------------
0 \to (\Lambda^2W^*\otimes W/W^*)\otimes\CO_{\BP}(-7) \to
\hspace{225pt} \\[2pt] \to
W^*\otimes\CO_{\BP}(-7) \oplus (\Lambda^3W\otimes W/\Lambda^4W)\otimes\CO_{\BP}(-6) \oplus \Lambda^2W^*\otimes W^*\otimes\CO_\BP(-6) \to
\quad \\[2pt] \to
\Lambda^3W\otimes W^*\otimes\CO_{\BP}(-5) \oplus W^*\otimes W\otimes\CO_{\BP}(-4) \to
\quad \\[2pt] \to
W^*\otimes W^*\otimes\CO_{\BP}(-3) \oplus \Lambda^2W^*\otimes\CO_{\BP}(-2) \to W\otimes\CO_{\BP}(-1) \to & g_*(F_0^*\otimes\CK)[1] \to 0.
\end{array}
$$
In particular, we have
$$
\eta_*\rho_*(F_0^*),\eta_*\rho_*(F_1^*) \in \D^{\le 1}(\BP), \qquad
\eta_*\rho_*(F_2^*) \in \D^{\le 0}(\BP),
$$
and for all $l = 0,1,2$ we have $g^*\eta_*\rho_*(F_l^*) \in \D^{\ge -3}(\TY)$.
\end{proposition}
\begin{proof}
The same arguments as in the proof of proposition~\ref{gsk} with lemma~\ref{phisok6} replaced
by lemma~\ref{phisok7}.
\end{proof}

\begin{proposition}\label{gskk}
If $n = 6$ then $\CH^{-1}(\eta_*\rho_*(F_2^*) \otimes \eta_*\rho^\opp_*(F_2)) \cong \eta_*\CR(-3)$.
Similarly, if $n = 7$ then $\CH^{-3}(\eta_*\rho_*(F_2^*) \otimes \eta_*\rho^\opp_*(F_2)) \cong \eta_*\CR(-7)$.
\end{proposition}
\begin{proof}
Note that by definition of the functors $\rho_*$ and $\rho^\opp_*$ we have
$$
\eta_*\rho_*(F_2^*) \otimes \eta_*\rho^\opp_*(F_2) \cong
g_*(\CO_\TY\oplus\CK) \otimes g_*(\CO_\TY\oplus\CK^*) \cong
g_*(g^*g_*(\CO_\TY\oplus\CK) \otimes (\CO_\TY\oplus\CK^*)).
$$
Using resolutions of propositions~\ref{gsk} and \ref{gsk7} we deduce that
$g^*g_*(\CO_\TY\oplus\CK) \in \D^{\ge -1}$, if $n = 6$, and
$g^*g_*(\CO_\TY\oplus\CK) \in \D^{\ge -3}$, if $n = 7$,
%$$
%\begin{array}{ll}
%g^*g_*(\CO_\TY\oplus\CK) \in \D^{\ge -1}, & \text{if $n = 6$}\\[2pt]
%g^*g_*(\CO_\TY\oplus\CK) \in \D^{\ge -3}, & \text{if $n = 7$}
%\end{array}
%$$
and that
$$
\begin{array}{ll}
\CH^{-1}(g^*g_*(\CO_\TY\oplus\CK)) \cong \CO_\TY(-3H_Y)\oplus\CK(-3H_Y), & \text{if $n = 6$}\\[2pt]
\CH^{-3}(g^*g_*(\CO_\TY\oplus\CK)) \cong \CO_\TY(-7H_Y)\oplus\CK(-7H_Y), & \text{if $n = 7$.}
\end{array}
$$
Since the functor $g_*$ is left-exact and
$g_*((\CO_\TY\oplus\CK) \otimes (\CO_\TY\oplus\CK^*)) \cong \eta_*\CR$
is a pure sheaf by theorem~\ref{NCR}, we deduce the desired isomorphisms.
\end{proof}

\begin{proposition}\label{gskkp}
If $n = 6$ then
$$
\eta_*\rho_*({F_2'}^*) \otimes \eta_*\rho^\opp_*(F_2) \in \D^{\ge -1\!}(\BP),\
\eta_*\rho_*({F_1'}^*) \otimes \eta_*\rho^\opp_*(F_1) \in \D^{\ge -2\!}(\BP),\
\eta_*\rho_*({F_0'}^*) \otimes \eta_*\rho^\opp_*(F_0) \in \D^{\ge -3\!}(\BP),
$$
and if $n = 7$ then
$$
\eta_*\rho_*({F_2'}^*) \otimes \eta_*\rho^\opp_*(F_2) \in \D^{\ge -3\!}(\BP),\
\eta_*\rho_*({F_1'}^*) \otimes \eta_*\rho^\opp_*(F_1) \in \D^{\ge -3\!}(\BP),\
\eta_*\rho_*({F_0'}^*) \otimes \eta_*\rho^\opp_*(F_0) \in \D^{\ge -4\!}(\BP).
$$
Moreover, if $n = 7$ then
%$\eta_*\rho_*({F_1'}^*) \otimes \eta_*\rho_*(F_1) \in \D^{\ge -3}(\BP)$,
%$\eta_*\rho_*({F_0'}^*) \otimes \eta_*\rho_*(F_0) \in \D^{\ge -4}(\BP)$.
%Moreover, the sheaves
$\CH^{-3}(\eta_*\rho_*({F_1'}^*) \otimes \eta_*\rho_*(F_1)) \cong \eta_*\FF_1$
and
$\CH^{-4}(\eta_*\rho_*({F_0'}^*) \otimes \eta_*\rho_*(F_0)) \cong \eta_*\FF_0$
where $\FF_1,\FF_0$ are torsion free sheaves on $Y$ and
$$
\begin{array}{l}
{\FF_1}_{|Y\setminus Z} \cong (\CK^*\otimes F_1\otimes\CEnd(\CO_\TY\oplus\CK)\otimes\CO_\TY(-7H_Y))_{|\TY\setminus\TZ},\\[2pt]
{\FF_0}_{|Y\setminus Z} \cong (\CO_\TY(H_Y-H_G)\otimes F_0\otimes\CEnd(\CO_\TY\oplus\CK)\otimes\CO_\TY(-7H_Y))_{|\TY\setminus\TZ}
\end{array}
$$
{\rm(}recall that the map $g_Y:\TY \to Y$ identifies $\TY \setminus \TZ$ with $Y\setminus Z${\rm)}.
\end{proposition}
\begin{proof}
%Note that
%$$
%\eta_*\rho_*({F_k'}^*) \otimes \eta_*\rho_*(F_k) \cong
%g_*({F_k'}^*\oplus {F_k'}^*\otimes\CK) \otimes g_*(F_k\oplus F_k\otimes\CK^*) \cong
%g_*(g^*g_*({F_k'}^*\oplus {F_k'}^*\otimes\CK) \otimes (F_k\oplus F_k\otimes\CK^*)).
%$$
Using resolutions of propositions~\ref{gsk} and~\ref{gsk7} and quasiisomorphisms
$$
F_2^* \cong {F_2'}^*,
\qquad
\{ F_1^* \to W^*\otimes F_2^* \} \cong {F_1'}^*,
\qquad
\{ F_0^* \to W^*\otimes F_1^* \to S^2W^*\otimes F_2^* \} \cong {F_0'}^*,
$$
it is easy to see that $\eta_*\rho_*({F_2'}^*)$, $\eta_*\rho_*({F_1'}^*)$ and $\eta_*\rho_*({F_0'}^*)$
have locally free resolutions concentrated in degrees
$-1 \dots 0$, $-2 \dots 0$ and $-3 \dots 0$ respectively if $n = 6$
and in degrees
$-3 \dots 0$, $-3 \dots 0$ and $-4 \dots 0$ respectively if $n = 7$.
Moreover, since $F_k$ is a vector bundle on $\TY$, it follows that
$\eta_*\rho^\opp_*(F_k) \in \D^{\ge 0}(\BP)$, and the first claims follow.
Further, $\CH^0(\rho^\opp_*(F_k)) = \CH^0(g_{Y*}(F_k\otimes(\CO_\TY \oplus \CK^*)))$
is a torsion free sheaf on~$Y$, hence the bottom cohomology of
$\eta_*\rho_*({F_k'}^*)\otimes\eta_*\rho^\opp_*(F_k)$ is a torsion free sheaf on~$Y$ as well.
Finally, over $Y\setminus Z$ the map $g_Y$ is an isomorphism, hence
$$
\rho_*({F'_k}^*)_{|Y\setminus Z} \cong ({F'_k}^*\otimes(\CO_\TY\oplus\CK))_{|\TY \setminus \TZ},
\qquad
\rho^\opp_*(F_k)_{|Y\setminus Z} \cong ({F_k}\otimes(\CO_\TY\oplus\CK^*))_{|\TY \setminus \TZ}.
$$
Taking into account that ${F'_1}^* \cong \CK^*$ is a vector bundle on $\TY$
and that ${F'_0}^*$ has two cohomology (in the standard t-structure),
$S^2\CK^*$ the 0-th and $\CO_\TY(H_Y-H_G)$ the $(-1)$-st,
we see that
$$
\begin{array}{l}
\eta_*{\FF_1}_{|Y \setminus Z} \cong
\Tor_3^{\CO_\BP}(\eta_*((\CK^*\otimes(\CO_\TY\oplus\CK))_{|\TY \setminus \TZ}),\eta_*(({F_1}\otimes(\CO_\TY\oplus\CK^*))_{|\TY \setminus \TZ})),
\\[2pt]
\eta_*{\FF_0}_{|Y \setminus Z} \cong
\Tor_3^{\CO_\BP}(\eta_*((\CO_\TY(H_Y-H_G)\otimes(\CO_\TY\oplus\CK))_{|\TY \setminus \TZ}),\eta_*(({F_0}\otimes(\CO_\TY\oplus\CK^*))_{|\TY \setminus \TZ})),
\end{array}
$$
and since the embedding $\eta:Y\setminus Z \to \BP\setminus Z$ is a regular embedding
of codimension 3, and moreover $\det\CN^*_{(Y\setminus Z)/(\BP\setminus Z)} \cong \CO_{Y\setminus Z}(-7H_P)$
(this follows, e.g., from proposition~\ref{gsk7}), we deduce the desired isomorphisms.
\end{proof}

\section{Further computations}\label{s_fd}

Let $\beta:\TY\to\BP\times\TY$ be the graph of the projection $g:\TY\to\BP$.
Denote by $\beta':\TY \to \TY\times\BP$ the composition of $\beta$
with the transposition $\BP\times \TY \to \TY\times\BP$. Consider the maps
$\id_\TY\times\beta,\beta'\times\id_\TY:\TY\times \TY \to \TY\times\BP\times \TY$
and denote them by $\beta$ and $\beta'$ for brevity. Let
\begin{equation}\label{fd}
\FD := \beta'_*\CO_{\TY\times \TY} \otimes \beta_*\CO_{\TY\times \TY}
\in \D^b(\TY\times\BP\times \TY).
\end{equation}
%Let also $\eta:Y \to \BP$ be the natural embedding and
%$\Delta^\BP:\BP \to \BP\times\BP$, $\Delta^Y:Y \to Y\times Y$ be the diagonals.

\begin{lemma}\label{fd1}
We have $\FD \cong (g\times\id_\BP\times g)^*\Delta^\BP_*\CO_\BP$,
where $(g\times\id_\BP\times g):\TY\times\BP\times\TY \to \BP\times\BP\times\BP$
is the projection and $\Delta^\BP:\BP \to \BP\times\BP\times\BP$ is the diagonal.
In particular, $\FD$ is supported on a infinitesimal neighborhood of
$\TY\times_\BP \TY = \TY\times_\BP \BP\times_\BP \TY \subset \TY\times\BP\times \TY$.
\end{lemma}
\begin{proof}
Consider the following diagram
$$
\xymatrix{
\TY \times\TY \ar[rr]^-{\beta'} \ar[d]^{g\times g} &&
\TY \times\BP \times \TY       \ar[d]_{g\times \id_\BP\times g} &&
\TY \times\TY \ar[ll]_-{\beta}  \ar[d]^{g\times g} \\
\BP \times \BP \ar[rr]^-{\id_\BP\times\Delta^\BP} &&
\BP \times \BP \times \BP &&
\BP \times \BP \ar[ll]_-{\Delta^\BP\times\id_\BP}
}
$$
Note that we have base change isomorphisms
$\beta_*(g\times g)^* \cong (g\times\id_\BP\times g)^*(\Delta^\BP\times\id_\BP)_*$ and
$\beta'_*(g\times g)^* \cong (g\times\id_\BP\times g)^*(\id_\BP\times\Delta^\BP)_*$
by Corollary~2.27 of~\cite{K1}. Therefore
\begin{multline*}
\FD = \beta'_*\CO_{\TY\times \TY} \otimes \beta_*\CO_{\TY\times \TY} \cong
\beta'_*(g\times g)^*\CO_{\BP\times \BP} \otimes \beta_*(g\times g)^*\CO_{\BP\times \BP} \cong
\\ \cong
(g\times\id_\BP\times g)^*(\id_\BP\times\Delta^\BP)_*\CO_{\BP\times \BP} \otimes
(g\times\id_\BP\times g)^*(\Delta^\BP\times\id_\BP)_*\CO_{\BP\times \BP} \cong
\\ \cong
(g\times\id_\BP\times g)^*((\id_\BP\times\Delta^\BP)_*\CO_{\BP\times \BP} \otimes
(\Delta^\BP\times\id_\BP)_*\CO_{\BP\times \BP}) \cong
(g\times\id_\BP\times g)^*\Delta^\BP_*\CO_\BP,
\end{multline*}
which gives us the desired isomorphism.
\end{proof}

Recall that in section~\ref{ncrpf} we have defined  a triangulated subcategory
$\TD_{\natural\BP} \subset \D^b(\TY\times\BP\times\TY)$ and introduced on it a t-structure.

\begin{lemma}\label{fkflfd}
For any $t\in\ZZ$ and $k,l\in\{0,1,2\}$ we have
$$
\begin{array}{ll}
(F_k^*\boxtimes\CO_\BP(tH_P)\boxtimes F_l)\otimes \FD \in \TD_{\natural\BP}^{\ge -1} \subset \D^b(\TY\times\BP\times \TY),
& \text{if $n=6$, and}\\[5pt]
(F_k^*\boxtimes\CO_\BP(tH_P)\boxtimes F_l)\otimes \FD \in \TD_{\natural\BP}^{\ge -3} \subset \D^b(\TY\times\BP\times \TY),
& \text{if $n=7$.}
\end{array}
$$
\end{lemma}
\begin{proof}
First of all note that
$(F_k^*\boxtimes\CO_\BP(tH_P)\boxtimes F_l)\otimes \FD \in \TD_{\natural\BP}$
by lemma~\ref{fd1} since $\TD_{\natural\BP}$ is a $(\BP\times\BP\times\BP)$-linear
subcategory in $\D^b(\TY\times\BP\times\TY)$.
Further, again by lemma~\ref{fd1} we have
\begin{multline*}
(\eta\times\id_\BP\times\eta)_*\rho^{\natural\BP}_*((F_k^*\boxtimes\CO_\BP(tH_P)\boxtimes F_l)\otimes \FD) \cong
\\ \cong
(g\times\id_\BP\times g)_*((F_k^*\otimes(\CO_\TY\oplus\CK))\boxtimes\CO_\BP(tH_P)\boxtimes(F_l\otimes(\CO_\TY\oplus\CK^*)))\otimes\Delta^\BP_*\CO_\BP \cong
\qquad\qquad \\ \qquad\qquad \cong
(g_*(F_k^*\otimes(\CO_\TY\oplus\CK))\boxtimes\CO_\BP(tH_P)\boxtimes g_*(F_l\otimes(\CO_\TY\oplus\CK^*)))\otimes\Delta^\BP_*\CO_\BP \cong
\\ \qquad\qquad\qquad\qquad\qquad \cong
\Delta^\BP_*(g_*(F_k^*\otimes(\CO_\TY\oplus\CK))\otimes\CO_\BP(tH_P)\otimes g_*(F_l\otimes(\CO_\TY\oplus\CK^*))).
%\cong
%\\ \cong
%\Delta^\BP_*(\eta_*\rho_*(F_k^*) \otimes \eta_*\rho^\opp_*(F_l) \otimes \CO_\BP(tH_P)).
\end{multline*}
Further by projection formula we have
$$
g_*(F_k^*\otimes(\CO_\TY\oplus\CK))\otimes g_*(F_l\otimes(\CO_\TY\oplus\CK^*)) \cong
g_*(g^*g_*(F_k^*\otimes(\CO_\TY\oplus\CK))\otimes (F_l\otimes(\CO_\TY\oplus\CK^*)).
$$
Note that by propositions~\ref{gsk} and \ref{gsk7} we have
$$
\begin{array}{ll}
g^*g_*(F_k^*\otimes(\CO_\TY\oplus\CK)) = g^*\eta_*\rho_*(F_k^*) \in \D^{\ge -1}(\TY),
& \text{if $n=6$, and}\\[2pt]
g^*g_*(F_k^*\otimes(\CO_\TY\oplus\CK)) = g^*\eta_*\rho_*(F_k^*) \in \D^{\ge -3}(\TY),
& \text{if $n=7$.}
\end{array}
$$
Since $F_l\otimes(\CO_\TY\oplus\CK^*)$ is a vector bundle on $\TY$ and the pushforward functor $g_*$ is right-exact
we see that $(\eta\times\id_\BP\times\eta)_*\rho^{\natural\BP}_*((F_k^*\boxtimes\CO_\BP(tH_P)\boxtimes F_l)\otimes \FD)$
is in $\D^{\ge -1}(\BP)$ for $n = 6$ and in $\D^{\ge -3}(\BP)$ for $n = 7$.
Hence by definition of t-structure on $\TD_{\natural\BP}$ the claim follows.
\end{proof}

Let $\pi:\TY\times\BP\times\TY \to \TY\times\TY$ be the projection.

\begin{lemma}\label{h0fd}
If $n = 6$ then $\TCH^{-1}((F_2^*\boxtimes F_2)\otimes\pi_*\FD) \cong \rho_\natural^*(\Delta^Y_*\CR(-3H_Y))$.
\end{lemma}
\begin{proof}
Consider the following diagram
$$
\xymatrix{
\TY \times\BP \times \TY \ar[rr]^-{\pi} \ar[d]_{g\times \id_\BP\times g} &&
\TY\times\TY \ar[d]_{g\times g} \\
\BP \times \BP \times \BP \ar[rr]^-{p_{13}} &&
\BP\times\BP
}
$$
Note that we have a base change isomorphism
$\pi_*(g\times\id_\BP\times g)^* \cong (g\times g)^*p_{13*}$
by Corollary~2.27 of~\cite{K1}, hence
$$
\pi_*\FD \cong
\pi_*(g\times\id_\BP\times g)^*\Delta^\BP_*\CO_\BP \cong
(g\times g)^*p_{13*}\Delta^\BP_*\CO_\BP \cong
(g\times g)^*\Delta^\BP_*\CO_\BP,
$$
where $\Delta^\BP$ in the RHS stands for the diagonal $\BP \to \BP\times\BP$.
Further,
\begin{multline*}
(\eta\times\eta)_*\rho^\natural_*((F_2^*\boxtimes F_2)\otimes\pi_*\FD) \cong
%\\ \cong
(g\times g)_*([(\CO_\TY\oplus\CK)\boxtimes(\CO_\TY\oplus\CK^*)]\otimes(g\times g)^*\Delta^\BP_*\CO_\BP) \cong
\\ \cong
[g_*(\CO_\TY\oplus\CK)\boxtimes g_*(\CO_\TY\oplus\CK^*)]\otimes\Delta^\BP_*\CO_\BP \cong
%\Delta^\BP_*(\Delta^{\BP*}[g_*(\CO_\TY\oplus\CK)\boxtimes g_*(\CO_\TY\oplus\CK^*)]) \cong
\Delta^\BP_*(g_*(\CO_\TY\oplus\CK)\otimes g_*(\CO_\TY\oplus\CK^*))
\end{multline*}
and we conclude by proposition~\ref{gskk}
since $g_*(\CO_\TY\oplus\CK)\otimes g_*(\CO_\TY\oplus\CK^*) \cong \eta_*\rho_*(F_2^*)\otimes\eta_*\rho^\opp_*(F_2)$.
\end{proof}

Finally, consider on $\TY\times \TY$ the objects $\CT$ and $\CT^*\in\D^b(\TY\times \TY)$
defined as the following complexes
\begin{equation}\label{tts1}
\arraycolsep = 1pt
\left\{
\begin{array}{ccccc}
&&&& F_0^*\boxtimes F_0\\
&& W\otimes F_0^*\boxtimes F_1 && \oplus \\
S^2W\otimes F_0^*\boxtimes F_2 & \to & \oplus & \to & F_1^*\boxtimes F_1 \\
&& W\otimes F_1^*\boxtimes F_2 && \oplus \\
&&&& F_2^*\boxtimes F_2
\end{array}
\right\} \cong \CT^*
\end{equation}
\begin{equation}\label{tts2}
\left\{
\begin{array}{ccccc}
F_0^*\boxtimes F_0 \\
\oplus && W^*\otimes F_1^*\boxtimes F_0 \\
F_1^*\boxtimes F_1 & \to & \oplus & \to & S^2W^*\otimes F_2^*\boxtimes F_0\\
\oplus && W^*\otimes F_2^*\boxtimes F_1 \\
F_2^*\boxtimes F_2
\end{array}
\right\} \cong \CT
\end{equation}

\begin{lemma}\label{tts}
We have $\pi^*\CT^*\otimes\FD \in \TD_{\natural\BP}^{\le 0}$ for $n = 6,7$.
Further, if $n = 6$ we have
$\pi^*\CT\otimes\FD \in \TD_{\natural\BP}^{\ge -3}$,
and if $n = 7$ we have
$\pi^*\CT\otimes\FD \in \TD_{\natural\BP}^{\ge -5}$.
%
%
%
%$$
%\pi^*\CT^*\otimes\FD \in \TD_{\natural\BP}^{\le 0},
%\quad\text{and}\quad
%\begin{array}{ll}
%\pi^*\CT\otimes\FD \in \TD_{\natural\BP}^{\ge -3}, & \text{if $n = 6$,}\\[2pt]
%\pi^*\CT\otimes\FD \in \TD_{\natural\BP}^{\ge -5}, & \text{if $n = 7$.}
%\end{array}
%$$
Moreover, if $n = 7$ then $\TCH^{-5}(\CT\otimes\pi_*(\FD)) \cong \rho_\natural^*(\Delta^Y_*\CR(-7H_Y))$,
where $\Delta^Y:Y \to Y\times Y$ is the diagonal.
\end{lemma}
\begin{proof}
First of all note, that using the definition of objects $F'_2$, $F'_1$ and $F'_0$
we can rewrite $\CT^*$ in the form
$$
\CT^* \cong \{F_0^*\boxtimes F_0' \to F_1^*\boxtimes F_1' \to F_2^*\boxtimes F_2\}.
$$
Arguing like in the proof of lemma~\ref{fkflfd} we can show that
$$
(\eta\times\id_\BP\times\eta)_*\rho^{\natural\BP}_*((F_k^*\boxtimes\CO_\BP\boxtimes F'_k)\otimes \FD) \cong
\Delta^\BP_*(\eta_*\rho_*(F_k^*) \otimes \eta_*\rho^\opp_*(F'_k)).
$$
Further, by lemmas~\ref{gsk} and \ref{gsk7} we have
$\eta_*\rho_*(F_0^*),\eta_*\rho_*(F_1^*) \in \D^{\le 1}(\BP)$, $\eta_*\rho_*(F_2^*) \in \D^{\le 0}(\BP)$.
On the other hand, by lemma~\ref{ffst} we have
$\eta_*\rho^\opp_*(F'_0) \in \D^{\le 1}(\BP)$, $\eta_*\rho^\opp_*(F'_1),\eta_*\rho^\opp_*(F'_2) \in \D^{\le 0}(\BP)$.
Combining, we deduce the first claim.

Similarly, $\CT$ can be rewritten as
$$
\CT \cong \{{F_2'}^*\boxtimes F_2 \to {F'_1}^*\boxtimes F_1 \to {F'_0}^*\boxtimes F_0\}.
$$
Again, we have
$
(\eta\times\id_\BP\times\eta)_*\rho^{\natural\BP}_*(({F_k'}^*\boxtimes\CO_\BP\boxtimes F_k)\otimes \FD) \cong
\Delta^\BP_*(\eta_*\rho_*({F_k'}^*) \otimes \eta_*\rho^\opp_*(F_k))
$
and we deduce the second claim from proposition~\ref{gskkp}.

Finally, put $n = 7$ and consider the spectral sequence
$$
\BE_1^{-p,q} = \TCH^{q}(({F_p'}^*\boxtimes F_p)\otimes \pi_*(\FD)) \Longrightarrow \BE_\infty^n = \TCH^{n}(\CT\otimes\pi_*(\FD)).
$$
Since
$(\eta\times\eta)_*\rho^\natural_*(({F_k'}^*\boxtimes F_k)\otimes \pi_*(\FD)) \cong
\Delta^\BP_*(\eta_*\rho_*({F_k'}^*) \otimes \eta_*\rho^\opp_*(F_k))$
it follows from proposition~\ref{gskkp} that $\BE_\infty^{-5} = \BE_3^{-2,-3}$ and
$$
\begin{array}{lllll}
\BE_2^{-2,-3} &=& \Ker (\TCH^{-3}(({F_2'}^*\boxtimes F_2)\otimes \pi_*(\FD)) & \stackrel{d_1}\to & \TCH^{-3}(({F_1'}^*\boxtimes F_1)\otimes \pi_*(\FD))),\\[2pt]
\BE_3^{-2,-3} &=& \Ker (\qquad\qquad \BE_2^{-2,-3} & \stackrel{d_2}\to & \TCH^{-4}(({F_0'}^*\boxtimes F_0)\otimes \pi_*(\FD))).
\end{array}
$$
It follows also from propositions~\ref{gskk} and~\ref{gskkp} that the first of the above differentials
is a morphism $d_1:\Delta^Y_*(\CR(-7H_Y)) \to \Delta^Y_*\FF_1$.
Let us show that this map is zero. For this we note that by proposition~\ref{gskkp}
the restriction of $d_1$ to the open subset $\TY \setminus \TZ = Y \setminus Z \subset Y$
takes form
$$
\CEnd(\CO_{\TY\setminus\TZ} \oplus \CK_{|\TY\setminus\TZ}) \to
\CK_{|\TY\setminus\TZ}^*\otimes (W/\CK)_{|\TY\setminus\TZ}\otimes\CEnd(\CO_{\TY\setminus\TZ} \oplus \CK_{|\TY\setminus\TZ})
$$
It is clear that this map is induced by the canonical section of the vector bundle
$\CK^*\boxtimes W/\CK$ on $\TY\times\TY$. But this section vanishes on the diagonal $\TY \subset \TY\times\TY$,
hence the above map is zero. Finally, since the sheaf $\FF_1$ on $\TY$ is torsion free, it follows
that the map $d_1:\Delta^Y_*(\CR(-7H_Y)) \to \Delta^Y_*\FF_1$ is zero on the whole~$\TY$,
hence $\BE_2^{-2,-3} = \BE_1^{-2,-3}$.
Similarly, the second differential of the spectral sequence
is a morphism $d_2:\Delta^Y_*(\CR(-7H_Y)) \to \Delta^Y_*\FF_0$. Its restriction to
$\TY \setminus \TZ = Y \setminus Z \subset Y$ takes form
$$
\CEnd(\CO_{\TY\setminus\TZ} \oplus \CK_{|\TY\setminus\TZ}) \to
\CO_{\TY\setminus\TZ}(H_Y - H_G) \otimes (\Lambda^2(W/\CK)/\CO_\TY(H_G - H_Y))_{|\TY\setminus\TZ}\otimes\CEnd(\CO_{\TY\setminus\TZ} \oplus \CK_{|\TY\setminus\TZ})
$$
In the same way as before, it is clear that this map is induced by the canonical section of the vector bundle
$\CO_\TY(H_Y - H_G) \boxtimes (\Lambda^2(W/\CK)/\CO_\TY(H_G - H_Y))$ on $\TY\times\TY$.
This section also vanishes on the diagonal $\TY \subset \TY\times\TY$,
hence the above map is zero, and since the sheaf $\FF_0$ on $\TY$ is torsion free, it follows
that the map $d_2:\Delta^Y_*(\CR(-7H_Y)) \to \Delta^Y_*\FF_0$ is zero on the whole~$\TY$,
hence $\BE_3^{-2,-3} = \BE_2^{-2,-3}$.
Thus we see that
$$
\TCH^{-5}(\CT\otimes\pi_*(\FD)) =
\BE_\infty^{-5} = \BE_3^{-2,-3} = \BE_2^{-2,-3} = \BE_1^{-2,-3} =
\TCH^{-3}(({F_2'}^*\boxtimes F_2)\otimes \pi_*(\FD)) =
\rho_\natural^*(\Delta^Y_*\CR(-7H_Y))
$$
which completes the proof.
\end{proof}

\section{The kernel}\label{s_ker}

Recall vector bundles $E_0,E_1,E_2$ on $X$ and $F_0,F_1,F_2$ on~$\TY$ defined in~\eqref{e012} and~\eqref{f012def}.
Note that the spaces $\Hom(E_2,E_1) \cong \Hom(E_1,E_0) \cong W^*$ are canonically dual
to $\Hom(F_2,F_1) \cong \Hom(F_1,F_0) \cong W$. This allows to construct canonical maps on $X\times \TY$
$$
E_k\boxtimes F_k \to W\otimes (E_{k-1}\boxtimes F_k) \to E_{k-1}\boxtimes F_{k-1},
$$
where the first map is the coevaluation $E_k \to \Hom(E_k,E_{k-1})^*\otimes E_{k-1}$ tensored by $F_k$,
and the second map is the evaluation $\Hom(F_k,F_{k-1})\otimes F_k \to F_{k-1}$ tensored by $E_{k-1}$.

\begin{lemma}
The composition of the morphisms $E_2\boxtimes F_2 \to E_1\boxtimes F_1 \to E_0\boxtimes F_0$ is zero.
\end{lemma}
\begin{proof}
The maps are given by the canonical elements $\id_W \in W^*\otimes W = \Hom(E_k,E_{k-1})\otimes\Hom(F_k,F_{k-1})$,
hence there composition is given by the image of $\id_W\otimes \id_W \in {W^*}^{\otimes 2}\otimes W^{\otimes 2}$
in $S^2W^*\otimes\Lambda^2W = \Hom(E_2,E_0)\otimes\Hom(F_2,F_0)$ which is easily seen to be zero.
\end{proof}

We consider the complex $\{E_2\boxtimes F_2 \to E_1\boxtimes F_1 \to E_0\boxtimes F_0\}$ as an object
of the derived category $\D^b(X\times\TY)$. Let $Q(X,\TY)$ denote the incidence quadric
$$
Q(X,\TY) = (X\times\TY)\times_{(\BP^\vee\times\BP)} Q,
$$
where $Q \subset \BP^\vee\times\BP$ is the usual incidence quadric.
Note that $Q(X,\TY)$ is a divisor on $X\times\TY$, and moreover
$$
Q(X,\TY) \in |H_X+H_Y|,
$$
where $H_X$ and $H_Y$ are the divisor classes of hyperplane sections of $X$ and $Y$ respectively.

Let $i:Q(X,\TY) \to X\times\TY$ denote the embedding.
Then we have the following resolution
\begin{equation}\label{qxhy}
0 \to \CO_{X\times\TY}(-H_X-H_Y) \to \CO_{X\times\TY} \to i_*\CO_{Q(X,\TY)} \to 0.
\end{equation}

\begin{lemma}\label{le}
There exists a coherent sheaf $\CE \in \Coh(Q(X,\TY))$ such that
\begin{equation}\label{equ_e}
0 \to E_2\boxtimes F_2 \to E_1\boxtimes F_1 \to E_0\boxtimes F_0 \to i_*\CE \to 0
\end{equation}
is an exact sequence. Moreover, $i_*\CE \in \TD_{\opp X}^0$.
\end{lemma}
\begin{proof}
Consider the locus of pairs of intersecting subspaces in $X\times \BG = \Gr(2,W)\times\Gr(n-4,W)$
and let $T \subset X\times\TY$ be its preimage. It is clear that $T$ is a subscheme in $Q(X,\TY)$
(if the kernel of a skew-form intersects with a given 2-dimensional subspace then this
subspace is isotropic for the skew-form). We are going to show that $\CE \cong \CJ_{T,Q(X,\TY)}(H_X+H_G)$,
the sheaf of ideals of $T$ in $Q(X,\TY)$ twisted by $\CO_{Q(X,\TY)}(H_X+H_G)$. For this we consider
the projectivization $\PP_X(\CU)$ and the following morphism of vector bundles on $\PP_X(\CU)\times\TY$:
$$
\CO_{\PP_X(\CU)/X}(-1)\boxtimes\CO_\TY \to \CU\boxtimes\CO_\TY \to W\otimes\CO_{\PP_X(\CU)}\boxtimes\CO_{\TY} \to \CO_{\PP_X(\CU)}\boxtimes(W/\CK)
$$
(the first two morphisms here are the natural embeddings and the third morphism is the natural projection).
The composition can be considered as a global section of the vector bundle $\CO_{\PP_X(\CU)/X}(1)\boxtimes(W/\CK)$ on $\PP_X(\CU)\times\TY$.
The zero locus of this section is a desingularization $\TT$ of $T$.
Consider the Koszul resolution of $\CO_\TT$ on $\PP_X(\CU)\times\TY$:
$$
\dots \to
\CO_{\PP_X(\CU)/X}(-2)\boxtimes \Lambda^2 \CK^\perp \to
\CO_{\PP_X(\CU)/X}(-1)\boxtimes \CK^\perp \to
\CO_{\PP_X(\CU)\times\TY} \to
\CO_\TT \to 0.
$$
Its pushforward to $X\times\TY$ gives the following resolution of $\CO_T$
$$
0 \to S^2\CU(-H_X)\boxtimes\Lambda^4\CK^\perp \to \CU(-H_X)\otimes\Lambda^3\CK^\perp \to \CO_X(-H_X)\boxtimes\Lambda^2\CK^\perp \to \CO_{X\times\TY} \to \CO_T \to 0.
$$
Twisting it by $\CO_X(H_X)\boxtimes\CO_\TY(H_G)$ and taking into account an isomorphism
$\Lambda^k\CK^\perp(H_G) \cong \Lambda^{4-k}(W/\CK)$ we obtain the following resolution of $\CO_T(H_X+H_G)$:
$$
0 \to S^2\CU\boxtimes\CO_\TY \to \CU\otimes(W/\CK) \to \CO_X\boxtimes\Lambda^2(W/\CK) \to \CO_X(H_X)\boxtimes\CO_\TY(H_G) \to \CO_T(H_X+H_G) \to 0.
$$
It remains to note that the composition of the canonical embedding $\CO_X\boxtimes\CO_\TY(H_G-H_Y) \to \CO_X\boxtimes\Lambda^2(W/\CK)$
with the map $\CO_X\boxtimes\Lambda^2(W/\CK) \to \CO_X(H_X)\boxtimes\CO_\TY(H_G)$ from this resolution
coincides with the embedding $\CO_X\boxtimes\CO_\TY(H_G-H_Y) \to \CO_X(H_X)\boxtimes\CO_\TY(H_G)$
obtained from~(\ref{qxhy}) by a $\CO_X(H_X)\boxtimes\CO_\TY(H_G)$-twisting, so that we obtain an exact triple of complexes
$$
\xymatrix@R=10pt{
&& \CO_X\boxtimes\CO_\TY(H_G-H_Y) \ar[r] \ar[d] & \CO_X(H_X)\boxtimes\CO_\TY(H_G) \ar@{=}[d] \\
S^2\CU\boxtimes\CO_\TY \ar[r] \ar@{=}[d] & \CU\otimes(W/\CK) \ar[r] \ar@{=}[d] & \CO_X\boxtimes\Lambda^2(W/\CK) \ar[r] \ar[d] & \CO_X(H_X)\boxtimes\CO_\TY(H_G) \\
E_2\boxtimes F_2 \ar[r] & E_1\boxtimes F_1 \ar[r] & E_0\boxtimes F_0
}
$$
The corresponding long exact sequence of the cohomology sheaves shows that we have a quasiisomorphism
$\{E_2\boxtimes F_2 \to E_1\boxtimes F_1 \to E_0\boxtimes F_0\} \cong i_*\CJ_{T,Q(X,\TY)}$.

To prove the second claim we first note that $i_*\CE \in \D^b(X)\boxtimes\TD_\opp = \TD_{\opp X}$.
Moreover, $i_*\CE \in \TD_{\opp X}^{\ge 0}$ since the functor
$\rho^{\opp X}_*$ is left-exact and $i_*\CE \in \D^{\ge 0}(X\times\TY)$.
On the other hand, $i_*\CE \in \TD_{\opp X}^{\le 0}$ by lemma~\ref{ffst}.
\end{proof}

\section{The proof}\label{uhps}

Recall the notation $\BP = \PP(\Lambda^2W^*)$.
Let $\CX \subset X\times\BP$ be the universal hyperplane section of $X$.

\begin{lemma}[\cite{K2}]\label{H_X1}
The universal hyperplane section $\CX$ is a smooth projective variety,
flat over $\BP$ and its relative dimension over $\BP$ equals $2n-5$.
\end{lemma}
\begin{proof}
It is easy to see that the projection $\CX \to X$ is smooth
(in fact it is a projectivization of a vector bundle), hence $\CX$ is smooth.
On the other hand, the fibers of the projection $\CX \to \BP$ are hyperplane sections
of $X$, hence all of them have dimension $\dim X - 1 = 2n - 5$. The flatness is evident.
\end{proof}

Recall that the incidence quadric $Q(X,\TY)$ can be identified
with the fiber product $\CX\times_{\BP}\TY$. Let
$j:Q(X,\TY) \to \CX\times\TY$ be the corresponding embedding. Then
the sheaf $j_*\CE \in \Coh(\CX\times\TY)$ gives a kernel functor
$\Phi_{j_*\CE}:\D^b(\TY) \to \D^b(\CX)$.
We are going to show that the functor $\Phi_{j_*\CE}:\D^b(\TY) \to \D^b(\CX)$
is fully faithful on the subcategory $\TD \subset \D^b(\TY)$.
This will take the most part of this section.
We will use appropriately modified arguments of \cite{K1}.

Let $\alpha$ denote the embedding $\CX \to X\times\BP$.
Note that we have the following resolution of $\CO_\CX$ on $X\times\BP$

\begin{equation}\label{ex1}
0 \to \CO_X(-H_X)\boxtimes\CO_\BP(-H_P) \to \CO_{X\times\BP} \to \alpha_*\CO_{\CX} \to 0.
\end{equation}

Let $\beta$ denote the embedding of $\TY$ to $\BP\times \TY$
given by the graph of $g$. Then we have a commutative square
\begin{equation}\label{abij}
\vcenter{\xymatrix{
Q(X,\TY) \ar[r]^i \ar[d]^j &
X\times \TY \ar[d]^{\beta} \\
\CX\times \TY \ar[r]^-{\alpha} &
X\times \BP \times \TY
}}
%\vcenter{\xymatrix{
%Q(X,Y) \ar[r]^i \ar[d]^j &
%X\times Y \ar[d]^{\beta} \\
%\CX\times Y \ar[r]^-{\alpha} &
%X\times \BP \times Y
%}}
\end{equation}
where we write $\alpha$ instead of $\alpha\times\id_Y$ and $\beta$
instead of $\id_X\times\beta$ for brevity.
Further, consider an object
$\CE^* = \RCHom_{Q(X,\TY)}(\CE,\CO_{Q(X,\TY)}) \in \D^b(Q(X,\TY))$.

\begin{lemma}\label{ces_sh}
We have $\CE^* \in \D^{[0,1]}(Q(X,\TY))$
and there is a quasiisomorphism on $X\times \TY$
\begin{equation}\label{equ_estar}
\{E_0^*(-H_X)\boxtimes F_0^*(-H_Y) \to
E_1^*(-H_X)\boxtimes F_1^*(-H_Y) \to
E_2^*(-H_X)\boxtimes F_2^*(-H_Y) \}
\cong
i_*\CE^*[1].
\end{equation}
Moreover, $i_*\CE^* \in \TD_{X}^{[0,1]}$.
\end{lemma}
\begin{proof}
Applying the functor
$\RCHom(-,\CO_X(-H_X)\boxtimes\CO_\TY(-H_Y))$ to~(\ref{equ_e}) we obtain a quasiisomorphism
\begin{multline*}
\{E_0^*(-H_X)\boxtimes F_0^*(-H_Y) \to
E_1^*(-H_X)\boxtimes F_1^*(-H_Y) \to
E_2^*(-H_X)\boxtimes F_2^*(-H_Y) \}
\cong \\ \cong
\RCHom(i_*\CE,\CO_X(-H_X)\boxtimes\CO_\TY(-H_Y))[2]
\end{multline*}
on $X\times \TY$.
On the other hand, we have
$\CO_{Q(X,\TY)} = i^*\CO_{X\times \TY} = i^!\CO_X(-H_X)\boxtimes\CO_\TY(-H_Y)[1]$,
since $i$ is a divisorial embedding, and by the duality theorem we have
\begin{multline*}
i_*\CE^* =
i_*\RCHom(\CE,\CO_{Q(X,Y)}) =
i_*\RCHom(\CE,i^!\CO_X(-H_X)\boxtimes\CO_\TY(-H_Y))[1] =
\\ =
\RCHom(i_*\CE,\CO_X(-H_X)\boxtimes\CO_\TY(-H_Y))[1].
\end{multline*}
Combining these two isomorphisms we deduce~(\ref{equ_estar}).
Now note that from~(\ref{equ_estar}) it follows that $i_*\CE^* \in \D^{[-1,1]}(X\times\TY)$.
On the other hand, $\CE^* \in \D^{\ge 0}(Q(X,\TY))$ since the dualization functor is left exact.
It remains to note that the functor $i_*$ is exact and conservative since $i$ is a closed embedding.

To prove the second claim we note that $i_*\CE^* \in \TD_X^{\ge 0}$ since the functor
$\rho^X_*$ is left-exact and $i_*\CE^* \in \D^{\ge 0}(X\times\TY)$.
On the other hand, $i_*\CE \in \TD_X^{\le 1}$ by lemma~\ref{ffst}.
\end{proof}

Denote
\begin{equation}\label{ce1d}
\CE_1 = j_*\CE,
\qquad
\CE^{\#t} = \CE^*((t-n+1)H_X+H_Y),
\qquad
%\quad\text{and}\quad
\CE_1^{\#t} = j_*\CE^{\#t}\otimes\CO_\BP(t)[2n-5].
\end{equation}

\begin{lemma}\label{eed}
We have $\CE_1 \in \TD_{\opp\CX}^0$, $\CE_1^{\#t} \in \TD_{\CX}^{[5-2n,6-2n]}$.
\end{lemma}
\begin{proof}
Use $\alpha_*j_* = \beta_*i_*$ and apply lemma~\ref{finitet}
together with lemmas~\ref{le} and~\ref{ces_sh}.
\end{proof}

\begin{lemma}\label{ed_la}
The functor $\Phi_{\CE_1^{\#0}}$ is left adjoint to $\Phi_{\CE_1}$.
\end{lemma}
\begin{proof}
By lemma~\ref{ladjoint} it suffices to check that
$\CE_1^{\#0} \cong \RCHom(\CE_1,\omega_{\CX}[\dim \CX])$.
Let $q_1$ and $q$ be the projections $\CX\times \TY\to \TY$ and
$X\times \TY \to  \TY$. It is clear that $q\circ i = q_1\circ j$.
Using the duality theorem and the functoriality of the twisted pullback
we deduce
\begin{multline*}
\RCHom(\CE_1,\omega_{\CX}[\dim \CX]) \cong
\RCHom(j_*\CE,q_1^!\CO_\TY) \cong
j_*\RCHom(\CE,j^!q_1^!\CO_\TY) \cong
j_*\RCHom(\CE,i^!q^!\CO_\TY) \cong \\ \cong
j_*\RCHom(\CE,i^!\omega_X[\dim X]) \cong
j_*\RCHom(\CE,\omega_X(H_X+H_Y)[\dim X-1]) \cong
j_*\CE^*(H_Y-(n-1)H_X)[2n-5],
\end{multline*}
since $\omega_X \cong \CO_X(-nH_X)$ and $\dim X = 2n-4$.
\end{proof}

Twisting (\ref{equ_estar}) by $\CO_{X\times\TY}((t-n+1)H_X+H_Y)$
and taking into account~(\ref{ce1d}),
we obtain a quasiisomorphism
\begin{equation}\label{equ_ediez}
\{E_0^*((t-n)H_X)\boxtimes F_0^* \to
E_1^*((t-n)H_X)\boxtimes F_1^* \to
E_2^*((t-n)H_X)\boxtimes F_2^* \}
\cong
i_*\CE^{\#t}[1].
\end{equation}

Consider the following diagram
\begin{equation}\label{diag_1}
\vcenter{\xymatrix{
\TY \times \CX \ar[d]^\alpha &
\TY \times \CX \times \TY \ar[d]^\alpha \ar[l]_{p_{12}} \ar[r]^{p_{23}} &
\CX \times \TY \ar[d]^\alpha \\
\TY \times (X \times \BP)  &
\TY \times (X \times \BP) \times \TY \ar[d]^q \ar[l]_{p_{12}} \ar[r]^{p_{23}} &
(X \times \BP) \times \TY \\
& \TY \times \BP \times \TY \ar[d]^\pi \\
& \TY \times \TY
}}
\end{equation}
where $q$ is the projection along $X$ and $\pi$ is the projection along $\BP$.
Consider the following objects in $\D(\TY\times \BP \times \TY)$:
\begin{equation}\label{cct}
\CC_t = q_*\alpha_*(p_{12}^*\CE_1^{\#t}\otimes p_{23}^*\CE_1),
\qquad\qquad
\TCC_t = q_*(\alpha_*p_{12}^*\CE_1^{\#t}\otimes \alpha_*p_{23}^*\CE_1)
\end{equation}

\begin{lemma}\label{k}
The convolution of kernels $\CE_1$ and $\CE_1^{\#t}$ is given by
$\CE_1\circ\CE_1^{\#t} \cong \pi_*\CC_t \in \D^b(\TY\times \TY)$.
\end{lemma}
\begin{proof}
Use the definition of the convolution (section~\ref{kf}) and note that
$\pi\circ q\circ \alpha = p_{13}:\TY\times\CX\times\TY \to \TY\times\TY$ is
the projection along $\CX$.
\end{proof}

\begin{lemma}\label{lbcc}
There exists an integer $N\in\ZZ$ such that
$\CC_t \in \TD_{\natural\BP}^{[-N,1]} \subset \D^b(\TY\times\BP\times \TY)$ for all $t\in\ZZ$.
\end{lemma}
\begin{proof}
By lemma~\ref{tensort} and lemma~\ref{eed} we have
$p_{12}^*\CE_1^{\#t} \otimes p_{23}^*\CE_1 \in
\TD_{\natural\CX}^{[-N,6-2n]} \subset \D^b(\TY\times\CX\times \TY)$
for some $N \in \ZZ$.
We conclude by lemma~\ref{mapt} and lemma~\ref{H_X1}.
\end{proof}

\begin{lemma}\label{cctt}
We have an exact triangle
\begin{equation}\label{tcccct}
\TCC_t \to
\CC_t \to
\CC_{t-1}[2]
\end{equation}
in $\TD_{\natural\BP} \subset \D^b(\TY\times\BP\times \TY)$.
\end{lemma}
\begin{proof}
Since $\alpha:\CX\to X\times\BP$ is a divisorial embedding,
and $\CX$ is a zero locus of a section of the bundle
$\CO_X(H_X)\boxtimes\CO_\BP(H_P)$, we have a distinguished triangle
$$
\alpha^*\alpha_*\CF \to \CF \to \CF\otimes(\CO_X(-H_X)\boxtimes\CO_\BP(-H_P))[2]
$$
for any object $\CF$ on $\TY\times\CX\times \TY$.
Taking $\CF = p_{23}^*\CE_1$, tensoring with $p_{12}^*\CE_1^{\#t}$,
applying $q_*\alpha_*$ and taking into account the projection formula
$\alpha_*(p_{12}^*\CE_1^{\#t}\otimes \alpha^*\alpha_*p_{23}^*\CE_1) \cong
\alpha_*p_{12}^*\CE_1^{\#t}\otimes \alpha_*p_{23}^*\CE_1$,
and the definition (\ref{cct}) of $\CC_t$ and $\TCC_t$,
we obtain (\ref{tcccct}).
\end{proof}

Recall the objects~$\FD \in \D^b(\TY\times\BP\times\TY)$ and $\CT,\CT^* \in \D^b(\TY\times\TY)$
defined in~\eqref{fd}, \eqref{tts1} and \eqref{tts2}.

\begin{lemma}\label{tcc}
If $n = 6$ then we have $\TCC_t \in \TD_{\natural\BP}^{\ge -11} \subset \D^b(\TY\times\BP\times \TY)$
for all $t\in\ZZ$. Moreover,
$$
\TCC_t = \begin{cases}
\pi^*\CT^*\otimes \FD,          & \text{for $t=0$}\\
F_2^*\boxtimes\CO_\BP(3H_P)\boxtimes F_2\otimes\FD[4],
                                & \text{for $t=3$}\\
\pi^*\CT\otimes \CO_\BP(6H_P)\otimes\FD[6],
                                & \text{for $t=6$}\\
0,                              & \text{for $t=1,2,4,5$}
\end{cases}
$$
In particular, $\TCC_0 \in \TD_{\natural\BP}^{\le 0}$, $\TCC_3 \in \TD_{\natural\BP}^{\ge -5}$,
and $\TCC_6 \in \TD_{\natural\BP}^{\ge -9}$.

If $n = 7$ then we have $\TCC_t \in \TD_{\natural\BP}^{\ge -15} \subset \D^b(\TY\times\BP\times \TY)$
for all $t\in\ZZ$. Moreover,
$$
\TCC_t = \begin{cases}
\pi^*\CT^*\otimes \FD,          & \text{for $t=0$}\\
\pi^*\CT\otimes \CO_\BP(7H_P)\otimes\FD[8],
                                & \text{for $t=7$}\\
0,                              & \text{for $t=1,2,3,4,5,6$}
\end{cases}
$$
In particular, $\TCC_0 \in \TD_{\natural\BP}^{\le 0}$ and $\TCC_7 \in \TD_{\natural\BP}^{\ge -13}$.
Moreover, $\TCC_8 \in \TD_{\natural\BP}^{\ge -13}$.
\end{lemma}
\begin{proof}
Consider the diagram
\begin{equation}\label{diag_2}
\vcenter{\xymatrix{
&
\TY \times \CX \ar[d]^\alpha &
\TY \times \CX \times \TY \ar[d]^\alpha \ar[l]_{p_{12}} \ar[r]^{p_{23}} &
\CX \times \TY \ar[d]^\alpha \\
Q(X,\TY) \ar[d]^i \ar[ur]^j &
\TY \times (X \times \BP)  &
\TY \times (X \times \BP) \times \TY \ar[l]_{p_{12}} \ar[r]^{p_{23}} \ar[d]^q &
(X \times \BP) \times \TY &
Q(X,\TY) \ar[d]^i \ar[ul]_j \\
\TY \times X \ar[ur]^{\beta'} &
\TY \times X \times \TY \ar[ur]^{\beta'} \ar[l]_{p_{12}}&
\TY\times\BP\times \TY \ar[d]^\pi &
\TY \times X \times \TY \ar[ul]_\beta    \ar[r]^{p_{23}}&
X \times \TY \ar[ul]_\beta \\
&& \TY\times \TY
}}
\end{equation}
and note that the maps $p_{12}$ and $p_{23}$ are flat.
Hence we have,
$$
\alpha_*p_{23}^*\CE_1 \cong
p_{23}^*\alpha_*\CE_1 =
p_{23}^*\alpha_*j_*\CE \cong
p_{23}^*\beta_*i_*\CE \cong
\beta_*p_{23}^*i_*\CE
$$
and similarly
$\alpha_*p_{12}^*\CE_1^{\#t} \cong
\beta'_*p_{12}^*i_*\CE^{\#t}\otimes\CO_\BP(tH_P)[2n-5]$.
Therefore,
$$
\TCC_t \cong
q_*(\beta'_*p_{12}^*i_*\CE^{\#t}\otimes\CO_\BP(tH_P)[2n-5] \otimes
\beta_*p_{23}^*i_*\CE).
$$
Applying $\beta_*p_{23}^*$ to~(\ref{equ_e}) and
$\beta'_*p_{12}^*$ to~(\ref{equ_ediez}) we obtain quasiisomorphisms
of $\beta_*p_{23}^*i_*\CE$ and $\beta'_*p_{12}^*i_*\CE^{\#t}[1]$
with the following complexes
%\begin{equation}\label{bpie}
$$
\arraycolsep=1pt
\begin{array}{lllllll}
\{ & \beta_*(\CO\boxtimes E_2\boxtimes F_2) & \to &
\beta_*(\CO\boxtimes E_1\boxtimes F_1) & \to &
\beta_*(\CO\boxtimes E_0\boxtimes F_0) & \}, \\[5pt]
\{ & \beta'_*(F_0^*\boxtimes E_0^*((t-n)H_X)\boxtimes \CO) & \to &
\beta'_*(F_1^*\boxtimes E_1^*((t-n)H_X)\boxtimes \CO) & \to &
\beta'_*(F_2^*\boxtimes E_2^*((t-n)H_X)\boxtimes \CO) & \}.
\end{array}
$$
%\end{equation}
So, $\TCC_t$ can be represented by the complex
\begin{equation}\label{tccres}
\left\lbrace \TCC_t^{0,2} \to
\TCC_t^{0,1} \oplus \TCC_t^{1,2} \to
\TCC_t^{0,0} \oplus \TCC_t^{1,1} \oplus \TCC_t^{2,2} \to
\TCC_t^{1,0} \oplus \TCC_t^{2,1} \to
\TCC_t^{2,0} \right\rbrace ,
\end{equation}
with the rightmost term placed in degree $1$, where
\begin{multline*}
\TCC_t^{k,l} :=
q_*(\beta'_*(F_k^*\boxtimes E_k^*((t-n)H_X)\boxtimes \CO)
\otimes\CO_\BP(tH_P)[2n-5] \otimes
\beta_*(\CO\boxtimes E_l\boxtimes F_l)) \cong  \\ \cong
q_*((F_k^* \boxtimes (E_k^*(t-n)\otimes E_l)
\boxtimes \CO_\BP(tH_P) \boxtimes F_l) \otimes
\beta'_*\CO_{\TY\times X\times \TY} \otimes
\beta_*\CO_{\TY\times X\times \TY})[2n-5]
\end{multline*}
for $k,l = 0,1,2$.
%(more rigorously, we should express each of the quasiisomorphisms for
%$\beta_*p_{23}^*i_*\CE$ and $\beta'_*p_{12}^*i_*\CE^{\#t}[1]$ as a pair of distinguished
%triangles, then tensoring these triangles we would obtain quite a lot of new triangles
%expressing $\TCC_t$ in terms of $\TCC_t^{k,l}$).
To compute $\TCC_t^{k,l}$ we consider the diagram
$$
\xymatrix{
\TY \times X \times \TY \ar[r]^-{\beta'} \ar[d]_q &
\TY \times (X \times \BP) \times \TY \ar[d]_q &
\TY \times X \times \TY \ar[l]_-{\beta} \ar[d]_q \\
\TY \times \TY \ar[r]^-{\beta'} &
\TY \times \BP \times \TY &
\TY \times \TY \ar[l]_-{\beta}
}
$$
and note that
$$
\beta'_*\CO_{\TY\times X\times \TY} \otimes \beta_*\CO_{\TY\times X\times \TY} \cong
\beta'_*q^*\CO_{\TY\times \TY} \otimes \beta_*q^*\CO_{\TY\times \TY} \cong
q^*\beta'_*\CO_{\TY\times \TY} \otimes q^*\beta_*\CO_{\TY\times \TY} \cong
%q^*(\beta'_*\CO_{\TY\times \TY} \otimes \beta_*\CO_{\TY\times \TY}) \cong
q^*\FD.
$$
since $q$ is flat.
Substituting this into the formula for $\TCC_t^{k,l}$ we deduce
\begin{equation}\label{tcctkl}
\TCC^{k,l}_t \cong
(F_k^* \boxtimes
(H^\bullet(X,E_k^*((t-n)H_X)\otimes E_l) \otimes \CO_\BP(tH_P))
\boxtimes F_l) \otimes \FD)[2n-5]
\end{equation}
By lemma~\ref{fkflfd} we have
$(F_k^*\boxtimes\CO_\BP(tH_P)\boxtimes F_l)\otimes\FD \in \TD_{\natural\BP}^{\ge -1}$ for $n=6$ and
$(F_k^*\boxtimes\CO_\BP(tH_P)\boxtimes F_l)\otimes\FD \in \TD_{\natural\BP}^{\ge -3}$ for $n=7$.
On the other hand, it is clear that
$H^\bullet(X,E_k^*((t-n)H_X)\otimes E_l) \in \D^{\ge 0}(\kk)$. Therefore
$$
%\TCC_t^{k,l} \in \TD_{\natural\BP}^{[-8,1]},  \quad \text{ if $n=6$} \qquad \text{and} \qquad
%\TCC_t^{k,l} \in \TD_{\natural\BP}^{[-12,1]}, \quad \text{ if $n=7$}
\begin{array}{ll}
\TCC_t^{k,l} \in \TD_{\natural\BP}^{\ge -8}, & \text{if $n=6$}\\[2pt]
\TCC_t^{k,l} \in \TD_{\natural\BP}^{\ge -12}, & \text{if $n=7$}
\end{array}
$$
for all $t\in\ZZ$, $k,l = 0,1,2$.
Hence, looking at~(\ref{tccres}) we see that
$\TCC_t \in \D_{\natural\BP}^{\ge -11}$ if $n = 6$ and
$\TCC_t \in \D_{\natural\BP}^{\ge -15}$ if $n = 7$
for all $t\in\ZZ$, which gives us the first claims of the lemma.

Now we can compute the cohomology groups $H^\bullet(X,E_k^*((t-n)H_X)\otimes E_l)$
for $0\le t\le n$ explicitly via the Borel--Bott--Weil thoerem on $X = \Gr(2,W)$.
We have
$$
H^\bullet(X,E_k^*((t-n)H_X)\otimes E_l) \cong
\Hom^\bullet(E_k,E_l((t-n)H_X)) \cong
\begin{cases}
S^{k-l}W^*, & \text{for $t=n$ and $k\ge l$} \\
S^{l-k}W[4-2n], & \text{for $t=0$ and $l\ge k$} \\
\kk[-4], & \text{for $t=3$, $n=6$, $k=l=2$} \\
0, & \text{otherwise if $0\le t\le n$}
\end{cases}
$$
Substituting this into~(\ref{tcctkl}) and using (\ref{tccres}) we deduce that
$$
\arraycolsep = 1pt
\left\{
\begin{array}{ccccc}
&&&& F_0^*\boxtimes\CO_\BP\boxtimes F_0 \\
&& W\otimes F_0^*\boxtimes\CO_\BP\boxtimes F_1 && \oplus \\
S^2W\otimes F_0^*\boxtimes\CO_\BP\boxtimes F_2 & \to & \oplus & \to & F_1^*\boxtimes\CO_\BP\boxtimes F_1 \\
&& W\otimes F_1^*\boxtimes\CO_\BP\boxtimes F_2 && \oplus \\
&&&& F_2^*\boxtimes\CO_\BP\boxtimes F_2
\end{array}
\right\}\otimes\FD \cong \TCC_0
%\{
%S^2W\otimes F_0^*\boxtimes\CO_\BP\boxtimes F_2\otimes\FD \to
%W\otimes F_0^*\boxtimes\CO_\BP\boxtimes F_1\otimes\FD \oplus
%W\otimes F_1^*\boxtimes\CO_\BP\boxtimes F_2\otimes\FD \to
%F_0^*\boxtimes\CO_\BP\boxtimes F_0\otimes\FD \oplus
%F_1^*\boxtimes\CO_\BP\boxtimes F_1\otimes\FD \oplus
%F_2^*\boxtimes\CO_\BP\boxtimes F_2\otimes\FD
%\} \cong \TCC_0
$$
$$
F_2^*\boxtimes\CO_\BP(3H_P)\boxtimes F_2\otimes\FD[4] \cong \TCC_3,
\qquad\text{if $n=6$}
$$
$$
\arraycolsep = 1pt
\left\{
\begin{array}{ccccc}
F_0^*\boxtimes\CO_\BP(nH_P)\boxtimes F_0 \\
\oplus && W^*\otimes F_1^*\boxtimes\CO_\BP(nH_P)\boxtimes F_0 \\
F_1^*\boxtimes\CO_\BP(nH_P)\boxtimes F_1 & \to & \oplus & \to & S^2W^*\otimes F_0^*\boxtimes\CO_\BP(nH_P)\boxtimes F_2\\
\oplus && W^*\otimes F_2^*\boxtimes\CO_\BP(nH_P)\boxtimes F_1 \\
F_2^*\boxtimes\CO_\BP(nH_P)\boxtimes F_2
\end{array}
\right\}\otimes\FD \cong \TCC_n[6-2n]
$$
and
$$
\begin{array}{ll}
\TCC_1 = \TCC_2 = \TCC_4 = \TCC_5 = 0, & \text{for $n=6$}\\[2pt]
\TCC_1 = \TCC_2 = \TCC_3 = \TCC_4 = \TCC_5 = \TCC_6 = 0, & \text{for $n=7$}
\end{array}
$$
The formulas and the cohomological bounds for $\TCC_t$ with $0\le t\le n$ evidently follow
from lemma~\ref{tts} and lemma~\ref{fkflfd}.

So, it remains to check that $\TCC_8 \in \TD_{\natural\BP}^{\ge -13}$ if $n = 7$.
For this we need to know $\TCC_8^{k,l}$ for $k\le l$.
Computing $H^\bullet(X,E_k^*(H_X)\otimes E_l) \cong \Hom^\bullet(E_k,E_l(H_X))$
via the Borel--Bott--Weil Theorem, we see that $\TCC_8$ takes form
$$
\arraycolsep = 3pt
\left\{
\begin{array}{ccccc}
&& \hphantom{W^*}\Lambda^2W^*\otimes F_0^*\boxtimes\CO_\BP(8H_P)\boxtimes F_0\\
W^*\otimes F_0^*\boxtimes\CO_\BP(8H_P)\boxtimes F_1 && \oplus \\
\oplus & \to & W^*\otimes W^*\otimes F_1^*\boxtimes\CO_\BP(8H_P)\boxtimes F_1 & \to & \dots  \\
W^*\otimes F_1^*\boxtimes\CO_\BP(8H_P)\boxtimes F_2 && \oplus \\
&& W^*\otimes W^*\otimes F_2^*\boxtimes\CO_\BP(8H_P)\boxtimes F_2
\end{array}
\right\}\otimes\FD \
$$
In particular, we see that $\TCC_8^{0,2} = 0$ which already implies $\TCC_8 \in \TD_{\natural\BP}^{\ge -14}$.
Finally, we see that the maps
$\{F_1^* \to W^*\otimes F_2^*\} = \{\CK^\perp \to W^*\otimes\CO_\TY\}$ and
$\{F_0^* \to W^*\otimes F_1^*\} = \{\Lambda^2\CK^\perp \to \CO_\TY(H_Y - H_G) \oplus W^*\otimes\CK^\perp\}$
are injective, hence $\TCC_8 \in \TD_{\natural\BP}^{\ge -13}$.
\end{proof}

\begin{lemma}\label{cc}
We have
$\CC_t \in \TD_{\natural\BP}^{[-10,1]}(\TY\times\BP\times \TY)$ if $n=6$, and
$\CC_t \in \TD_{\natural\BP}^{[-14,1]}(\TY\times\BP\times \TY)$ if $n=7$
for all $t$.
%
%\begin{equation}\label{lbcc1}
%\begin{array}{ll}
%\CC_t \in \TD_{\natural\BP}^{[-10,1]}(\TY\times\BP\times \TY), & \text{if $n=6$}\\[2pt]
%\CC_t \in \TD_{\natural\BP}^{[-14,1]}(\TY\times\BP\times \TY), & \text{if $n=7$}
%\end{array}
%\end{equation}
\end{lemma}
\begin{proof}
We know already that $\CC_t \in \TD_{\natural\BP}^{\le 1}$ by lemma~\ref{lbcc}.
So, it remains to establish the left cohomological bound.
It follows from triangle~(\ref{tcccct}) and lemma~\ref{tcc} that
for any $l \le -11$ in case $n=6$ and for any $l \le -15$ in case $n=7$
we have an exact sequence
$$
0 = \TCH^{l-2}(\TCC_t) \to
\TCH^{l-2}(\CC_t) \to
\TCH^{l-2}(\CC_{t-1}[2]) \to
\TCH^{l-1}(\TCC_t) = 0.
$$
Therefore,
$\TCH^l(\CC_{t-1}) = \TCH^{l-2}(\CC_{t-1}[2]) \cong \TCH^{l-2}(\CC_t)$.
So, in case $n=6$ if $\TCH^l(\CC_t) \ne 0$ for $l \le -11$ then
$\TCH^{l-2s}(\CC_{t+s}) \ne 0$ for all $s$ which
contradicts the claim of lemma~\ref{lbcc},
and similarly in case $n=7$.
\end{proof}

\begin{lemma}\label{cc3}
If $n=6$ then $\CC_3[-6] \in \TD_{\natural\BP}^{\ge 2}$.
\end{lemma}
\begin{proof}
Triangles~\eqref{tcccct} for $t = 5$ and $t = 4$ together with lemma~\ref{tcc}
give isomorphisms $\CC_5 \cong \CC_4[2] \cong \CC_3[4]$.
Then triangle~\eqref{tcccct} for $t = 6$ shifted by $[-12]$ takes form
$$
\TCC_6[-12] \to \CC_6[-12] \to \CC_3[-6].
$$
Now, $\TCC_6[-12] \in \TD_{\natural\BP}^{\ge 3}$ by lemma~\ref{tcc},
while $\CC_6[-12] \in \TD_{\natural\BP}^{\ge 2}$ by lemma~\ref{cc}
and the claim follows.
\end{proof}

\begin{lemma}\label{cc0}
If $n=6$ then $\CC_0 \in \TD_{\natural\BP}^{[0,1]}$ and $\TCH^0(\CC_0) \cong \TCH^{-5}(\TCC_3)$.
\end{lemma}
\begin{proof}
Triangles~\eqref{tcccct} for $t = 2$ and $t = 1$ together with lemma~\ref{tcc}
give isomorphisms $\CC_2 \cong \CC_1[2] \cong \CC_0[4]$.
Then triangle~\eqref{tcccct} for $t = 3$ shifted by $[-6]$ takes form
$$
\TCC_3[-6] \to \CC_3[-6] \to \CC_0.
$$
Now, $\TCC_3[-6] \in \TD_{\natural\BP}^{\ge 1}$ by lemma~\ref{tcc},
while $\CC_3[-6] \in \TD_{\natural\BP}^{\ge 2}$ by lemma~\ref{cc3}.
Therefore $\CC_0 \in \TD_{\natural\BP}^{\ge 0}$ and
$\TCH^0(\CC_0) = \TCH^1(\TCC_3[-6]) = \TCH^{-5}(\TCC_3)$.
\end{proof}

\begin{lemma}\label{cc7}
If $n=7$ then $\CC_7 \in \TD_{\natural\BP}^{\ge -12}$.
\end{lemma}
\begin{proof}
Triangle~~\eqref{tcccct} for $t = 8$ takes form
$$
\TCC_8 \to \CC_8 \to \CC_7[2].
$$
Note that $\TCC_8 \in \TD_{\natural\BP}^{\ge -13}$ by lemma~\ref{tcc},
and $\CC_8 \in \TD_{\natural\BP}^{\ge -14}$ by lemma~\ref{cc}.
Therefore $\CC_7[2] \in \TD_{\natural\BP}^{\ge -14}$
and the lemma follows.
\end{proof}

\begin{lemma}\label{cc07}
If $n=7$ then $\CC_0 \in \TD_{\natural\BP}^{[0,1]}$ and $\TCH^0(\CC_0) \cong \TCH^{-13}(\TCC_7)$.
\end{lemma}
\begin{proof}
Triangles~\eqref{tcccct} for $t = 1,\dots,6$ together with lemma~\ref{tcc}
give isomorphisms $\CC_6 \cong \CC_5[2] \cong \dots \cong \CC_0[12]$.
Then triangle~\eqref{tcccct} for $t = 7$ shifted by $[-14]$ takes form
$$
\TCC_7[-14] \to \CC_7[-14] \to \CC_0.
$$
Now, $\TCC_7[-14] \in \TD_{\natural\BP}^{\ge 1}$ by lemma~\ref{tcc},
while $\CC_7[-14] \in \TD_{\natural\BP}^{\ge 2}$ by lemma~\ref{cc7}.
Therefore $\CC_0 \in \TD_{\natural\BP}^{\ge 0}$ and
$\TCH^0(\CC_0) = \TCH^1(\TCC_7[-14]) = \TCH^{-13}(\TCC_7)$.
\end{proof}

\begin{proposition}\label{c0}
$\CC_0 \in \TD_{\natural\BP}^0$ and
$\pi_*\CC_0 \cong \rho_\natural^*({\Delta_Y}_*\CR) \in \TD_{\natural}^0$.
\end{proposition}
\begin{proof}
Triangle~\eqref{tcccct} for $t = 0$ takes form
$$
\TCC_0 \to \CC_0 \to \CC_{-1}[2].
$$
Now, $\TCC_0 \in \TD_{\natural\BP}^{\le 0}$ by lemma~\ref{tcc},
while $\CC_{-1}[2] \in \TD_{\natural\BP}^{\le -1}$ by lemma~\ref{cc}.
Therefore $\CC_0 \in \TD_{\natural\BP}^{\le 0}$.
Combining with lemma~\ref{cc0} we conclude that
$\CC_0 \in \TD_{\natural\BP}^0$ and
$\CC_0 \cong \TCH^0(\CC_0) \cong \TCH^{-5}(\TCC_3)$ if $n=6$ and
combining with lemma~\ref{cc07} we conclude that
$\CC_0 \cong \TCH^0(\CC_0) \cong \TCH^{-13}(\TCC_7)$ if $n=7$.
Using again lemma~\ref{tcc} we see that
$$
\begin{array}{ll}
\CC_0 \cong \TCH^{-1}((F_2^*\boxtimes\CO_\BP(3H_P)\boxtimes F_2)\otimes\FD), & \text{if $n=6$,}\\[2pt]
\CC_0 \cong \TCH^{-5}(\pi^*\CT\otimes\CO_\BP(7H_P)\otimes\FD), & \text{if $n=7$.}
\end{array}
$$
Since $\FD$ is supported on $\TY\times_\BP\BP\times\BP\TY = \TY\times_\BP\TY$
and $\TY\times_\BP\times\TY$ is finite over $\TY\times\TY$
it follows from lemma~\ref{finitet} that in case $n = 6$ we have
\begin{multline*}
\pi_*\CC_0 \cong
\pi_*\TCH^{-1}((F_2^*\boxtimes\CO_\BP(3H_P)\boxtimes F_2)\otimes\FD) \cong
\TCH^{-1}(\pi_*((F_2^*\boxtimes\CO_\BP(3H_P)\boxtimes F_2)\otimes\FD)) \cong
\\ \cong
\TCH^{-1}((F_2^*\boxtimes\CO_\BP(3H_P)\boxtimes F_2)\otimes\pi_*\FD) \cong
\rho_\natural^*({\Delta_Y}_*\CR)
\end{multline*}
by lemma~\ref{h0fd}. Similarly, in case $n = 7$ we use lemma~\ref{tts} and deduce the claim.
\end{proof}

\begin{corollary}
The functor
$\Phi_{\CE_1^{\#0}}\circ\Phi_{\CE_1}:\D^b(\TY) \to \D^b(\TY)$ is isomorphic to
the functor $\rho^*\circ\rho_*$.
\end{corollary}
\begin{proof}
By lemma~\ref{k} and proposition~\ref{c0} the functor $\Phi_{\CE_1^{\#0}}\circ\Phi_{\CE_1}$
is given by the kernel $\rho_\natural^*({\Delta_Y}_*\CR)$. On the other hand,
by lemma~\ref{rhophirho} the same kernel gives the functor $\rho^*\circ\rho_*$.
\end{proof}

\begin{corollary}\label{tdff}
The restriction of the functor $\Phi_{\CE_1}:\D^b(\TY) \to \D^b(\CX)$
to the subcategory $\TD \subset \D^b(\TY)$ is fully faithful.
\end{corollary}
\begin{proof}
We have
$\Hom(\Phi_{\CE_1}(F),\Phi_{\CE_1}(F')) =
\Hom(\Phi_{\CE_1^{\#0}}(\Phi_{\CE_1}(F)),F')$.
But if $F,F'\in\TD$ then
$\rho^*(\rho_*(F)) \cong F$ by theorem~\ref{NCR} hence
$\Hom(\rho^*(\rho_*(F)),F') = \Hom(F,F')$.
\end{proof}

Now we are almost ready to give a proof of the main result of the paper, theorem~\ref{themain}.
The final preparatory step is the following

\begin{proposition}\label{ldcomp}
For every $0\le k\le 2$ the functor $\Phi_{\CE_1}^*\circ\pi^*:\D^b(X) \to \D^b(\TY)$
is full and faithful on the subcategory
$\langle E_2,E_1,E_0\rangle = \CA_0 \subset \D^b(X)$
and takes it to the subcategory $\langle F_0^*,F_1^*,F_2^*\rangle = \CB_0 \subset \TD \subset \D^b(\TY)$.
Moreover, it takes the bundle $\Lambda^2(W/\CU) \in \CA_0$ to (a shift of) $F_2^* = \CO_\TY$.
\end{proposition}
\begin{proof}
The functor $\Phi_{\CE_1}^*\circ\pi^*$ is given by the kernel
$\pi_*\CE_1^{\# 0} \cong
p_{1*}\alpha_*j_*\CE^{\# 0}[7] \cong
p_{1*}\beta_*i_*\CE^{\# 0}[7] \cong
i_*\CE^{\# 0}[7]$.
Using the resolution~\eqref{equ_ediez} and taking into account that
$H^\bullet(X,E_k^*(-n)\otimes E_l) = \Hom^\bullet(E_k,E_l(-n)) \cong \Hom^\bullet(E_l,E_k[2n-4])^*$
we deduce that up to a shift it takes
$$
E_0 \mapsto {}'F_0 := F_0^*,
\quad
E_1 \mapsto {}'F_1 := \{ W\otimes F_0^* \to F_1^* \},
\quad
E_2 \mapsto {}'F_2 := \{ S^2W\otimes F_0^* \to W\otimes F_1^* \to F_2^* \}.
$$
Since $\Hom(E_k,E_l) = S^{k-l}W^* = \Hom({}'F_k,{}'F_l)$
and the above functor evidently induces an isomorphism on these $\Hom$-spaces,
it gives an equivalence $\CA_0 \cong \CB_0$.
Finally, note that
$$
\Lambda^2(W/\CU) = \{ E_2 \to W\otimes E_1 \to \Lambda^2W\otimes E_0 \}
\qquad \mapsto \qquad
\{ {}'F_2 \to W\otimes {}'F_1 \to \Lambda^2W\otimes {}'F_0 \}
$$
which up to a shift coincides with $F_2^*$.
\end{proof}

\noindent{\bf Proof of theorem~\ref{themain}.}
\nopagebreak

\noindent
We are going to use theorem~\ref{hpd_crit}.
We already have proved in corollary~\ref{tdff} that the functor $\Phi_{\CE_1} = \Phi_{j_*\CE}$
embeds the category $\TD = \D^b(Y,\CR)$ fully and faithfully into $\D^b(\CX)$, the derived category
of the universal hyperplane section of $X$.
Let us check that the image of the functor $\Phi_{j_*\CE}$ is contained in the subcategory
$\CC = [\langle \CA_1(1)\boxtimes\D^b(\BP),\dots,
\CA_{n-1}(n-1)\boxtimes\D^b(\BP)\rangle]^\perp$.

Indeed, for this it suffices to check that $\Phi_{j_*\CE}^*(\CA_t(t)\boxtimes\D^b(\BP)) = 0$
for $1\le t\le n-1$. But the category $\CA_t(t)\boxtimes\D^b(\BP)$ is generated by objects
of the form $\alpha^*(E \boxtimes G)$ with $E \in \CA_t(t)$, $G \in \D^b(\BP)$, and
it is clear that
$$
\Phi_{j_*\CE}^*(\alpha^*(E \boxtimes G)) \cong
\Phi_{j_*\CE}^*(\pi^*E) \otimes g^*G
$$
since the functor $\Phi_{j_*\CE}^*$ is $\BP$-linear.
Thus it suffices to check that $\Phi_{j_*\CE}^*(\pi^*E) = 0$ for $E \in \CA_t(t)$ and $1\le t\le n-1$.
Arguing like in proposition~\ref{ldcomp} and taking into account that
$$
H^\bullet(X,E_k^*(-n)\otimes E) = \Hom(E_k,E(-n)) = \Hom(E,E_k[2n-4])^* = 0
$$
since $E \in \CA_t(t)$, $E_k \in \CA_0$, and~\eqref{ldx} is a semiorthogonal decomposition,
we conclude that $\Phi_{j_*\CE}^*(\pi^*E) = 0$.

Further, we note that by proposition~\ref{ldcomp} we have
$\Phi_{j_*\CE}^*(\CB'_t) = \CB_t$, hence due to $\BP$-linearity
of the functor $\Phi_{j_*\CE}^*$ we have $\Phi_{j_*\CE}^*(\CB'_t(-tH_P)) = \CB_t(-tH_Y)$.
Finally, the subcategories $\CB_t(-tH_Y)$ form a dual Lefschetz collection in $\TD = \D^b(Y,\CR)$
by proposition~\ref{ldcb}. Therefore by theorem~\ref{hpd_crit} the noncommutative variety
$(Y,\CR)$ is Homologically Projectively Dual to $X$ and its dual Lefschetz decomposition
is given by~\eqref{ldtd}.
%
%Moreover, the kernel $\CE_1$
%is a pushforward of the object $\CE \in (Y,\CR^\opp)\times_{\BP} \CX$
%(actually $\CE \in \D^b(\TY\times_\BP\CX) \cap \TD_{\opp X} = \D^b(Q(X,\TY)) \cap \TD_{\opp X}$
%by lemma~\ref{le} and the latter category naturally identifies with
%$\D^b((Y,\CR^\opp)\times_{\BP} \CX)$).
%
%It remains to check that $\D^b(\CX) = \langle \Phi_{\CE_1}(\TD),\CA_1(1)\boxtimes\D^b(\BP),\dots,\CA_5(5)\boxtimes\D^b(\BP)\rangle$.
%For this we note that actually it gets as an outcome of section~6 of~\cite{K2} if we only check that
%$\Phi_{\CE_1}^*$ takes the dual Lefschetz collection $\langle\CB_{11}(-11),\dots,\CB_1(-1),\CB_0\rangle$
%of the subcategory $[\langle \CA_1(1)\boxtimes\D^b(\BP),\dots,\CA_5(5)\boxtimes\D^b(\BP)\rangle]^\perp$
%of $\D^b(\CX)$, which was constructed in~\cite{K2}, section~5, to a dual Lefschetz collection in $\TD$.
%But this is true by proposition~\ref{ldcomp} and proposition~\ref{ldcb}.

\section{Applications to linear sections of $\Gr(2,6)$}\label{s_gr26}

\nc{\thecase}[1]{\medskip\noindent{$\mathbf r = #1.$}}

In this section $W$ is a vector space, $\dim W = 6$,
$X = \Gr(2,W) = \Gr(2,6)$,
$Y = \Pf(4,W^*)$ is a cubic hypersurface in $\PP(\Lambda^2W^*)$,
$Z = \Pf(2,W^*) = \Gr(2,W^*) \subset \PP(\Lambda^2W^*)$,
and $\CR$ is a sheaf of $\CO_Y$-algebras on~$Y$ defined in section~\ref{ncrpf}.
As a consequence of Homological Projective Duality between $X$ and $(Y,\CR)$ we have the following

\begin{corollary}
Let $L \subset \Lambda^2W^*$ be a vector subspace, $\dim L = r$,
and put $L^\perp \subset \Lambda^2W$ for its orthogonal.
Assume that $X_L = X\cap \PP(L^\perp)$, $Y_L = Y\cap\PP(L)$, and $Z_L = Z\cap \PP(L)$
have expected dimension, i.e.
\begin{equation}\label{adm}
\dim X_L = 8 - r, \quad
\dim Y_L = r - 2, \quad
\dim Z_L = r - 7.
\end{equation}
Then we have the following semiorthogonal decompositions
\begin{equation}\label{xlyl}
\begin{array}{lll}
\D^b(X_L) &=& \langle \CC_L,\CA_{r}(1),\dots,\CA_{5}(6-r)\rangle\\
\D^b(Y_L,\CR) &=& \langle \CB_{11}(3-r),\dots,\CB_{15-r}(-1),\CC_L\rangle.
\end{array}
\end{equation}
\end{corollary}

Let us write down explicitly semiorthogonal decompositions~\eqref{xlyl} for different values of $r$.
Note that for $r\le 3$ we have $\CC_L = \D^b(Y_L,\CR)$ and for $r\ge 6$ we have $\CC_L = \D^b(X_L)$.
Note also that $Z_L$ is empty for $r \le 6$ by~\eqref{adm}, hence the algebra $\CR$ on $Y_L$
is a matrix algebra, so that $\D^b(Y_L,\CR) \cong \D^b(Y_L)$.

\thecase1
In this case conditions~\eqref{adm} mean that $Y_L$ is empty, hence $L$ is spanned by a nondegenerate
skew-form $\omega \in \Lambda^2W^*$, so $X_L$ is a Lagrangian Grassmannian $\LGr(2,W)$ with respect
to $\omega$. Decompositions~\eqref{xlyl} then give
$$
\D^b(X_L) = \langle \CA_1(1),\dots,\CA_5(5) \rangle,
$$
reproving theorem~3.1 of \cite{K4}.

\thecase2
In this case conditions~\eqref{adm} mean that $Y_L = \{y_1,y_2,y_3\}$ is a scheme of length $3$, corresponding
to a pencil of skew-forms $\omega$ in $\Lambda^2W^*$ with degenerations of rank $4$.
In this case $\CC_L = \D^b(Y_L)$, and the first decomposition of~\eqref{xlyl} gives
$$
\D^b(X_L) = \langle \D^b(Y_L),\CA_2(1),\dots,\CA_5(4) \rangle.
$$
In particular, we deduce the following

\begin{corollary}
A smooth linear section $X_L$ of $X = \Gr(2,W) = \Gr(2,6)$ of codimension $2$
admits an exceptional collection of length $12$. Explicitly,
$$
\D^b(X_L) = \langle \Phi_{\CE}(y_1),\Phi_{\CE}(y_2),\Phi_{\CE}(y_3),\CA_2(1),\dots,\CA_5(4) \rangle.
$$
\end{corollary}

\thecase3
In this case conditions~\eqref{adm} mean that $Y_L$ is an elliptic curve, corresponding
to a net of skew-forms $\omega$ in $\Lambda^2W^*$ with degenerations of rank $4$.
In this case $\CC_L = \D^b(Y_L)$, and the first decomposition of~\eqref{xlyl} gives
$$
\D^b(X_L) = \langle \D^b(Y_L),\CA_3(1),\CA_4(2),\CA_5(3) \rangle.
$$

\thecase4
In this case conditions~\eqref{adm} mean that $Y_L$ is a del Pezzo surface of degree $3$.
In this case decompositions~\eqref{xlyl} give
$$
\D^b(X_L) = \langle \CC_L,\CA_4(1),\CA_5(2) \rangle,
\qquad
\D^b(Y_L) = \langle \CO_{Y_L}(-1),\CC_L \rangle
$$
Since $\D^b(Y_L)$ admits an exceptional collection
(of length $9$) when $Y_L$ is smooth, we deduce the following

\begin{corollary}
A smooth linear section $X_L$ of $X = \Gr(2,W) = \Gr(2,6)$ of codimension $4$
admits an exceptional collection of length $12$.
\end{corollary}

\thecase5
In this case conditions~\eqref{adm} mean that $Y_L$ is a cubic $3$-fold,
and $X_L$ is a Fano threefold $V_{14}$. In this case $\CC_L = \D^b(Y_L)$,
and decompositions~\eqref{xlyl} give
$$
\D^b(X_L) = \langle \CC_L,\CU^*,\CO(1) \rangle,
\qquad
\D^b(Y_L) = \langle \CO_{Y_L}(-2),\CO_{Y_L}(-1),\CC_L \rangle,
$$
reproving theorem~3.1 of \cite{K0}.

\thecase6
In this case conditions~\eqref{adm} mean that $Y_L$ is a Pfaffian cubic $4$-fold,
and $X_L$ is a K3-surface of degree $14$. In this case $\CC_L = \D^b(X_L)$,
and the second decomposition of~\eqref{xlyl} gives the following

\begin{theorem}
Let $Y_L$ be a Pfaffian cubic $4$-fold and let $X_L$ be the orthognal K{\rm3}-surface of degree~$14$.
Then we have a semiorthogonal decomposition
$$
\D^b(Y_L) = \langle \CO_{Y_L}(-3),\CO_{Y_L}(-2),\CO_{Y_L}(-1),\D^b(X_L) \rangle.
$$
\end{theorem}

\thecase7
In this case conditions~\eqref{adm} mean that $X_L$ is a curve of genus $8$,
$Y_L$ is a Pfaffian cubic $5$-fold (which is necessarily singular), and $Z_L$
is a scheme of length $14$.
Certainly, $Z_L$ is contained in the singular locus of $Y_L$, but actually,
the latter can be strictly bigger then $Z_L$.
In this case $\CC_L = \D^b(X_L)$,
and $\D^b(Y_L,\CR)$ is a noncommutative (partial) resolution of $Y_L$ (at $Z_L$).
The second decomposition of~\eqref{xlyl} gives the following
$$
\D^b(Y_L,\CR) = \langle \CB_{11}(-4),\dots,\CB_8(-1),\D^b(X_L) \rangle.
$$

\thecase{8,\dots,14}
In these cases conditions~\eqref{adm} mean that $Y_L$ is a Pfaffian cubic of dimension
$(r-2)$. In this case $\CC_L = \D^b(X_L)$, which is either zero (for $r\ge 9$),
or admits an exceptional collection (for $r = 8$ and $X_L$ smooth),
and $\D^b(Y_L,\CR)$ is a noncommutative (partial) resolution of $Y_L$ (at $Z_L$).
The second decomposition of~\eqref{xlyl} gives the following
$$
\D^b(Y_L,\CR) = \langle \CB_{11}(3-r),\dots,\CB_{15-r}(-1),\D^b(X_L) \rangle.
$$
In particular, if $X_L$ is smooth we deduce that $\D^b(Y_L,\CR)$
admits a full exceptional collection.

\section{Applications to linear sections of $\Gr(2,7)$}\label{s_gr27}

In this section $W$ is a vector space, $\dim W = 7$,
$X = \Gr(2,W) = \Gr(2,7)$,
$Y = \Pf(4,W^*)$ has codimension 3 in $\PP(\Lambda^2W^*)$,
$Z = \Pf(2,W^*) = \Gr(2,W^*) \subset \PP(\Lambda^2W^*)$,
and $\CR$ is a sheaf of $\CO_Y$-algebras on~$Y$ defined in section~\ref{ncrpf}.
As a consequence of Homological Projective Duality between $X$ and $(Y,\CR)$ we have the following.

\begin{corollary}
Let $n = 7$.
Let $L \subset \Lambda^2W^*$ be a vector subspace, $\dim L = r$,
and put $L^\perp \subset \Lambda^2W$ for its orthogonal.
Assume that $X_L = X\cap \PP(L^\perp)$, $Y_L = Y\cap\PP(L)$, and $Z_L = Z\cap \PP(L)$
have expected dimension, i.e.
\begin{equation}\label{adm1}
\dim X_L = 10 - r, \quad
\dim Y_L = r - 4, \quad
\dim Z_L = r - 11.
\end{equation}
Then we have the following semiorthogonal decompositions
\begin{equation}\label{xlyl1}
\begin{array}{lll}
\D^b(X_L) &=& \langle \CC_L,\CA_{r}(1),\dots,\CA_{6}(7-r)\rangle\\
\D^b(Y_L,\CR) &=& \langle \CB_{13}(7-r),\dots,\CB_{21-r}(-1),\CC_L\rangle.
\end{array}
\end{equation}
\end{corollary}

Let us write down explicitly semiorthogonal decompositions~\eqref{xlyl1} for different values of $r$.
Note that for $r\le 7$ we have $\CC_L = \D^b(Y_L,\CR)$ and for $r\ge 7$ we have $\CC_L = \D^b(X_L)$.
Note also that $Z_L$ is empty for $r \le 10$ by~\eqref{adm1}, hence the algebra $\CR$ on $Y_L$
is a matrix algebra, so that $\D^b(Y_L,\CR) \cong \D^b(Y_L)$.

\thecase{1,2,3}
In this case conditions~\eqref{adm1} mean that $Y_L$ is empty.
It follows that $X_L$ is a smooth codimension~$r$ linear section of $X = \Gr(2,7)$.
Decompositions~\eqref{xlyl1} then give
$$
\D^b(X_L) = \langle \CA_r(1),\dots,\CA_6(7-r) \rangle.
$$
In particular, we deduce the following

\begin{corollary}
A smooth linear section $X_L$ of $X = \Gr(2,W) = \Gr(2,7)$ of codimension $\le 3$
admits an exceptional collection of length $21 - 3r$ consisting of bundles
$S^2\CU(t),\CU(t),\CO(t)$ with $1\le t\le 7-r$.
\end{corollary}

\thecase4
In this case conditions~\eqref{adm1} mean that $Y_L = \{y_1,y_2,\dots,y_{42}\}$ is a scheme of length $\deg Y = 42$,
and the first decomposition of~\eqref{xlyl1} gives
$$
\D^b(X_L) = \langle \D^b(Y_L),\CA_4(1),\CA_5(2),\CA_6(3) \rangle.
$$
In particular, we deduce the following

\begin{corollary}
A smooth linear section $X_L$ of $X = \Gr(2,W) = \Gr(2,7)$ of codimension $4$
admits an exceptional collection of length $51$.
Explicitly,
$$
\D^b(X_L) = \langle\Phi_\CE(y_1),\dots,\Phi_\CE(y_{42}),S^2\CU(1),\CU(1),\CO(1),S^2\CU(2),\CU(2),\CO(2),S^2\CU(3),\CU(3),\CO(3)\rangle.
$$
\end{corollary}

\thecase5
In this case conditions~\eqref{adm1} mean that $X_L$ is a Fano 5-fold of index 2,
$Y_L$ is a curve of degree $\deg Y = 42$ in a half-canonical embedding (so that $g(Y_L) = 43$),
and the first decomposition of~\eqref{xlyl1} gives
$$
\D^b(X_L) = \langle \D^b(Y_L),\CA_5(1),\CA_6(2) \rangle.
$$

\thecase6
In this case conditions~\eqref{adm1} mean that $X_L$ is a Fano 4-fold of index 1,
$Y_L$ is a canonically embedded surface of degree $\deg Y = 42$,
and the first decomposition of~\eqref{xlyl1} gives
$$
\D^b(X_L) = \langle \D^b(Y_L),\CA_6(1) \rangle,
$$

\thecase7
This case is the most interesting.
Conditions~\eqref{adm1} mean that both $X_L$ and $Y_L$ are Calabi-Yau 3-folds,
and decompositions~\eqref{xlyl1} boil down to just an equivalence of $\D^b(X_L)$
and $\D^b(Y_L)$. Note that we don't need $X_L$ and $Y_L$ to be smooth.
Let us formulate this result explicitly.

\begin{theorem}
Assume that $\dim W = 7$, $\dim L = 7$ and let $L \subset \Lambda^2 W^*$.
If\/ $\PP(L)$ doesn't intersect $\Gr(2,W^*) \subset \PP(\Lambda^2W^*)$
and the corresponding linear sections of the Grassmannian
$X_L = \Gr(2,W) \cap \PP(L^\perp)$ and of the Pfaffian variety
$Y_L = Y \cap \PP(L)$ are $3$-dimensional, then we have an equivalence
of categories
$$
\D^b(X_L) \cong \D^b(Y_L).
$$
\end{theorem}

As it was explained in the Introduction, this theorem is a generalization of the result of~\cite{BC}.

\thecase8
In this case conditions~\eqref{adm1} mean that $X_L$ is a canonically embedded surface
of degree $\deg X = 14$, $Y_L$ is a Fano 4-fold of index $1$,
and the second decomposition of~\eqref{xlyl1} gives
$$
\D^b(Y_L) = \langle \CB_{13}(-1),\D^b(X_L) \rangle.
$$

\thecase9
In this case conditions~\eqref{adm1} mean that $X_L$ is a
curve of degree $\deg X = 14$ in a half-canonical embedding
(so that $g(X_L) = 15$), $Y_L$ is a Fano 5-fold of index 2,
and the second decomposition of~\eqref{xlyl1} gives
$$
\D^b(X_L) = \langle \CB_{13}(-2),\CB_{12}(-1),\D^b(X_L) \rangle.
$$

\thecase{10}
In this case conditions~\eqref{adm1} mean that $X_L$ is a scheme of length $\deg X = 14$,
$Y_L$ is a Fano 6-fold of index $3$, and the second decomposition of~\eqref{xlyl1} gives
$$
\D^b(Y_L) = \langle \CB_{13}(-3),\CB_{12}(-2),\CB_{11}(-1),\D^b(X_L) \rangle.
$$
In particular, we deduce the following

\begin{corollary}
A smooth $6$-fold linear section $Y_L$ of the Pfaffian variety $Y$
admits an exceptional collection of length $23$ consisting of sheaves
$\Phi^*_{\CE}(x_i)$, $i = 1,\dots, 14$, where $X_L = \{x_i\}$,
and of the bundles $F_0^*(t),F_1^*(t),F_2^*(t)$ with $-3\le t\le -1$.
\end{corollary}

\thecase{11,\dots,20}
In these cases conditions~\eqref{adm1} mean that $X_L$ is empty, and
$Y_L$ is a Fano variety of dimension $r-4$ and of index $r-7$,
and the second decomposition of~\eqref{xlyl1} shows that the objects
and of the bundles $F_0^*(t),F_1^*(t),F_2^*(t)$ with $7-r\le t\le -1$
form an exceptional collection in the derived category
of a noncommutative crepant resolution $(Y_L,\CR)$ of $Y_L$.

%\end{document}

\end{document}